\documentclass{article}
\usepackage[utf8]{inputenc}
\usepackage[numbers]{natbib}
\usepackage{manfnt} 
\newcommand{\dnote}[1]{%
    \noindent 
    \begin{tabular}{@{}m{0.13\textwidth}@{}m{0.87\textwidth}@{}}%
        \huge\textdbend &#1%
    \end{tabular}%
    \par 
}
\usepackage{amsthm}

\newtheorem{theorem}{Theorem}
\newtheorem{lemma}{lemma}
\newtheorem{proposition}{Proposition}
\newtheorem{Cor}{Corollary}
\newtheorem{definition}{Definition}
\usepackage[top=0.8in, bottom=1in, left=1.25in, right=1.25in]{geometry}
\usepackage{scrextend}
\usepackage[utf8]{inputenc} 
\usepackage[T1]{fontenc}    
\usepackage[english]{babel}
\usepackage{lmodern}
\usepackage[bitstream-charter]{mathdesign}
\usepackage[table]{xcolor}
 \usepackage{hyperref}
 \definecolor{burgundy}{rgb}{0.5, 0.0, 0.13}
\definecolor{camel}{rgb}{0.76, 0.6, 0.42}
\definecolor{chamoisee}{rgb}{0.63, 0.47, 0.35}
\definecolor{grey1}{RGB}{128,128,128}
 \hypersetup{colorlinks=true,linkcolor=burgundy,citecolor=chamoisee,urlcolor=burgundy,linktoc=page}
\usepackage{url}            
\usepackage{booktabs}       
\usepackage{nicefrac}       
\usepackage{microtype}      
\usepackage{dsfont}         
\renewcommand{\epsilon}{\varepsilon}
\renewcommand{\phi}{\varphi}

\title{Three Rates of Convergence or Separation via U-Statistics in a Dependent Framework}

\author{%
  Quentin Duchemin\thanks{This work was supported by grants from Région Ile-de-France.} \\
 LAMA, Univ Gustave Eiffel, CNRS, Marne-la-Vallée, France.\\
  \texttt{quentin.duchemin@univ-eiffel.fr} \\
  $\And$\\
  Yohann De Castro\\
  Institut Camille Jordan, École Centrale de Lyon, Lyon, France \\
  \texttt{yohann.de-castro@ec-lyon.fr} \\
    $\And$\\
  Claire Lacour\\
  LAMA, Univ Gustave Eiffel, CNRS, Marne-la-Vallée, France. \\
  \texttt{claire.lacour@univ-eiffel.fr} 
}

\date{}

\usepackage{amsmath}
\usepackage{mathtools}
\DeclarePairedDelimiter\ceil{\lceil}{\rceil}
\DeclarePairedDelimiter\floor{\lfloor}{\rfloor}
\usepackage{diagbox}
\usepackage{subcaption}

\usepackage{algorithm,algorithmic,refcount}
\usepackage{wrapfig}
\usepackage{tikz}
\usepackage{makecell}

\usepackage{tikz}
\usetikzlibrary{%
  arrows,%
  calc,%
  shapes.geometric,%
  shapes.misc,%
  shapes.symbols,%
  shapes.arrows,%
  automata,%
  through,%
  positioning,%
  scopes,%
  decorations.shapes,%
  decorations.text,%
  decorations.pathmorphing,%
  shadows}
\usetikzlibrary{positioning}
\usepackage{multirow}
\usepackage{makecell}
\usepackage{mdframed}
\newtheorem{assumption}{Assumption}
\usepackage{savesym}
\savesymbol{checkmark}
\usepackage{dingbat}
\usepackage{wasysym}
\makeatletter
\g@addto@macro\appendix{%
  \addtocontents{toc}{\protect\setcounter{tocdepth}{1}}%
}
\makeatother

\usepackage{enumitem}
\newcolumntype{C}[1]{>{\centering\arraybackslash}m{#1}}

\begin{document}

\maketitle

\begin{abstract}
Despite the ubiquity of U-statistics in modern Probability and Statistics, their non-asymptotic analysis in a dependent framework may have been overlooked. In a recent work, a new concentration inequality for U-statistics of order two for uniformly ergodic discrete time Markov chains has been proved. In this paper, we put this theoretical breakthrough into action by pushing further the current state of knowledge in three different active fields of research. 
First, we establish a new exponential inequality for the estimation of spectra of integral operators with MCMC methods. The novelty is that this result holds for kernels with positive and negative eigenvalues, which is new as far as we know. 

In addition, we investigate generalization performance of online algorithms working with pairwise loss functions and Markov chain samples. We provide an online-to-batch conversion result by showing how we can extract a low risk hypothesis from the sequence of hypotheses generated by any online~learner.

We finally give a non-asymptotic analysis of a goodness-of-fit test on the density of the stationary measure of a Markov chain. We identify some classes of alternatives over which our test based on the~$L^2$ distance has a prescribed power.
\end{abstract}

\section{Introduction}
\label{sec:intro}

For the last twenty years, the phenomenon of the concentration of measure has received much attention. The main interesting
feature of concentration inequalities is that, unlike central limit theorems or large deviations inequalities, they are nonasymptotic. Among others, Pascal Massart, Michel Ledoux and Gabor Lugosi produced a series of works that led to a large span of powerful inequalities.  Their results have found application in model selection \citep[cf.][]{massart07,Lerasle16}, statistical learning \citep[cf.][]{clemencon17}, online learning \citep[cf.][]{online-pairwise} or random graphs \citep[cf.][]{CL18,DCD20}. Most of the concentration inequalities are formulated for U-statistics of order~$m$ \citep[cf.][Chapter 12]{van2000asymptotic}, which are defined as a sum of the form \[\sum_{1\leq i_1<\dots<i_m\leq n} h_{i_1, \dots,i_m}(X_{i_1}, \dots,X_{i_m}),\] where~$X_1, \dots,X_n$ are random variables taking values in a measurable space~$(E,\Sigma)$ (with~$E$ Polish) and where~$h_{i_1,\dots,i_m}$ are measurable functions of~$m$ variables~$h_{i_1,\dots,i_m}:E^{m} \rightarrow \mathds R$. The pioneering works considered independent random variables $(X_i)_{i\geq 1}$, an assumption that can be prohibitive for practical applications which often involve some dependence structure. To cope with this issue, some researchers left the independent setting by working with Markov chains or by adopting some mixing conditions and we refer for example to \cite{FJS18,jiang2018bernsteins,paulin15,adamczak:ustats,clemencon17}. The previous mentioned papers considered U-statistics of order $m=1$ and the non-asymptotic behaviour of tails of U-statistics of order $m\geq2$ in a dependent framework remains so far understudied. Recently, the two papers \cite{duchemin} and \cite{shen20} made a first step to fill this gap. While \cite{shen20} consider U-statistics of arbitrary order with smooth and symmetric kernels and work with mixing conditions, \cite{duchemin} are focused on U-statistics of order two for uniformly ergodic discrete time Markov chains and bounded kernels. 
Let us highlight that we work with the result of the former paper rather than the one from~\cite{shen20} since we need a concentration inequality valid for any initial distribution of the chain. We give further details at the beginning of Section~\ref{sec:concentration}.

Our paper is in the line of work of \cite{massart2000some} where concentration of measure is applied to tackle problems arising from model selection. In this work, we shed light on the large number of potential theoretical breakthroughs allowed by a better understanding of the non-asymptotic tail behavior of U-statistics of order two in a dependent framework. We present new theoretical results in three different branches of Statistics ranging from online learning to goodness-of-fit tests.  In Section~\ref{sec:contributions} we describe in details our three contributions, highlighting their applications to learning theory and the proof innovations.

\subsection{Our contributions}
\label{sec:contributions}

Our new results - that we referred to as applications
for brevity - push further the current state of knowledge in three different active areas of research in Probability, Statistics and Machine Learning. Although the recent progress in concentration inequality for U-statistics with dependent random variables is a key element in our proofs, our contributions are not a direct consequence of it. The purpose of this section is threefold: $(i)$ we present concisely our main results, $(ii)$ we highlight the proof innovations of our approach compared to previous works and $(iii)$ we propose relevant connections between our work and other important topics in learning theory and Statistics.

\begin{itemize}
    \item {\bf Estimation of spectra of signed integral operator with MCMC algorithms}
    (Section~\ref{sec:integral-ope})\\
        \noindent We study the convergence of sequence of spectra of kernel matrices towards the spectrum of some integral operator. Previous important works may include~\cite{adamczak2015} and, as far as we know, they all assume that the kernel is of positive-type (i.e. giving an integral operator with non-negative eigenvalues). Getting counterpart of those results for signed integral operators is of great interest since they arise for example in random graphs with latent space~\citep[cf.][]{CL18} which can be characterized by the so-called graphon~\citep[cf.][]{lovasz2012large}. For the first time, this paper proves a non-asymptotic result of convergence of spectra for kernels that are not of positive-type. We further prove that {\it independent Hastings algorithms} are valid sampling schemes to apply our result. \\
        {\bf Proof innovations.} In Section~\ref{sec:appli1-comparison}, we propose a comparison between our result and the one from~\cite{adamczak2015}. We explain why working with integral operators of positive-type allows Adamczak and Berdnorz to make use of a powerful decoupling technique. Thanks to this elegant argument, they are reduced to prove a concentration inequality for a sum of Banach space valued random variables where the~$i$-th summand depends only on the~$i$-th visited state of the Markov chain. By considering signed integral operators, the approach of the former paper cannot be adapted. Our proof relies on a low rank approximation of the kernel and on a concentration result for U-statistics with dependent random variables.
        \\
        {\bf Application to learning theory.}  A large number of learning algorithms aim at estimating the eigenvalues and/or the eigenvectors of data-dependent matrices. This is for example the case for Principal Component Analysis (PCA) or some manifold methods~\citep[cf.][]{rosasco2010learning}. It appears that these matrices can often be interpreted as the empirical versions of continuous objects such as integral operators. As highlighted in~\cite{rosasco2010learning}, the theoretical analysis of the above mentioned learning algorithms requires to quantify the difference between the eigen-structure of the empirical operators and their continuous counterparts. Specific examples coming from the Machine Learning and the Statistics communities where our result may find an echo include the estimation of the entire spectrum of a Markov operator~\citep[cf.][]{chakraborty2019}, estimation procedures in random graphs~\citep[cf.][]{DCD20} or the analysis of the generalization properties of neural networks~\citep[cf.][]{zhang2021understanding}.
    \item {\bf Online learning with pairwise loss functions} (Section~\ref{sec:online-pairwise})\\
        \noindent
In Machine Learning, several important problems involve a pairwise loss function, i.e. a loss
function which depends on a pair of examples.
One typical example is the problem of metric learning~\citep[cf.][]{jin2009} where one aims to learn a metric so that instances with the same labels are close while ones with different labels are far away from each other. Other pairwise learning tasks include preference learning~\citep[cf.][]{xing2002distance}, ranking~\citep[cf.][]{agarwal2012}, gradient
learning~\citep[cf.][]{meir2003} and AUC maximization~\citep[cf.][]{zhao2011online}. Batch learning algorithms with pairwise loss functions have been extensively studied and their generalization
properties have been well established. However, batch algorithms have some limitations especially when data becomes available in a sequential order or for large scale learning problems where their computational cost can be prohibitive. Online algorithms have been designed to efficiently solve learning problems in such situations: they deal with data coming on fly and try to improve the learned model along time based on the new observations. 
The performance of online learning algorithms is typically analyzed through the notion of {\it regret} which compares the payoff obtained by the algorithm along time with the one that would have been obtained by taking the optimal decision at each time step~\citep[cf.][]{bubeck}. The regret quantifies the number of mistakes made by the algorithm without requiring assumptions on the way the training sequence is generated. When the sequence of observations is the realization of some stochastic process, one can analyze online algorithms through a different lens by wondering how they generalize on future examples. More precisely, we would like to convert a regret bound of an online learner into a control of the excess risk. In the online learning research community, these types of results are called {\it online-to-batch conversion} and we refer to~\cite[Section 3.7]{onlinesurvey} for a comprehensive introduction to this topic. Online-to-batch conversion results for online learning with univariate or pairwise loss functions working with i.i.d. samples have been considered for quite a while in both Machine Learning and Statistics literature~\citep[cf.][]{online-pairwise,ying2017,guo2017,chen2018}. For dependent data sequences, generalization bounds for online algorithms have also been proved in the last decades with univariate loss functions~\citep[cf.][]{agarwal2012}. However, theoretical guarantees for the generalization performance of online algorithms with pairwise loss functions with non i.i.d. data have been so far little studied. Inspired by~\cite{online-pairwise}, our work is one of the first to bring results regarding this problem. In Section~\ref{sec:appli2-literature}, we establish clear connections with the existing literature.\\
        {\bf Proof innovations.} 
        \cite{online-pairwise} was a pioneering work for the study of generalization performance of online learning algorithms with pairwise loss functions and worked with i.i.d. observations. In this paper, we extend the result of \cite{online-pairwise} by considering a dependent framework that makes the theoretical analysis more challenging. In our proofs, we bypass the additional issues arising from data dependency using properties of uniformly ergodic Markov chains, concentration inequalities for U-statistics (of order one and two) of dependent random variables and reversibility of Markov chains by considering the time-reversed sequence. Using the marker {\tiny $\textdbend$}, we shed light in Section B on the specific parts of the proof where the arguments used in the i.i.d. framework fail, requiring a specific theoretical work handling a sequence of dependent observations.
    \item {\bf Adaptive goodness-of-fit tests in a density model} (Section~\ref{sec:adaptivegof})\\
    \noindent
    Several works have already proposed goodness-of-fit tests for the density of the stationary distribution of a sequence of dependent random variables. In \cite{Li}, a test based on an $L^2$-type distance between the nonparametrically estimated conditional density and its model-based parametric counterpart is proposed. In \cite{BAI} a Kolmogorov-type test is considered. \cite{chwialkowski2016kernel} derive a test procedure for $\tau$-mixing sequences using Stein discrepancy computed in a reproducing kernel Hilbert space. In all the above mentioned papers, asymptotic properties of the test statistic are derived but no non-asymptotic analysis of the methods is conducted. As far as we know, this paper is the first to provide a non-asymptotic condition on the classes of alternatives ensuring that the statistical test reaches a prescribed power working in a dependent framework.
\end{itemize}

\clearpage

\subsection{Outline}

In Section~\ref{sec:notations-concentration}, we introduce useful notations for our paper and we present the concentration inequality for U-statistics that is an important argument of our proofs. The next three sections are dedicated to our main results.
We start by providing a convergence result for the estimation of spectra of integral operators with MCMC algorithms (see Section~\ref{sec:integral-ope}). We show that independent Hastings algorithms satisfy under mild conditions the assumptions of Section~\ref{sec:concentration} and we illustrate our result with the estimation of the spectra of some Mercer kernels. For the second application of our concentration inequality, we investigate the generalization performance of online algorithms with pairwise loss functions in a Markovian framework (see Section~\ref{sec:online-pairwise}). We motivate the study of such problems and we provide an online-to-batch conversion result. In a third and final application, we propose a goodness-of-fit test for the density of the stationary measure of a Markov chain (see Section~\ref{sec:adaptivegof}). We give an explicit condition on the set of alternatives to ensure that the statistical test proposed reaches a prescribed power. 
The proofs related to the three applications are given in Section~\ref{sec:proof-appli1}, Section~\ref{sec:proof-appli2} and Sections~\ref{proof-hypo-test-appli3}-\ref{proof-cor:appli3} respectively.

\section{Notations and Concentration inequality for U-statistics with dependent random variables}
\label{sec:notations-concentration}

\subsection{Notations}
Let us consider an arbitrary measurable space $(F,\mathcal F)$. 
For any measure~$\omega$ on~$(F,\mathcal F)$, the total variation norm of~$\omega$ is defined by $\|\omega\|_{\mathrm{TV}}:= \sup_{A \in \mathcal F}|\omega(A)|$. The space of square summable functions on $F$ with respect to the measure $\omega$ defined by \[L^2(\omega):=\{ f:F\to \mathds R \text{ measurable } | \int_F f(x)^2d\omega(x)<\infty\},\]endowed with the inner product \[(f,g) \in L^2(\omega) \times L^2(\omega) \mapsto \langle f,g\rangle:=\int_F f(x)g(x)d\omega(x),\] is a Hilbert space and we denote by $\|\cdot\|_2$ the norm induced by $\langle \cdot,\cdot\rangle$. For any function $h :F \to \mathds R$, we define the supremum norm of $h$ by~$\|h\|_{\infty}:=\sup_{x\in F} |h(x)|$. We denote by $\mathcal B(\mathds R)$ the Borel algebra on $\mathds R$ and we set $\mathds N^*:=\mathds N\backslash \{0\}$. For any $x \in \mathds R_+$, we denote by~$\floor{x}$ (resp. $\ceil{x}$) the largest integer that is less than or equal to $x$ (resp. the smallest integer greater than or equal to $x$). For any $x,y \in \mathds R$, we set~$x\vee y:=\max(x,y)$ and $x\wedge y:=\min(x,y).$ Given a sequence of real valued random variables~$(X_n)_{n\in \mathds N}$ and a sequence of positive reals~$(a_n)_{n\in\mathds N}$, the notation~$X_n=\mathcal O_{\mathds P}(a_n)$ means that $(X_n/a_n)_{n\in \mathds N}$ converges to zero in probability as~$n\to \infty$.

\subsection{Concentration inequality for U-statistics of uniformly ergodic Markov chains}
\label{sec:concentration}

In this section, we present the concentration result from \cite{duchemin} for U-statistics of uniformly ergodic discrete time Markov chains that will be an essential tool in our proofs. Let us mention that we do not work with the concentration inequality from~\cite{shen20} since it only holds for stationary chains if the kernel $h$ is $\pi$-canonical (see Assumption~\ref{assumption3}). Stationarity may be seen as a strong assumption which would make our main results from Section~\ref{sec:integral-ope} of little interest since MCMC methods are used when we are not able to directly sample from the distribution $\pi$. Regarding Sections~\ref{sec:online-pairwise} and~\ref{sec:adaptivegof}, the concentration inequality for U-statistics used in the proofs of our results needs to hold for any initial distribution of the chain.
\bigskip

Let $(E,\Sigma)$ be a measurable space. We consider a Markov chain $(X_i)_{i\geq 1}$ on $(E,\Sigma)$ with transition kernel $P:E\times E\to [0,1]$ and with a unique stationary distribution $\pi$. We consider a measurable function $h:(E \times E,\Sigma \otimes \Sigma) \to (\mathds R,\mathcal B(\mathds R))$ and we are interested in the following U-statistic
\[U_{\mathrm{stat}}(n):= \sum_{1\leq i \neq j \leq n} \left(h(X_i,X_j) - \mathds E_{(X,Y) \sim \pi \otimes \pi}[h(X,Y)]\right).\] We will work under the following set of assumptions.

\begin{assumption}
\label{assumption1} The Markov chain~$(X_i)_{i\geq1}$ is~$\psi$-irreducible~\citep[cf.][Section 4.2]{tweedie} for some maximal irreducibility measure~$\psi$ on~$\Sigma$. Moreover, there exist some natural number~$m$ and a constant~$\delta_m>0$ such that \begin{equation}\forall x \in E, \; \forall A \in \Sigma, \quad \delta_m \mu(A) \leq P^m(x,A) .\label{uni-ergodicity-small}\end{equation}
 for some probability measure~$\mu$.
 \end{assumption}
 
\pagebreak[3]

A Markov chain satisfying Assumption~\ref{assumption1} is called uniformly ergodic~\citep[cf.][Chapter 16]{tweedie} and admits a unique stationary distribution denoted by~$\pi$. Assumption~\ref{assumption1} also implies that the regeneration times associated to the split chain are exponentially integrable, meaning that their Orlicz norm with respect to the function $\psi_1(x) =\exp(x)-1$ are bounded by some constant $\tau>0.$ We refer to \cite[Section 2.3]{duchemin} for details.
\medskip

Assumption~\ref{assumption2} can be read as a reverse Doeblin's condition and is used in \cite{duchemin} as a decoupling tool. In their paper, the authors give several natural examples for which this condition holds.  
\begin{assumption}
 \label{assumption2} There exist~$\delta_M>0$ and some probability measure~$\nu$ such that \[\forall x \in E,\; \forall A \in \Sigma, \quad  P(x,A) \leq \delta_M  \nu(A).\] 
\end{assumption}

The last assumption introduces the notion of {\it $\pi$-canonical} kernel, which is the counterpart in the Markovian setting of the canonical (or degenerate) property of the independent framework.
\begin{assumption} \label{assumption3} Denoting by $\pi$ the stationary distribution of the Markov chain $(X_i)_{i\geq1}$, we assume that~$h:(E\times E,\Sigma \otimes \Sigma) \to (\mathds R,\mathcal B(\mathds R))$ is measurable, bounded and is~$\pi$-canonical, namely \[\forall x,y \in E, \quad \mathds E_{\pi}[h(X,x)]=\mathds E_{\pi}[h(X,y)]=\mathds E_{\pi}[h(x,X)]=\mathds E_{\pi}[h(y,X)].\]
 This common expectation will be denoted~$\mathds E_{\pi}[h]$.
\end{assumption}
\noindent Let us mention that several important kernels are~$\pi$-canonical. This is the case of translation-invariant kernels which have been widely studied in the Machine Learning community (cf. \cite{Lerasle16}). Another example of~$\pi$-canonical kernel is a rotation invariant kernel when~$E=\mathds S^{d-1}:=\{x \in \mathds R^d : \|x\|_2=1\}$ with~$\pi$ also rotation invariant \citep[cf.][]{CL18,DCD20}. Note also that if the kernel $h$ is not~$\pi$-canonical, the U-statistic decomposes into a linear term and a~$\pi$-canonical U-statistic. This is called the \textit{Hoeffding decomposition} \citep[cf.][p.176]{ginenick} and takes the following form 
\begin{align*}
&\sum_{i\neq j}\left( h(X_i,X_j)- \mathds{E}_{(X,Y) \sim \pi \otimes\pi}[h(X,Y)]\right)=\sum_{i\neq j}\left(\widetilde h(X_i,X_j) - \mathds E_{ \pi } \left[\widetilde h(X,\cdot) \right]\right) \\
&\qquad\qquad+  \sum_{i\neq j} \left(\mathds E_{X \sim \pi} \left[ h(X,X_j) \right] - \mathds E_{(X,Y) \sim \pi \otimes\pi} \left[ h(X,Y) \right]\right)\\
&\qquad\qquad+  \sum_{i\neq j} \left(\mathds E_{X \sim \pi} \left[ h(X_i,X) \right] - \mathds E_{(X,Y) \sim \pi \otimes\pi} \left[ h(X,Y) \right]\right),
\end{align*}
where the kernel~$\widetilde h$ is~$\pi$-canonical with 
\[\forall x,y \in E, \quad\widetilde h(x,y) = h(x,y)- \mathds E_{X \sim \pi} \left[ h(x,X) \right]-\mathds E_{X \sim \pi} \left[ h(X,y) \right].\]
We will use this method several times in our proofs (for example in Eq.\eqref{hoeffding-decomposition}).

\pagebreak[3]

We are now ready to state the result from \cite{duchemin} that is one key theoretical tool to derive our three contributions presented in the next section.
\begin{theorem}
\label{mainthm2} Suppose that Assumptions~\ref{assumption1},~\ref{assumption2} and~\ref{assumption3} are satisfied. Then there exist constants~$\beta,\kappa>0$ (depending on the Markov chain $(X_i)_{i\geq1}$) such that for any~$u \geq 1$ and any~$n\geq2$, with probability at least~$1-\beta e^{-u}\log n$,
\[\frac{2}{n(n-1)} U_{\mathrm{stat}}(n) \leq \kappa \|h\|_{\infty}\log n\,\, \bigg\{   \frac{u}{n} + \left(\frac{u}{n}\right)^{2}  \bigg\}.\]
\end{theorem}

\section{Estimation of spectra of signed integral operator with MCMC algorithms}
\label{sec:integral-ope}

\subsection{MCMC estimation of spectra of signed integral operators}

Let us consider a Markov chain~$(X_n)_{n \geq 1}$ on~$E$ satisfying the assumptions of Theorem~\ref{mainthm2} with stationary distribution~$\pi$, and some symmetric kernel~$h:E\times E \to \mathds R$ such that $h \in L^2(\pi \otimes \pi)$. We can associate to~$h$ the kernel of a linear operator~$\mathbf H$ defined by \begin{equation} \label{def:integral-ope}\mathbf H f(x) := \int_{E} h(x,y) f(y) d\pi(y).\end{equation}

This is a Hilbert-Schmidt operator on~$L^2(\pi)$ and thus it has a real spectrum consisting of a square summable sequence of eigenvalues~\cite[cf.][p.267]{conway2019course}. In the following, we will denote the eigenvalues of~$\mathbf H$ by~$\lambda(\mathbf H):=(\lambda_1, \lambda_2, \dots)$. For some~$n \in \mathds N^*$, we consider \begin{equation}\mathbf{\widetilde H}_n := \frac{1}{n} \left( h(X_i,X_j) \right)_{1\leq i,j\leq n}\text{ and }  \mathbf{ H}_n := \frac{1}{n} \left( (1-\delta_{i,j})h(X_i,X_j) \right)_{1\leq i,j\leq n},\label{eq:def-Hn} \end{equation}
with respective eigenvalues~$\lambda(\mathbf{ \widetilde H}_n)$ and~$\lambda(\mathbf{  H}_n)$. Following~\cite[Section 2]{gine2000}, we introduce in Definition~\ref{def:delta2} the rearrangement distance~$\delta_2$ which measures closeness of spectra.

\begin{definition}
\label{def:delta2}
Given two sequences~$x,y$ of reals -- completing finite sequences by zeros -- such that
\[
\sum_i x_i^2+y_i^2<\infty\,,
\]we define the~$\ell_2$ rearrangement distance~$\delta_2(x,y)$ as 
\[\delta_2^2(x,y):= \inf_{\sigma \in \mathfrak{S}}\sum_i (x_i-y_{\sigma(i)})^2\,,
\]
where~$\mathfrak{S}$ is the set of permutations of natural numbers. $\delta_2$ is a pseudometric
on $\ell_2$, where $\ell_2$ is the Hilbert space of all square summable sequences.
\end{definition}

Theorem~\ref{cor-integral-operator-bounded} gives conditions ensuring that both the spectrum of~$\mathbf{H}_n$ and the one of~$\mathbf{\widetilde H}_n$ converge towards the spectrum of the integral operator~$\mathbf{H}$ as~$n\to \infty.$ Theorem~\ref{cor-integral-operator-bounded} holds under Assumption~\ref{assumption:kernel-appli1} that we discuss in details in Section~\ref{sec:appli1-comparison}. The proof of Theorem~\ref{cor-integral-operator-bounded} is postponed to Section~\ref{sec:proof-appli1}.

\begin{assumption}\label{assumption:kernel-appli1}
$h:E \times E \to \mathds R$ is a bounded and symmetric function square integrable with respect to~$\pi \otimes \pi$. 
Moreover there exist continuous functions~$\phi_r:E \to \mathds R$,~$r\in I$ (where~$I=\mathds N$ or~$I=1, \dots,N$) that form an orthonormal basis of~$L^2(\pi)$ and a sequence of real numbers $(\lambda_r)_{r\in I} \in \ell_2$ such that we have pointwise \[h(x,y) = \sum_{r \in I} \lambda_r \phi_r(x) \phi_r(y),\] 
with $\displaystyle \Upsilon:= \sup_{r \in I } \|\phi_r\|_{\infty}^2 <\infty$ and $S:=\sup_{x \in E} \sum_{r\in I} |\lambda_r|\phi_r(x)^2<\infty.$\\
We further denote $\displaystyle \Lambda:= \sup_{r \in I} |\lambda_r|$.
\end{assumption}

\begin{theorem} \label{cor-integral-operator-bounded} Let $(X_i)_{i\geq 1}$ be a Markov chain on~$E$ satisfying Assumptions~\ref{assumption1} and~\ref{assumption2} described in Section~\ref{sec:concentration} with stationary distribution~$\pi$. Suppose that Assumption~\ref{assumption:kernel-appli1} is satisfied. Then for any~$t>0$,
\begin{align*}
\mathds P &\left(\frac14 \delta_2(\lambda(\mathbf{H}),\lambda( \mathbf{ H}_n))^2  \geq  \frac{S^2(1+\kappa) \log n}{n}+2\sum_{i>\ceil{ n ^{1/4}}, i\in I} \lambda_i^2+ t \right)\\
&\leq  32 \sqrt{n}\exp\left(- \mathcal C \min\left( n t^2 , \sqrt n t \right)\right)+\beta \log(n) \exp\left( - \frac{n}{\log n}\min\left( \mathcal B t, \left(\mathcal B t\right)^{1/2}\right) \right),
\end{align*}
where for some universal constant $K>0$, we have $\mathcal B = \left( K \kappa S\right)^{-1}$,~$\mathcal C =   \left( K^{1/2}m\tau (S + \Lambda \Upsilon)  \right)^{-2}$. $\kappa>0$ and~$\beta>0$ are the constants from Theorem~\ref{mainthm2} and depend on the Markov chain.  We refer to Assumption~\ref{assumption1} and the following remark for the definitions of the constants $m$ and $\tau$.
\end{theorem}
{\bf Remark} The same bound holds for~$\delta_2(\lambda(\mathbf{H}),\lambda( \mathbf{\widetilde H}_n))^2$.

\subsection{Comparison with the existing literature} \label{sec:appli1-comparison}

\paragraph{Previous works.}
In \cite{adamczak2015}, the authors studied the convergence properties of MCMC methods to estimate the spectrum of integral operators with bounded \textit{positive} kernels (i.e. such that $\mathbf H$ has non-negative eigenvalues). They show a sub-exponential tail behavior for the~$\delta_2$ distance between the spectrum of~$\mathbf H$ and the one of the random matrix~$\mathbf H_n$. Their result has the merit to hold for geometrically ergodic Markov chains, but they work with the restrictive assumption that the eigenvalues of~$\mathbf H$ are non-negative. This makes their kernel of positive-type~\citep[cf. Eq.(14)][]{adamczak2015}, allowing them to use of powerful decoupling argument. They are reduced to study a sum of rank one operators on~$L^2(\pi)$ and the proof is concluded by using the splitting technique~\citep[cf.][Section 5.1]{tweedie} with a Bernstein-type inequality for sums of independent Banach space valued random variables.

\paragraph{Comparing Theorem~\ref{cor-integral-operator-bounded} with previous works.}
In this paper we consider signed integral operators and we cannot adapt the proof proposed by~\cite{adamczak2015}. Working with stronger conditions on the Markov chain~$(X_i)_{i\geq1}$ compared to the former paper, Theorem~\ref{cor-integral-operator-bounded} proves a new concentration inequality for the~$\delta_2$ distance between~$ \lambda(\mathbf{H})$ and~$\lambda( \mathbf{ H}_n)$ that holds for \textbf{arbitrary signs of the eigenvalues} of~$\mathbf H.$ We provide a synthetic description of our proof at the beginning of Section~\ref{sec:proof-appli1}
. Note that the set of operators handled
by Theorem~\ref{cor-integral-operator-bounded} is not a superset of the ones handled by~\cite[Theorem 2.2]{adamczak2015}. The difference lies in the fact that we ask the family of functions $(\phi_r)_{r\in I}$ to be uniformly bounded (cf. Assumption~\ref{assumption:kernel-appli1}). Let us mention that the set of assumptions considered in~\cite{adamczak2015} implies that $S<\infty$. In the following, we comment our extra assumption $\Upsilon<\infty$ with more details.
\begin{enumerate}
\item The basis functions $(\phi_r)_{r\in I}$ are continuous and Assumption~\ref{assumption2} typically holds for a compact space $E$. Hence, by considering that $E$ is compact and that the sequence $(\lambda_r)_{r\in I}$ has finite support (i.e. $I=[N]$ for some natural number $N$), it holds $\Upsilon<\infty$. 
\item Asking for $\Upsilon<\infty$ is only useful to apply a concentration inequality for Markov chains at one specific step of our proof (cf. Eq.\eqref{1st-approach}). Hence this assumption might be weakened by applying other exponential tail control for empirical processes of Markov chains. Nevertheless we point out that Theorem~\ref{cor-integral-operator-bounded} holds for uniformly ergodic Markov chains which is equivalent to the standard {\it drift condition} where the drift function $V$ is bounded \citep[cf.][Chap.16]{tweedie}. Hence, the assumptions needed for the exponential inequality from~\cite[Theorem 1.1]{adamczak2015-concentration} or the one from~\cite[Theorem 5]{durmus2021probability} imply that $\Upsilon<\infty$. Hence, weakening the condition~$\Upsilon<\infty$ seems challenging but we believe that it might be done in some specific settings using for instance the work from~\cite[Section~3.2]{bertail2018new}. 
\end{enumerate}

\subsection{Admissible sampling schemes: Independent Hastings algorithm}
One can use the previous result to estimate the spectrum of the integral operator~$\mathbf H$ using MCMC methods. To do so, we need to make sure that the Markov chain used for the MCMC method satisfies the conditions of Theorem~\ref{mainthm2}. It is for example well known that Metropolis random walks on~$\mathds R$ are not uniformly ergodic \citep[cf.][]{tweedie}. In the following, we show that an independent Hastings algorithm can be used on bounded state space to generate a uniformly ergodic chain with the desired stationary distribution.

\subsubsection{Independent Hastings algorithm on bounded state space.}

Let~$E\subset \mathds R^k$ be a bounded subset of~$\mathds R^k$ equipped with the Borel~$\sigma$-algebra~$\mathcal B(E)$. We consider a density~$\pi~$ which is only known up to a factor and a probability density~$q$ with respect to the Lebesgue measure~$\lambda_{Leb}$ on~$E$, satisfying~$\pi(y),q(y)>0~$ for all~$y \in E$. In the independent Hastings algorithm, a candidate transition generated according to the law~$q\lambda_{Leb}$ is then accepted with probability~$\alpha(x,y)$ given by
\[\alpha(x,y) = \min \left( 1, \frac{ \pi(y)q(x)}{\pi(x)q(y)}  \right).\]
With an approach similar to Theorem 2.1 from \cite{MT96}, Proposition~\ref{prop:inde-hasting} shows that under some conditions on the densities~$\pi$ and~$q$, the independent Hastings algorithm satisfies the Assumptions~\ref{assumption1} and~\ref{assumption2}. 

\begin{proposition} \label{prop:inde-hasting} 
Let us assume that~$\displaystyle\sup_{x \in E} q(x) < \infty$ and that there exists~$\beta>0$ such that \[\frac{q(y)}{\pi(y)}>\beta, \quad \forall y \in E.\]
Then, the independent Hastings algorithm satisfies the Assumptions~\ref{assumption1} and~\ref{assumption2}. 
\end{proposition}

\noindent
\begin{proof}
We denote~$P$ the transition kernel of the Markov chain generated with the independent Hastings algorithm. For any~$x \in E$, the density with respect to~$\lambda_{Leb}$ of the absolutely continuous part of~$P(x,dy)$ is~$p(x,\cdot)=q(\cdot) \alpha(x,\cdot)$, while the singular part is given by~$\mathds 1_x(\cdot) \left( \int_{z \in E} q(z) \alpha(x,z)d\lambda_{Leb}(z)  \right)$. For fixed~$x\in E$, we have either~$\alpha(x,y)=1$ in which case~$p(x,y)=q(y)\geq \beta \pi(y)$, or else \[p(x,y) = q(y)\frac{ \pi(y)q(x)}{\pi(x)q(y)} =  q(x) \frac{\pi(y)}{\pi(x)} \geq \beta \pi(y).\]
We deduce that for any~$x \in E$, it holds \[P(x,A) \geq \beta \int_{y\in A} \pi(y) d\lambda_{Leb}(y),\] which proves that the chain is uniformly ergodic (cf. Eq.\eqref{uni-ergodicity-small}). Hence Assumption~\ref{assumption1} is satisfied. Assumption~\ref{assumption2} trivially holds since~$E$ is bounded and~$\sup_{ y \in E} q(y) < \infty.$ 
\end{proof}

\medskip
From Proposition~\ref{prop:inde-hasting} and Theorem~\ref{cor-integral-operator-bounded}, we deduce that one can use a MCMC approach to estimate the spectrum of a signed integral operator $\mathbf H$ (that satisfies Assumption~\ref{assumption:kernel-appli1}) as defined in Eq.\eqref{def:integral-ope} where $E$ is a bounded subset of $ \mathds R^k$. More precisely, if the density $\pi$ of Eq.\eqref{def:integral-ope} is known up to a factor and if there exists some probability density $q$ with respect to $\lambda_{Leb}$ satisfying the assumptions of Proposition~\ref{prop:inde-hasting}, the Independent Hastings algorithm provides a Markov chain that satisfies Assumptions~\ref{assumption1} and~\ref{assumption2}. Hence the non-asymptotic bound from Theorem~\ref{cor-integral-operator-bounded} holds. We put this methodology into action in the new section by estimating the spectrum of some Mercer kernels on the~$d$-dimensional sphere.

\subsection{Estimation of the spectrum of Mercer kernels}

In this example, we illustrate Theorem~\ref{cor-integral-operator-bounded} by computing the eigenvalues of an integral operator naturally associated with a Mercer kernel using a MCMC algorithm. 
A function~$h:E \times E \to \mathds R$ is called a Mercer kernel if~$E$ is a compact metric space and~$h: E \times E \to \mathds R$ is a continuous symmetric and positive definite function. It is well known that if~$h$ is a Mercer kernel, then the integral operator~$L_h$ associated with $h$ is a compact and bounded linear operator, self-adjoint and semi-definite positive. The spectral theorem implies that if~$h$ is a Mercer kernel, then there is a complete orthonormal system~$(\phi_1,\phi_2,\dots)$ of eigenvectors of~$L_h$. The eigenvalues~$(\lambda_1,\lambda_2,\dots)$ are real and non-negative. The Mercer Theorem \citep[see e.g.][Theorem~4.49]{Mercer} shows that the eigen-structure of~$L_h$ can be used to get a representation of the Mercer kernel~$h$ as a sum of a convergent sequence of product functions for the uniform norm. In this context, Theorem~\ref{cor-integral-operator-bounded} allows to derive the convergence rate in the $\delta_2$ metric of the estimated spectrum towards the one of the integral operator $\mathbf H$ as presented in Proposition~\ref{prop:rate-mcmc}.

\begin{proposition} \label{prop:rate-mcmc} We keep the notations and the assumptions of Theorem~\ref{cor-integral-operator-bounded}. We assume further that there exists $s >0$, a (Sobolev) regularity parameter, such that for some constant $C(s)>0$,
\[ \forall R>1, \quad \sum_{i>R} \lambda_i^2 \leq C(s) R^{-2s}.\]
Then it holds\[
 \delta_2(\lambda(\mathbf{H}),\lambda( \mathbf{ H}_n))^2 =\left\{
    \begin{array}{ll}
       \mathcal O_{\mathds P}\left(\sqrt{\frac{\log n }{n}}\right) & \mbox{if } s\geq 1 \\
         \mathcal O_{\mathds P}\left(\frac{1 }{n^{s/2}}\right) & \mbox{if } s \in (0,1)
    \end{array}
\right..
\]
\end{proposition}
\begin{proof}Proposition~\ref{prop:rate-mcmc} directly follows from Theorem~\ref{cor-integral-operator-bounded} by choosing $t=\sqrt{\frac{\log n }{n}}$. \end{proof}
\medskip

\noindent
To illustrate our purpose, we consider the~$d$-dimensional sphere~$\mathds S^{d-1} =\{ x \in \mathds R^d : \|x\|_2=1\}$. We consider a positive definite kernel on~$\mathds S^{d-1}$ defined by~$\forall x,y \in \mathds S^{d-1}, \quad h(x,y) = \psi(x^{\top}y)$ where~$\psi:[-1,1] \to \mathds R$ is continuous. From the Funk-Hecke Theorem \citep[see e.g.][p.30]{muller2012analysis}, we know that the eigenvalues of the Mercer kernel~$h$ are \begin{equation}\label{mercer-eigs}\lambda_k = |\mathds S^{d-2}| \int_{-1}^1 \psi(t)P_k(d;t)\left(1-t^2\right)^{\frac{d-3}{2}}dt,\end{equation}
where~$|\mathds S^{d-2}|$ is the surface area of $\mathds S^{d-2}$ and~$P_k(d;t)$ is the Legendre polynomial of degree~$k$ in dimension~$d$. For any $k \in \mathds N$, the multiplicity of the eigenvalue $\lambda_k$ is the dimension of the space of spherical harmonics of degree $k$. To build the Markov chain~$(X_i)_{i\geq 1}$, we start by sampling randomly~$X_1$ on~$\mathds{S}^{d-1}$. Then, for any~$i \in \{2, \dots,n\}$, we sample
\begin{itemize}
\item a unit vector~$Y_i \in \mathds{S}^{d-1}$ uniformly, orthogonal to~$X_{i-1}$.
\item a real~$r_i \in [-1,1]$ encoding the distance between~$X_{i-1}$ and~$X_i$.~$r_i$ is sampled from a distribution~$f_{\mathcal{L}}:[-1,1] \to [0,1]$.
\end{itemize}
then~$X_i$ is defined by 
\begin{equation*}
X_i = r_i \times X_{i-1} + \sqrt{1-r_i^2}\times  Y_i\,.
\end{equation*}
By assuming that~$ \min_{r \in [-1,1]} f_{\mathcal L}(r) >0$ and~$\|f_{\mathcal L}\|_{\infty} <\infty$, Assumptions~\ref{assumption1} and~\ref{assumption2} hold and the stationary distribution of the chain $(X_i)_{i\geq1}$ is the Haar measure on $\mathds S^{d-1}$ \citep[cf.][]{DCD20}.

In Figure~\ref{fig:spectra}, we plot the non-zero eigenvalues using function~$\psi:t \mapsto (1+t)^2$ and taking~$f_{\mathcal{L}}$ proportional to~$r\mapsto f_{(5,1)}(\frac{r+2}{4})$ where~$f_{(5,1)}$ is the pdf of the Beta distribution with parameter~$(5,1)$. We plot both the true eigenvalues and the ones computed using a MCMC approach.

\begin{minipage}{0.5\textwidth}
\centering
\captionsetup{type=figure}
\includegraphics[scale=0.48]{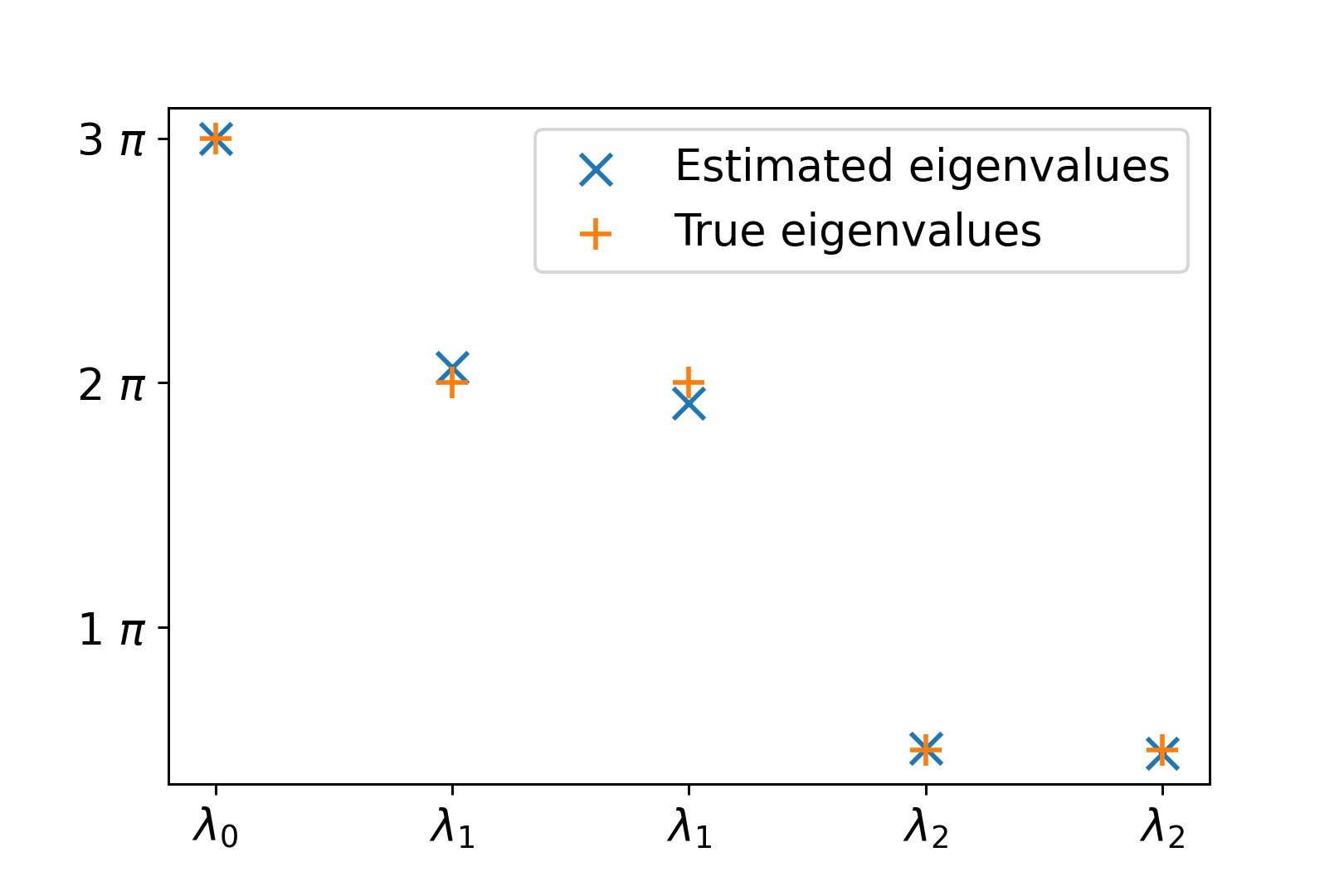}
\end{minipage} 
\begin{minipage}{0.45\textwidth}
\captionof{figure}{Consider function~$\psi:t \mapsto (1+t)^2$, $d=2$ and $n=1000$. The true eigenvalues can be computed using \eqref{mercer-eigs}, but in this case, we know the exact values of the three non-zero eigenvalues namely $\lambda_0=3\pi$, $\lambda_1=2\pi$ and $\lambda_2= \pi/2$. Their respective multiplicities are $1$, $2$ and $2$. The estimated eigenvalues are the eigenvalues of the matrix~$\mathbf H_n = \frac{1}{n} \left( (1-\delta_{i,j})\psi(X_i^{\top}X_j) \right)_{1\leq i,j\leq n}$ where the~$n$ points~$X_1, X_2, \dots, X_n$ are sampled on the Euclidean sphere~$\mathds{S}^{d-1}$ using a Markovian dynamic.}
\label{fig:spectra}
\end{minipage}

\section{Online Learning with Pairwise Loss Functions}
\label{sec:online-pairwise}
\subsection{Brief introduction to online learning and motivations}

\subsubsection{Presentation of the traditional online learning setting}

 Online learning is an active field of research in Machine Learning in which data becomes available in a sequential order and is used to update the best predictor for future data at each step. This method aims at learning some function~$f:E \to \mathcal Y$ where~$E$ is the space of inputs and~$ \mathcal Y$ is the space of outputs. At each time step~$t$, we observe a new example~$(x_t,y_t) \in E\times \mathcal Y$. Traditionally, the random variables~$(x_t,y_t)$ are supposed i.i.d. with common joint probability distribution~$(x,y)\mapsto p ( x , y )$ on~$E\times \mathcal Y$. In this setting, the loss function is given as~$\ell : \mathcal Y \times \mathcal Y \to \mathds R$, such that~$\ell( f ( x ) , y )$ measures the difference between the predicted value~$f ( x )$ and true value~$y$. The goal is to select at each time step~$t$ a function~$h_t:E \to \mathcal Y$ in a fixed set~$\mathcal H$ based on the observed examples until time~$t$ (namely~$(x_i,y_i)_{1\leq i \leq t}$) such that~$h_t$ has ‘‘small'' 
 risk~$\mathcal R$ defined by
 \[
 \mathcal R (h)= \mathds E_{(X,Y) \sim p}\big[ \ell(h(X),Y) \big]\,,
 \]
 where~$h$ is any measurable mapping from~$E$ to~$\mathcal Y$.

 Online learning is used when data is coming {\it on the fly} and we do not want to wait for the acquisition of the complete dataset to take a decision. In such cases, online learning algorithms allow to  dynamically adapt to new patterns in the data.
 
 \subsubsection{Online learning with pairwise loss functions}
 In some cases, the framework provided in the previous paragraph is not appropriated to solve the task at stake. Consider the example of ranking problems. The state space is~$E$ and there exists a function~$f:E \to \mathds R$ which assigns to each state~$x \in E$ a label~$f(x) \in \mathds R$.~$f$ naturally defines a partial order on~$E$. At each time step~$t$, we observe an example~$x_t \in E$ together with its label~$f(x_t)$ and we suppose that the random variables~$(x_t)_t$ are i.i.d. with common distribution~$p$. Our goal is to learn the partial order of the items in~$E$ induced by the function~$f$. More precisely, we consider a space~$\mathcal H \subset \{ h:E\times E \to \mathds R\}$, called the set of hypotheses. An {\it ideal} hypothesis~$h \in \mathcal H$ would satisfy 
 \[
 \forall x,u \in E, \quad f(x)\geq f(u) \Leftrightarrow \left( h(x,u)\geq 0 \text{ and }  h(u,x)\leq 0\right).
 \] 
 We consider a loss function~$\ell:\mathcal H \times E \times E \to \mathds R$ such that~$\ell(h,x,u)$ measures the ranking error induced by~$h$ and a typical choice is the~$0$-$1$ loss 
 \[
 \ell(h,x,u)=\mathds 1_{\{(f(x)-f(u))h(x,u)<0\}}.
 \]
U-statistics naturally arise in such settings as for example in \cite{clemencon2008} where Clémençon and al. study the consistency of the empirical risk minimizer of ranking problems using the theory of U-processes in an i.i.d. framework.

\medskip

\pagebreak[3]

\noindent
\textbf{Example}: {\tt Bipartite ranking problems}

{\it 
We describe the concrete problem of bipartite ranking. We consider that we have as input a training set of examples. Each example is described by some feature vector and is associated with a binary label. Typically one can consider that we have access to health data of an individual along time. We know at each time step her/his health status~$x_t$ and her/his label which is 0 if the individual is healthy and 1 if she/he is sick. In the bipartite ranking problem, we want to learn a \texttt{scorer} which maps any feature vector describing the health status of the individual to a real number such that sick states have a higher score than healthy ones. Following the health status of individuals is time-consuming and we cannot afford to wait for the end of the data acquisition process to understand the relationship between the feature vector describing the health status of the individual and her/his sickness. In such settings where data is coming on the fly, online algorithms are common tools that allow to learn a scorer function along time. At each time step the scorer function is updated based on the new measurement provided.}

\subsubsection{Generalization bounds for online learning}

The performance of online learning algorithms is often analyzed with the notion of {\it regret} which compares the payoff obtained by the algorithm along time with the one that would have been obtained by taking the optimal decision at each time step \citep[cf.][]{onlinesurvey,bubeck}. It is natural to wonder if stronger theoretical guarantees can be obtained when some probabilistic
structure underlies the sequence of examples, or loss functions, presented to the online algorithm.
As asked in~\cite{agarwal2012}, {\it ‘if the sequence of examples are generated by a stochastic process, can the online learning algorithm output a good predictor for future samples from the same process?'} In other words, we want to study the generalization ability of some online learner that generates a sequence of hypothesis~$(h_t)_{t\geq1}$ by bounding with high probability the excess risk defined as
\[\frac{1}{n}\sum_{t=1}^n\mathcal R(h_t)-\min_{h\in \mathcal H} \mathcal R(h).\]Generalization bounds for online learning with pairwise loss functions working with i.i.d. samples have been considered for
quite a while in both Machine Learning and Statistics literature~\citep[cf.][]{kar13,ying2017,guo2017,chen2018}. For dependent data sequences, generalization bounds for online algorithms have also been proved in the last decades with univariate loss functions~\citep[cf.][]{zhang2005data,xu2014generalization,agarwal2012}. However, theoretical guarantees for the generalization performance of online algorithms with pairwise loss functions with non i.i.d. data have been so far little studied. A quick and incomplete review of the literature is presented in Table~\ref{table:lit-online}.
\begin{table}[!ht]
\centering
\begin{tabular}{c||C{5cm}|C{5cm}}
& \cellcolor{gray!40} Univariate loss function & \cellcolor{gray!40}  Pairwise loss function \\\hline\hline
\cellcolor{gray!40} i.i.d. data &   \citeauthor[Section 3.7]{onlinesurvey} and references therein &  \citeauthor{kar13,ying2017,guo2017,chen2018}\\\hline
\cellcolor{gray!40} Dependent data &  \citeauthor{zhang2005data,xu2014generalization,agarwal2012}  & Our work
\end{tabular}
\caption{Overview of the literature providing generalization bounds for online learning algorithms.}
\label{table:lit-online}
\end{table}

\subsubsection{Generalization bounds for pairwise online learning with dependent data}

\paragraph{Connection with the existing literature.}
\label{sec:appli2-literature}

As far as we know, the few papers that investigate the generalization performance of pairwise online learning algorithms with non i.i.d. data have studied specific algorithms and/or specific learning tasks~\citep[cf.][]{qin2021online,zeng2021}. In~\cite{zeng2021}, the authors analyze online pairwise support vector machine while the work~\cite{qin2021online} is focused on online regularized pairwise learning algorithm with least squares loss function. One possible reason explaining this gap in the literature is that {\it ‘for pairwise learning [where] pairs of training examples are not i.i.d., [...] standard techniques can not be directly applied.'} \citep[cf.][]{zeng2021}. 

With the upcoming application, we are the first - as far as we know - to provide a generalization bounds for online algorithms with pairwise loss functions and Markov chain samples that hold for an arbitrary online learner, covering a large span of settings.

\paragraph{Online learning with a Markovian dynamic.}

The theoretical analysis of Machine Learning algorithms with an underlying Markovian distribution of the data has become a very active field of research. The first papers to study online learning with samples drawn from non-identical distributions were~\cite{smale2009online} and~\cite{steinwart2009learning} where online learning for least square regression and off-line support vector machines are investigated. In~\cite{zou2009learning}, the generalization performance of the empirical risk minimization algorithm is studied with uniformly ergodic Markov chain samples. Hence the analysis of online algorithms with dependent samples is recent and several works make the assumption that the sequence is a uniformly ergodic Markov chain. 
We motivate the Markovian assumption on the example of the previous paragraph.

\medskip

\noindent
\textbf{Example (continued)}: {\tt Interest in online algorithms with Markovian dynamic}

{\it 
The health status of the individual at time~$n+1$ is not independent from the past and a simple way to model this time evolution would be to consider that it only depends on the last measured health status namely the feature vector~$x_n$. This is a Markovian assumption on the sequence of observed health status of the individual. } 

\medskip

We have explained why pairwise loss functions capture ranking problems and naturally arise in several Machine Learning problems such as metric learning or bipartite ranking (cf. \cite{clemencon2008}). We have shown the interest to provide a generalization bounds for online learning with pairwise loss functions with a Markovian assumption on the distribution of the sequence of examples and this is the goal of the next section.

\subsection{Online-to-batch conversion for pairwise loss functions with Markov chains}

We consider a reversible Markov chain~$(X_i)_{i\geq 1}$ with state space~$E$ satisfying Assumption~\ref{assumption1} with stationary distribution~$\pi$. Using~\cite[Theorem 16.0.2]{tweedie}, we deduce that there exist constants~$0<\rho<1$ and~$L>0$ such that \begin{equation} \|P^n(x,\cdot)-\pi \|_{\mathrm{TV}} \leq L \rho^n, \qquad \forall n \geq0, \; \pi\mathrm{-a.e}\;  x \in E. \label{uni-ergodicity} \end{equation}
We assume that we have a function~$f : E \to \mathds R$ which defines the ordering of the objects in~$E$. We aim at finding a relevant approximation of the ordering of the objects in~$E$ by selecting a function~$h$ (called a \textit{hypothesis} function) in a space~$\mathcal H$ based on the observation of the random sequence~$(X_i,f(X_i))_{1 \leq i \leq n}$. To measure the performance of a given hypothesis~$h: E\times E\to \mathds R$, we use a pairwise loss function of the form~$\ell(h,X,U)$. Typically, one could use the \textit{misranking loss} defined by \[\ell(h,x,u)=\mathds 1_{\{(f(x)-f(u))h(x,u)<0\}},\] which is~$1$ if the examples are ranked in the wrong order and~$0$ otherwise. The goal of the learning problem is to find a hypothesis~$h$ which minimizes the \textit{expected misranking risk}
\[
\mathcal R(h) := 
\mathds E_{(X,X') \sim \pi\otimes\pi}\big[ \ell(h,X,X') \Big]. 
\]

\pagebreak[3]

We show that the investigation of the generalization performance of online algorithms with pairwise loss functions provided by \cite{online-pairwise} can be extended to a Markovian framework. 
Our contribution is two-fold. 
\begin{itemize}
\item Firstly, we prove that with high probability, the average risk of the sequence of hypotheses generated by an arbitrary online learner is bounded by some easily computable statistic. 
\item This first technical result is then used to show how we can extract a low risk hypothesis from a given sequence of hypotheses selected by an online learner. This is an {\it online-to-batch} conversion for pairwise loss functions with a Markovian assumption on the distribution of the observed states.
\end{itemize}

\noindent
Given a sequence of hypotheses~$(h_i)_{1\leq i \leq n} \in \mathcal H^n$ generated by any online algorithm, we define the {\it average paired empirical risk} $\mathcal M^n(h_1,\dots,h_{n-1-b_n})$ (see Eq.\eqref{eq:aer}) averaging the {\it paired empirical risks} $M_t$ (see Eq.\eqref{eq:er}) of hypotheses $h_{t-b_n}$ when paired with $X_t$ as follows
\begin{align}
    \mathcal M^n(h_1,\dots,h_{n-1-b_n})&:= \frac{1}{n-c_n}\sum_{t=c_n}^{n-1}M_t,\label{eq:aer}\\
    \text{and}\quad M_t&:=\frac{1}{t-b_n} \sum_{i=1}^{t-b_n}\ell(h_{t-b_n},X_t,X_i),\label{eq:er}
\end{align}
where 
\begin{equation}
c_n = \ceil{c \times n} \text{ for some }c \in (0,1)\quad\text{ and }\quad b_n  = \floor{q \log(n)}, \label{q-bn}
\end{equation} 
for an arbitrarily chosen~$q> \frac{1}{\log(1/\rho)}$ where~$\rho$ is a constant related to the uniform ergodicity of the Markov chain, see Eq.\eqref{uni-ergodicity}. In the following, we will simply denote~$\mathcal M^n(h_1,\dots,h_{n-1-b_n})$ by~$\mathcal M^n$ when the sequence of considered hypotheses is clear from the context.

$M_t$ is the {\it paired empirical risk} of hypothesis $h_{t-b_n}$ with $X_t$. It measures the performance of the hypothesis~$h_{t-b_n}$ on the example~$X_t$ when paired with examples seen before time~$t-b_n$. $\mathcal M^n$ is the mean value of a proportion $1-c$ of these paired empirical risks. Hence the parameter $c\in (0,1)$ controls the proportion of hypotheses~$h_{t-b_n}$ whose paired empirical risk $M_t$ does not appear in the average paired empirical risk value $\mathcal M^n$. The parameter $b_n$ controls the time gaps between elements of pairs $(X_t,X_i)$ appearing in Eq.\eqref{eq:er} in such way that their joint law is close to the product law $\pi\otimes\pi$ (mixing of the chain is met). The use of the burning parameter $b_n$ is the main difference with the work~\cite{online-pairwise} when defining $\mathcal M^n$ and $M_t$ in Eq.\eqref{eq:aer} and Eq.\eqref{eq:er}. From a pragmatic point of view,

\begin{itemize}
\item we discard the first hypotheses that are not reliable, namely we do not consider hypothesis~$h_{i}$ for~$i \leq c_n-b_n$. These first hypotheses are considered as not reliable since the online learner selected them based on a too small number of observed examples.
\item since~$h_{t-b_n}$ is learned from~$X_1, \dots ,X_{t-b_n}$, we test the performance of~$h_{t-b_n}$ on~$X_t$ (and not on some~$X_i$ with~$t-b_n+1\leq i <t~$) to ensure that the distribution of~$X_t$ conditionally on~$\sigma(X_1, \dots, X_{t-b_n})$ is approximately the stationary distribution of the chain~$\pi$ (see Assumption~\ref{assumption1} and Equation~\eqref{uni-ergodicity}). Stated otherwise, this ensures that sufficient mixing has occurred.
\end{itemize}
Note that we assume~$n$ large enough to ensure that~$c_n-b_n \geq 1.$ For any~$\eta >0$, we denote~$\mathcal N(\mathcal H, \eta)$ the~$L^{\infty}$~$\eta$-covering number for the hypothesis class~$\mathcal H$ (see Definition~\ref{def:covering-number}).

\begin{definition}  \citep[cf.][Chapter~5.1]{wainwright2019high}
\label{def:covering-number}
Let us consider some $\eta>0$. A $L^{\infty}$~$\eta$-cover of a set $\mathcal H$ is a set $\{g_1, \dots , g_N\} \subset \mathcal H$ such that for any $h \in \mathcal H$, there exists some $i \in \{1, \dots ,N\}$ such that $\|g_i-h\|_{\infty} \leq \eta.$  The~$L^{\infty}$~$\eta$-covering number~$\mathcal N(\mathcal H, \eta)$ is the cardinality of the smallest~$L^{\infty}$~$\eta$-cover of the set~$\mathcal H$.
\end{definition}

\noindent
Theorem~\ref{pairwise-thm1} bounds the average risk of the sequence of hypotheses in terms of its empirical counterpart~$\mathcal M^n$ and is proved in Section~\ref{proof-pairwise-thm1}.

\begin{theorem} \label{pairwise-thm1}
Assume that the Markov chain~$(X_i)_{i\geq 1}$ is reversible and satisfies Assumption~~\ref{assumption1}. Assume the hypothesis space~$(\mathcal H, \|\cdot \|_{\infty})$ is compact. Let~$h_0, h_1,\dots, h_n\in  \mathcal H$ be the ensemble of hypotheses generated by an arbitrary online algorithm working with a pairwise loss function~$\ell$ such that, \[\ell(h,x_1,x_2) =\phi(f(x_1)-f(x_2), h(x_1,x_2)),\] where~$\phi:\mathds R \times \mathds R \to [0,1]$ is a Lipschitz function w.r.t. the second variable with a finite Lipschitz constant~$Lip(\phi)$. Let~$\xi>0$ be an arbitrary positive number and let us consider~$q= \frac{\xi+1}{\log(1/\rho)}$ for the definition of~$b_n$ (see Eq.\eqref{q-bn}). Then for all~$c>0$ and for all~$\epsilon >0$ such that~$\epsilon \underset{n\to \infty}{=} o\left( n^{\xi}\right)$, we have for sufficiently large~$n$
\begin{align*}
\mathds P\left(
\Big|
\frac{1}{n-c_n}\sum_{t=c_n}^{n-1}\mathcal R(h_{t-b_n}) - \mathcal M^n 
\Big|
\geq \epsilon \right) \leq
&
2
\left[  32\mathcal N\left(\mathcal H, \frac{\epsilon}{8 Lip(\phi)}\right) +1 \right] b_n \\
&\times \exp \left(  - \frac{(c_n-b_n)  C(m,\tau)\epsilon^2}{16b_n^2} \right),
\end{align*}
where~$C(m,\tau)^{-1}=7\times 10^3\times m^2\tau^2$. We refer to Assumption~\ref{assumption1} and the following remark \citep[or to][Section 2]{duchemin} for the definitions of the constants $m$ and $\tau$ that depend on the Markov chain $(X_i)_{i\geq1}$.
\end{theorem}

\noindent

Theorem~\ref{pairwise-thm1} shows that average paired empirical risk $\mathcal M^n$ (see Eq.\eqref{eq:aer}) is close to average risk given by 
\[
\frac{1}{n-c_n}\sum_{t=c_n}^{n-1}\mathcal R(h_{t-b_n})\,.
\]
Quantitative errors bounds can be given assuming that the $L_\infty$-metric entropy (l.h.s of the next equation) satisfies 
\begin{equation}\label{eq:covering-nb}
\log \mathcal N(\mathcal H, \eta)=\mathcal O(\eta^{-\theta})\,,
\end{equation}
where $\theta>0$ is an exponent, depending on the dimension of state space~$E$ and the regularity of hypotheses of~$\mathcal H$, that can be computed in some situations (Lipschitz function, higher order smoothness classes, see \cite[Chapter~5.1]{wainwright2019high} for instance). Theorem~\ref{pairwise-thm1:regret} made this statement rigorous (cf. Eq.\eqref{regret-bound1}).

As previously mentioned, online learning algorithms are often studied through the lens of {\it regret}. The definition of a regret bound in our context is provided in Definition~\ref{def:regret-bound}.
\begin{definition}\label{def:regret-bound}
An online learning algorithm will be said to have a regret bound~$\mathfrak R_n$ if it presents an ensemble
$h_1,\dots , h_{n-1}$ such that
\[\mathcal M^n  \leq \min_{h \in \mathcal H} \big\{\mathcal M^n(h,\dots,h)  \big\}+ \mathfrak R_n.\]
\end{definition}
In the literature of learning theory~\cite{cucker2007learning}, we are often interested in the averaged excess generalization error
\[\frac{1}{n-c_n}\sum_{t=c_n}^{n-1}\mathcal R(h_{t-b_n}) - \mathcal R(h^*),\]
where~$h^*$ is the population risk minimizer and is given by~$h^* \in \underset{h\in\mathcal H}{\arg \min} \; \mathcal R(h)$.

As a consequence, most of works focused on online-to-batch conversion are interested in the overall convergence rate of the excess generalization error for online learners that achieve a given regret bound. Examples can be found with~\cite[Corollary 4]{guo2017} or with~\cite[Theorem 5]{kar13} where both papers work with pairwise loss functions with i.i.d.  observations. In Theorem~\ref{pairwise-thm1:regret} (cf. Eq.\eqref{regret-bound2}) we provide the overall rate for the averaged excess generalization error for an online learning satisfying a given regret bound. Theorem~\ref{pairwise-thm1:regret} is proved in Section~\ref{proof-pairwise-thm1:regret} and should be understood as an extension of the above mentioned results from~\cite{kar13} and~\cite{guo2017}.
\begin{theorem} \label{pairwise-thm1:regret}We keep the notations and assumptions of Theorem~\ref{pairwise-thm1}. Assume further that $\mathcal H$ satisfies Eq.\eqref{eq:covering-nb}. Then it holds 
\begin{equation}\label{regret-bound1}
\Bigg|
\frac{1}{n-c_n}\sum_{t=c_n}^{n-1}\mathcal R(h_{t-b_n}) - \mathcal M^n 
\Bigg|
=\mathcal O_{\mathds P}\Bigg[\frac{\log(n)\log(\log n)}{n^{\frac1{2+\theta}}}\Bigg]\,.
\end{equation}
Moreover, if the online learner has a regret bound $\mathfrak R_n$ (cf. Definition~\ref{def:regret-bound}), it holds
\begin{equation}\label{regret-bound2}
\frac{1}{n-c_n}\sum_{t=c_n}^{n-1}\mathcal R(h_{t-b_n}) - \mathcal R(h^*)
=\mathcal O_{\mathds P}\Bigg[\frac{\log(n)\log(\log n)}{n^{\frac1{2+\theta}}}+\mathfrak R_n\Bigg]\,.
\end{equation}
\end{theorem}

\subsection{Batch hypothesis selection}
Theorems~\ref{pairwise-thm1} and~\ref{pairwise-thm1:regret} are results on the performance of online learning algorithms. We will use these results to study the generalization performance of such online algorithms in the batch setting (see Theorem~\ref{pairwise-thm9}). Hence we are now interested in {\it selecting a good hypothesis from the ensemble of hypotheses generated by the online learner} namely that has a small empirical risk. 

We measure the risk for~$h_{t-b_n}$ on the last~$n-t$ examples of the sequence~$X_1, \dots, X_n$, and penalize each~$h_{t-b_n}$ based on the number of examples on which it is evaluated. More precisely, let us define the empirical risk of hypothesis~$h_{t-b_n}$ on~$\{X_{t+1},\dots,X_n\}$ as \[\widehat{\mathcal R}(h_{t-b_n}, t+ 1):= \binom{n-t}{2}^{-1} \sum_{k>i, i \geq t+1}^{n} \ell(h_{t-b_n},X_i,X_k).\]

For a confidence parameter~$\gamma \in (0,1)$ that will be specified in Theorem~\ref{pairwise-thm9}, the hypothesis~$\widehat h$ is chosen to minimize the following penalized empirical risk,
\begin{equation} \label{h-selection}
\widehat h= h_{\widehat t-b_n}\quad\text{and}\quad  \widehat t\in \arg \min_{c_n\leq t\leq n-1}\left( \widehat{\mathcal R}(h_{t-b_n}, t+ 1) +c_{\gamma}(n-t)\right),\end{equation}
\noindent
where \[c_{\gamma}(x) = \sqrt{\frac{C(m,\tau)^{-1}}{x} \log \frac{64(n-c_n)(n-c_n+1)}{\gamma}},\]
with~$C(m,\tau)^{-1}=7\times 10^3\times m^2\tau^2$.

 Theorem~\ref{pairwise-thm9} proves that the model selection mechanism previously described select a hypothesis~$\widehat h$ from the hypotheses of an arbitrary online learner whose risk is bounded relative to~$\mathcal M^n$. The proof of Theorem~\ref{pairwise-thm9} is postponed to Section~\ref{proof-pairwise-thm9}.

\begin{theorem} \label{pairwise-thm9}
Assume that the Markov chain~$(X_i)_{i\geq 1}$ is reversible and satisfies Assumptions~\ref{assumption1} and~\ref{assumption2}.\\Let~$h_0,\dots, h_n$ be the set of hypotheses generated by an arbitrary online algorithm~$\mathcal A$ working with a pairwise loss~$\ell$ which satisfies the conditions given in Theorem~\ref{pairwise-thm1}. 
Let~$\xi>0$ be an arbitrary positive number and let us consider~$q= \frac{\xi+1}{\log(1/\rho)}$ for the definition of~$b_n$ (see Eq.\eqref{q-bn}). For all~$ \epsilon>0$ such that~$\epsilon \underset{n\to \infty}{=} o\left( n^{\xi}\right)$, if the hypothesis is selected via Eq.\eqref{h-selection} with the confidence~$\gamma$ chosen as \[\gamma= 64(n-c_n+ 1) \exp\left(-(n-c_n)\epsilon^2 C(m,\tau)/128\right),\]then, when~$n$ is sufficiently large, we have
\[  \mathds P\left(  \mathcal R(\widehat h) \geq \mathcal M^n+\epsilon \right) \leq 32\left[  \mathcal N\left(\mathcal H, \frac{\epsilon}{16 \mathrm{Lip}(\phi)}\right)+1 \right]  \exp \left(  - \frac{(c_n-b_n)  C(m,\tau)\epsilon^2}{(16b_n)^2}+2\log n \right).\]
\end{theorem}

Analogously to the previous section, we can derive from Theorem~\ref{pairwise-thm9} a bound for the excess risk of the selected hypothesis~$\widehat h$.
\begin{Cor} \label{pairwise-cor9}
We keep the notations and assumptions of Theorem~\ref{pairwise-thm9}. Assume further that $\mathcal H$ satisfies Eq.\eqref{eq:covering-nb}. Then it holds 
\[
\Bigg|
\mathcal R(\widehat h) - \mathcal M^n 
\Bigg|
=\mathcal O_{\mathds P}\Bigg[\frac{\log^{2}n}{n^{\frac1{2+\theta}}}\Bigg]\,.
\]
Moreover, if the online learner has a regret bound $\mathfrak R_n$ (cf. Definition~\ref{def:regret-bound}), it holds
\[
\mathcal R(\widehat h) - \mathcal R(h^*)
=\mathcal O_{\mathds P}\Bigg[\frac{\log^{2}n}{n^{\frac1{2+\theta}}}+\mathfrak R_n\Bigg]\,.
\]
\end{Cor}

\clearpage
\section{Adaptive goodness-of-fit tests in a density model}
\label{sec:adaptivegof}

\subsection{Goodness-of-fit tests and review of the literature}
In its original formulation, the goodness-of-fit test aims at determining if a given distribution~$q$ matches some unknown distribution~$p$ from samples~$(X_i)_{i\geq 1}$ drawn independently from~$p$. Classical approaches to solve the goodness-of-fit problem use the empirical process theory. Most of the popular tests such as the Kolmogorov-Smirnov, Cramer-von Mises, and Anderson-Darling statistics are based on the empirical distribution function of the samples. Other traditional approaches may require space partitioning or closed-form integrals \cite{baringhaus1988consistent}, \cite{Lugosi08}. In \cite{rudzkis2013goodness}, a non-parametric method is proposed with a test based on a kernel density estimator. In the last decade, a lot of effort has been put into finding more efficient goodness-of-fit tests. The motivation was mainly coming from graphical models where the distributions are known up to a normalization factor that is often computationally intractable. To address this problem, several tests have been proposed based on Reproducing Kernel Hilbert Space (RKHS) embedding. A large span of them use classes of Stein transformed RKHS functions \citep{liu2016kernelized,gorham2020measuring}. For example in
 \cite{chwialkowski2016kernel}, a goodness-of-fit test is proposed for both i.i.d or non i.i.d samples. The test statistic uses the squared Stein discrepancy, which is naturally estimated by a V-statistic. One drawback of such approach is that the theoretical results provided are only asymptotic. This paper is part of a large list of works that proposed a goodness-of-fit test and where the use of U-statistics naturally emerge \citep[cf.][]{liu2016kernelized,FAN199736,butucea2007goodness,fan1999goodness,fernandez2019maximum,fromont2006}. To conduct a non-asymptotic analysis of the goodness-of-fit tests proposed for non i.i.d samples, a concentration result for U-statistics with dependent random variables is much needed.

 \subsection{Goodness-of-fit test for the density of the stationary measure of a Markov chain}
 In this section, we provide a goodness-of-fit test for Markov chains whose stationary distribution has density with respect to the Lebesgue measure~$\lambda_{Leb}$ on~$\mathds R$. Our work is inspired from \cite{fromont2006} where Fromont and Laurent tackled the goodness-of-fit test with i.i.d samples. Conducting a non-asymptotic theoretical study of our test, we are able to identify the classes of alternatives over which our method has a prescribed power.

Let~$X_1,\dots ,X_n$ be a Markov chain with stationary distribution~$\pi$ with density~$f$ with respect to the Lebesgue measure on~$\mathds R$. Let~$f_0$ be some given density in~$L^2(\lambda_{Leb})$ and let~$\alpha$ be in~$]0,1[$. Assuming that~$f$ belongs to~$L^2(\lambda_{Leb})$, we construct a level~$\alpha$ test of the null hypothesis~$"f=f_0"$ against the alternative~$"f \neq f_0"$ from the observation~$(X_1,\dots ,X_n)$. The test is based on the estimation of~$\|f-f_0\|^2_2$ that is~$\|f\|^2_2+\|f_0\|_2^2-2\langle f, f_0\rangle$.~$\langle f, f_0 \rangle$ is usually estimated by the empirical estimator~$\sum_{i=1}^nf_0(X_i)/n$ and the cornerstone of our approach is to find a way to estimate~$\|f\|_2^2$. We follow the work of \cite{fromont2006} and we introduce a set~$\{S_m,m\in \mathcal M \}$ of linear subspaces of~$L^2(\lambda_{Leb})$. For all~$m$ in~$\mathcal M$, let~$\{p_l,l\in \mathcal L_m\}$ be some orthonormal basis of~$S_m$. The variable
\[\widehat \theta_m = \frac{1}{n(n-1)} \sum_{l \in \mathcal L_m} \sum_{i\neq j =1}^n p_l(X_i)p_l(X_j)\]
estimates~$\|\Pi_{S_m}(f)\|^2_2$ where~$\Pi_{S_m}$ denotes the orthogonal projection onto~$S_m$. Then~$\|f-f_0\|^2_2$ can be approximated by \[\widehat T_m=\widehat \theta_m+\|f_0\|_2^2-\frac{2}{n} \sum_{i=1}^n f_0(X_i),\]
for any~$m$ in~$\mathcal M$. Denoting by~$t_m(u)$ the~$(1-u)$ quantile of the law of~$\widehat T_m$ under the hypothesis~$"f=f_0"$ and considering

\[u_{\alpha}=\sup_{u\in]0,1[} \mathds P_{f_0}\left(\sup_{m\in \mathcal M}(\widehat T_m-t_m(u))>0\right)\leq \alpha,\]
we introduce the test statistic~$T_{\alpha}$ defined by \begin{equation}T_{\alpha}=\sup_{m\in \mathcal M}(\widehat T_m-t_m(u_{\alpha})). \label{test-statistic}\end{equation}

The test consists in rejecting the null hypothesis if~$T_{\alpha}$ is positive. This approach can be read as a multiple testing procedure. Indeed, for each~$m$ in~$\mathcal M$, we construct a level~$u_{\alpha}$ test of the null hypothesis~$"f=f_0"$ by rejecting this hypothesis if~$\widehat T_m$ is larger than its~$(1-u_{\alpha})$ quantile under the hypothesis~$"f=f_0"$. We thus obtain a collection of tests and we decide to reject the null hypothesis if for some of the tests of the collection this hypothesis is rejected.

Now we define the different collection of linear subspaces~$\{S_m,m\in \mathcal M\}$ that we will use in the following. We will focus on constant piecewise functions, scaling functions and, in the case of compactly supported densities, trigonometric polynomials. \label{page-subspaces}

\begin{itemize}
\item For all~$D$ in~$\mathds N^*$ and~$k\in \mathds Z$, let \[I_{D,k}=\sqrt{D} \mathds 1_{[k/D,(k+1)/D[}.\]
For all~$D \in \mathds N^*$, we define~$S_{(1,D)}$ as the space generated by the functions~$\{I_{D,k},k\in \mathds Z\}$ and \[\widehat \theta_{(1,D)}=\frac{1}{n(n-1)}\sum_{k \in \mathds Z} \sum_{i\neq j=1}^n  I_{D,k}(X_i) I_{D,k}(X_j).\]
\item Let us consider a pair of compactly supported orthonormal wavelets~$(\phi, \psi)$ such that for all~$J \in \mathds N$,~$\{\phi_{J,k}=2^{J/2}\phi(2^J \cdot-k), k\in \mathds Z\}\cup \{\psi_{j,k}=2^{j/2}\psi(2^j \cdot-k), j\in \mathds N,j\geq J,k\in \mathds Z\}$ is an orthonormal basis of~$L^2(\lambda_{Leb})$. For all~$J \in \mathds N$ and~$D=2^J$, we define~$S_{(2,D)}$ as the space generated by the scaling functions~$\{\phi_{J,k},k\in \mathds Z\}$ and
\[\widehat \theta_{(2,D)}=\frac{1}{n(n-1)}\sum_{k \in \mathds Z} \sum_{i\neq j=1}^n \phi_{J,k}(X_i)\phi_{J,k}(X_j).\]
\item Let us consider the Fourier basis of~$L^2([0,1])$ given by
\begin{align*}
g_0(x)&= \mathds 1_{[0,1]}(x),\\
g_{2p-1}(x)&=\sqrt 2 \cos(2\pi px) \mathds 1_{[0,1]}(x) \quad \forall p \geq 1, \\
g_{2p}(x)&=\sqrt 2 \sin(2 \pi p x)\mathds 1_{[0,1]}(x) \quad \forall p \geq 1.
\end{align*}

For all~$D \in \mathds N^*$, we define~$S_{(3,D)}$ as the space generated by the functions~$\{g_l,l=0,\dots,D\}$ and
\[\widehat \theta_{(3,D)}=\frac{1}{n(n-1)}\sum_{l=0}^D \sum_{i\neq j=1}^n  g_l(X_i)g_l(X_j).\]
\end{itemize}

We denote~$\mathds D_1=\mathds D_3=\mathds N^*$ and~$\mathds D_2=\{2^J,J\in \mathds N \}$. For~$l$ in~$\{1,2,3\}$,~$D$ in~$\mathds D_l$,~$\Pi_{S_{(l,D)}}$ denotes the orthogonal projection onto~$S_{(l,D)}$ in~$L^2(\lambda_{Leb})$. For all~$l$ in~$ \{1,2,3\}$, we take~$\mathcal D_l \subset \mathds D_l$ with~$\cup_{l\in \{1,2,3\}} \mathcal D_l \neq \emptyset~$ and~$\mathcal D_3=\emptyset$ if the~$X_i$'s are not included in~$[0,1]$. Let~$\mathcal M=\left\{(l, D), l\in \{1,2,3\},D \in \mathcal D_l\right\}.$

Theorem~\ref{hypo-test-appli3} describes classes of alternatives over which the corresponding test has a prescribed power. We work under the additional Assumption~\ref{assumption0}. We refer to Section~\ref{proof-hypo-test-appli3} for the proof of Theorem~\ref{hypo-test-appli3}.

\begin{assumption}
\label{assumption0} The initial distribution of the Markov chain~$(X_i)_{i\geq1}$, denoted~$\chi$, is absolutely continuous with respect to the stationary measure~$\pi$ and its density, denoted by~$\frac{d\chi}{d\pi}$, has finite~$p$-moment for some~$p\in(1,\infty]$, i.e
\[\infty > \left\|\frac{d\chi}{d\pi}\right\|_{\pi,p} := \left\{
    \begin{array}{ll}
      \left[ \int \left| \frac{d\chi}{d\pi} \right|^p d\pi \right]^{1/p} & \mbox{if } p<\infty, \\
    \mathrm{ess }\sup \; \left| \frac{d\chi}{d\pi}\right|& \mbox{if }p=\infty.
    \end{array}
\right.\]
In the following, denote~$q = \frac{p}{p-1}\in [1, \infty)$ (with~$q=1$ if~$p=+\infty$) which satisfies~$\frac{1}{p}+\frac{1}{q}=1$.
\end{assumption}

\clearpage
\begin{theorem} \label{hypo-test-appli3}
Let~$X_1,\dots ,X_n$ a Markov chain on~$\mathds R$ satisfying the Assumptions~\ref{assumption1},~\ref{assumption2} and~\ref{assumption0} with stationary measure~$\pi$. We assume that~$\pi$ has density~$f$ with respect the Lebesgue measure on~$\mathds R$ and let~$f_0$ be some given density. Let~$T_{\alpha}$ be the test statistic defined by Eq.\eqref{test-statistic}. Assume that~$f_0$ and~$f$ belong to~$L^{\infty}(\mathds R)$ (the space of essentially bounded measurable functions on $\mathds R$) and that there exist~$p_1,p_2 \in (1,+\infty]$ such that \[C_{\chi} :=\left\|\frac{1}{f}\frac{d\chi}{d\lambda_{Leb}}\right\|_{f\lambda_{Leb},p_1}\vee \left\|\frac{1}{f_0}\frac{d\chi}{d\lambda_{Leb}}\right\|_{f_0\lambda_{Leb},p_2} <\infty, \]where we used the notations of Assumption~\ref{assumption0}. We fix some~$\gamma$ in~$]0,1[$. For any~$\epsilon \in ]0,2[$, there exist some positive constants~$C_1,C_2,C_3$ such that, setting for all~$m=(l, D)$ in~$\mathcal M$,
\begin{align*}V_m(\gamma)= C_1&\|f\|_{\infty}\frac{\log(3C_{\chi}/\gamma)  }{\epsilon n}
 +C_2\left(\|f\|_{\infty}  \log(D+1)  + \|f_0\|_{\infty}\right)  \frac{\log(3C_{\chi}/\gamma)  }{n}\\
 &+C_3 \left( \|f\|_{\infty}+1\right)D R\left(n,\log \left\{\frac{3 \beta \log n}{\gamma}\right\}\right),\end{align*}
 with \[R(n,u) =\log n \left\{   \frac{ u}{n}  +\left(\frac{u}{n}\right)^{2}  \right\},\]if~$f$ satisfies 
\begin{equation} \label{condition-thm-appli3} \|f-f_0\|^2_2>(1+\epsilon) \inf_{m\in \mathcal M} \left\{ \|f-\Pi_{S_m}(f )\|^2_2+t_m(u_{\alpha})+V_m(\gamma)\right\},\end{equation}
then \[\mathds P_f\left( T_{\alpha}\leq 0 \right) \leq \gamma.\]

\end{theorem}

In order to make the condition \eqref{condition-thm-appli3} more explicit and to study its sharpness, we define the uniform separation rate which provides for any~$\gamma \in (0,1)$ the smallest distance between the set of null hypotheses and the set of alternatives to ensure that the power of our statistic test with level~$\alpha$ is at least~$1-\gamma$.

\begin{definition}
Given~$\gamma \in ]0,1[$ and a class of functions~$\mathcal B \subset L^2(\lambda_{Leb})$, we define the uniform separation rate~$\rho(\Phi_{\alpha},\mathcal B,\gamma)$ of a level~$\alpha$ test~$\Phi_{\alpha}$ of the null hypothesis~$"f\in \mathcal F"$ over the class~$\mathcal B$ as the smallest number~$\rho$ such that the test guarantees a power at least equal to~$(1-\gamma)$ for all alternatives~$f\in \mathcal B$ at a distance~$\rho$ from~$\mathcal F$.
Stated otherwise, denoting by~$d_2(f,\mathcal F)$ the~$L^2$-distance between~$f$ and~$\mathcal F$ and by~$\mathds P_f$ the distribution of the observation~$(X_1,\dots,X_n)$,
\[\rho(\Phi_{\alpha},\mathcal B,\gamma)=\inf \left\{\rho>0,\forall f \in \mathcal B,d_2(f,\mathcal F)\geq \rho \implies \mathds P_f(\Phi_{\alpha} \text{ rejects})\geq 1- \gamma \right\}.\]
\end{definition}
In the following, we derive on explicit upper bound on the uniform separation rates of the test proposed above over several classes of alternatives. For~$s>0,P>0,M>0$ and~$l\in \{1,2,3\}$, we introduce \[\mathcal B^{(l)}_s(P, M)=\left\{f\in L^2(\lambda_{Leb}) \; |\; \forall D\in \mathcal D_l,\quad \|f-\Pi_{S_{(l,D)}}(f )\|_2^2\leq P^2D^{-2s},\|f\|_{\infty}\leq M\right\}.\]

These sets of functions include some Hölder balls or Besov bodies with smoothness~$s$, as highlighted in \cite[Section 2.3]{fromont2006}.
Corollary~\ref{cor:appli3} gives an upper bound for the uniform separation rate of our testing procedure over the classes~$\mathcal B^{(l)}_s(P, M)$ and is proved in Section~\ref{proof-cor:appli3}.

\begin{Cor} \label{cor:appli3} Let~$T_{\alpha}$ be the test statistic defined by \eqref{test-statistic}. Assume that for~$l \in \{1,2,3\}$,~$\mathcal D_l$ is~$\{2^J,0\leq J\leq \log_2\left( n/ (\log(n) \log \log n)^2  \right)\}$ or~$\emptyset$.
For all~$s>0,\; M>0, P>0$ and~$l\in \{1,2,3\}$ such that~$\mathcal D_l \neq \emptyset$, there exists some positive constant~$C=C(s,\alpha,\gamma,M,\|f_0\|_{\infty })$ such that the uniform separation rate of the test~$\mathds 1_{T_{\alpha}>0}$ over~$\mathcal B^{(l)}_s(P, M)$ satisfies for~$n$ large enough
\[\rho\left(\mathds{1}_{T_{\alpha}>0},\mathcal B^{(l)}_s(P, M), \gamma \right) \leq C' P^{\frac{1}{2s+1}}  \left( \frac{\log(n) \log \log n}{n}  \right)^{\frac{s}{2s+1}}.\]
\end{Cor}
{\bf Remark}
In Corollary~\ref{cor:appli3}, the condition~$n$ \textit{large enough} corresponds to 
\[ \left( \log (n) \frac{\log \log n}{n}\right)^{1/2} \leq P \leq  \frac{n^s}{(\log(n)\log \log n)^{2s+1/2}}.\]
For the problem of testing the null hypothesis~$"f=\mathds 1_{[0,1]}"$ against the alternative~$f=\mathds 1_{[0,1]}+g$ with~$g \neq 0$ and~$g\in B_s(P)$ where~$B_s(P)$ is a class of smooth functions (like some Hölder, Sobolev or Besov ball in~$L^2([0,1])$) with unknown smoothness parameter~$s$, Ingster in \cite{ingster} established in the case where the random variables~$(X_i)_{i\geq 1}$ are i.i.d. that the adaptive minimax rate of testing is of order~$(\sqrt{\log \log n}/n)^{2s/(4s+1)}$. From Corollary~\ref{cor:appli3}, we see that our procedure leads to a rate which is close (at least for sufficiently large smoothness parameter~$s$) to the one derived by Ingster in the i.i.d. framework since the upper bound on the uniform separation rate from Corollary~\ref{cor:appli3} can be read (up to a log factor) as~$\left(  [\log \log n]/n  \right)^{\frac{2s}{4s+2}}$.

\subsection{Simulations}

We propose to test our method on three practical examples.\footnote{The code is available at https://github.com/quentin-duchemin/goodness-of-fit-MC.}
In all our simulations, we use Markov chains of length~$n=100$. We choose different alternatives to test our method and we use i.i.d. samples from these distributions. We chose a level~$\alpha = 5 \%$ for all our experiments. All tests are conducted as follows.

\begin{enumerate}
\item We start by the estimation of the~$(1-u)$ quantiles~$t_m(u)$ of the variables~$\widehat T_m=\widehat \theta_m+\|f_0\|_2^2-\frac{2}{n}\sum_{i=1}^n f_0(X_i)$ under the hypothesis~$"f=f_0"$ for~$u$ varying on a regular grid of~$]0,\alpha[$. We sample~$5,000$ sequences of length~$n=100$ with i.i.d. random variables with distribution~$f_0$. We end up with an estimation~$\widehat t_m(u)$ of~$t_m(u)$ for any~$u$ in the grid and any~$m \in \mathcal M$.
\item Then, we estimate the value of~$u_{\alpha}$. We sample again~$5,000$ sequences of length~$n=100$ with i.i.d. random variables with distribution~$f_0$. We use them to estimate the probabilities~$\mathds P_{f_0}(\sup_{m \in \mathcal M}(\widehat T_m-\widehat t_m(u)) >0)$ for any~$u$ in the grid and we keep the larger value of~$u$ such that the corresponding probability is still larger than~$\alpha.$ The selected value of the grid is called~$u_{\alpha}$. Thanks to the first step, we have the estimates~$\widehat t_m(u_{\alpha})$ of~$t_m(u_{\alpha})$ for any~$m \in \mathcal M$.
\item Finally, we sample~$5,000$ Markov chains with length~$n=100$ with stationary distribution~$f$. For each sequence, we can compute~$\widehat T_m$. Dividing by~$5,000$ the number of sequences for which~$\sup_{m \in \mathcal M}(\widehat T_m-\widehat t_m(u_{\alpha})) >0$, we get an estimation of the power of the test.
\end{enumerate}

\begin{figure}[ht]
\begin{minipage}{0.3\linewidth}
\includegraphics[width=\textwidth]{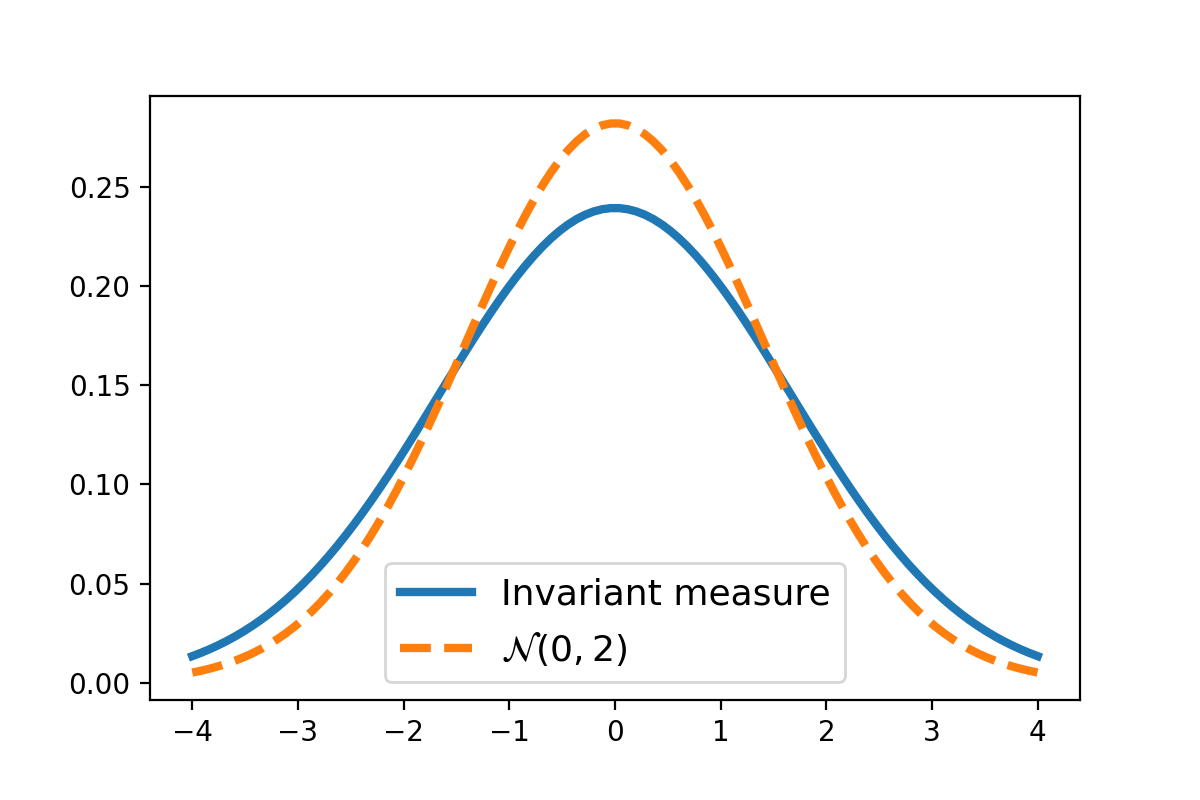}
\subcaption{Example 1}
\label{fig:figure1}
\end{minipage}%
\hfill
\begin{minipage}{0.3\linewidth}
\includegraphics[width=\textwidth]{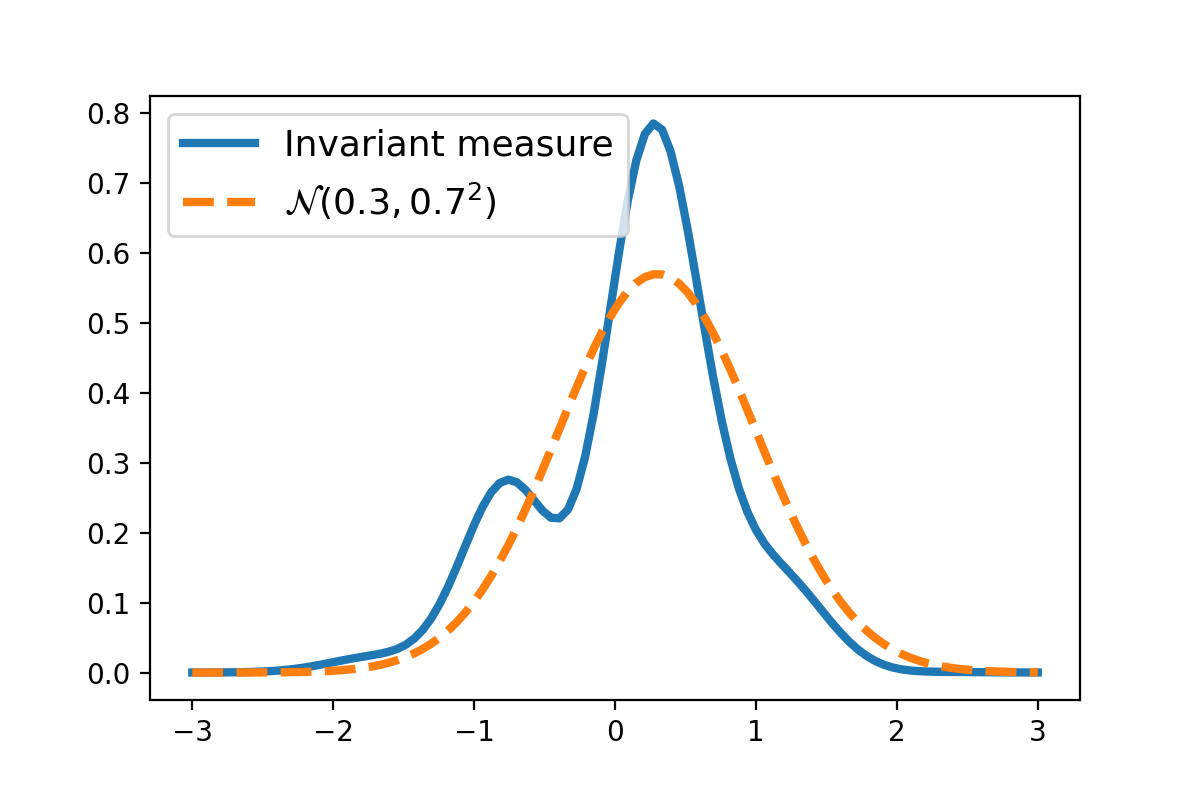}
\subcaption{Example 2}
\label{fig:figure2}
\end{minipage}%
\hfill
\begin{minipage}{0.3\linewidth}
\includegraphics[width=\textwidth]{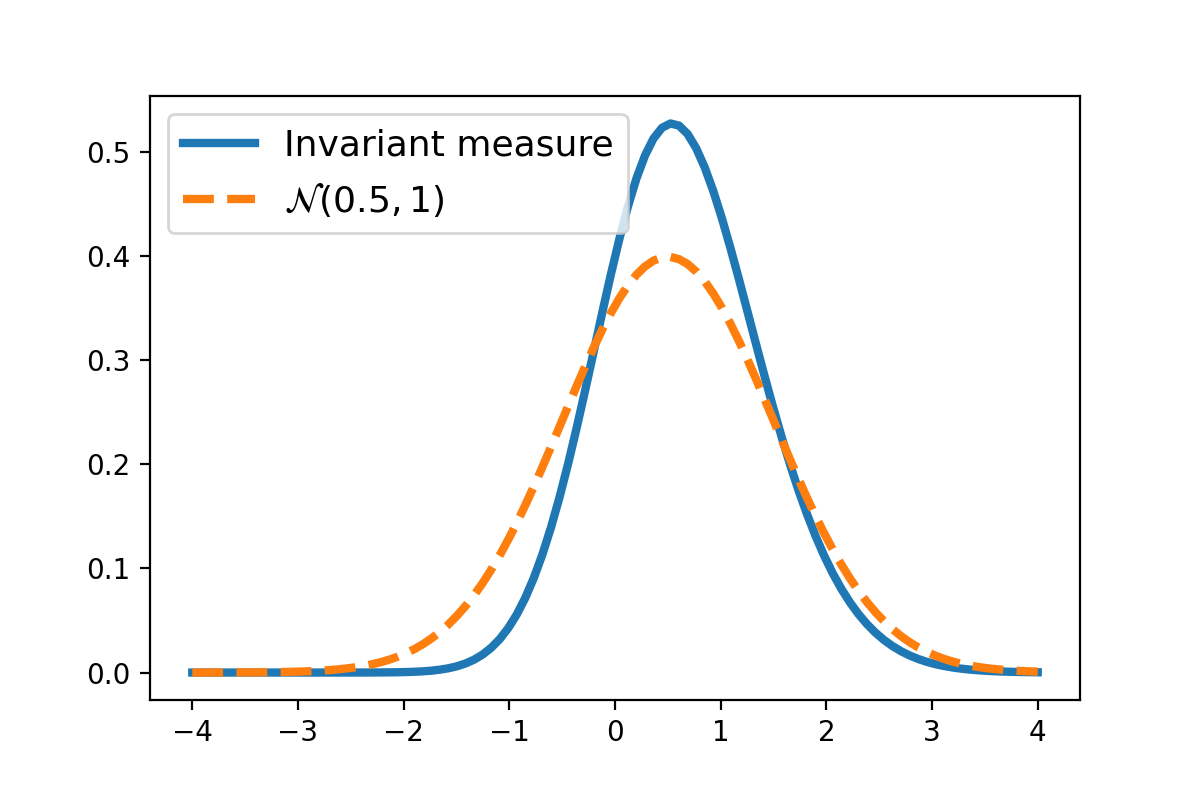}
\subcaption{Example 3}
\label{fig:figure3}
\end{minipage}
\caption{In solid line, we plot the density of the stationary measure of the Markov chain for the three examples of our simulations. In dotted line, we plot the density of the alternative that gives the smaller power on our experiments.}
\label{fig:simu}
\end{figure}

To define comparison points, we compare the power of our test with the classical Kolmogorov-Smirnov test (KS test) and the Chi-squared test ($\chi^2$ test). The rejection region associated with a test of level~$5\%$ is set by \textit{sampling under the null} for both the KS test and the~$\chi^2$ test. With Figure~\ref{fig:simu}, we provide a visualization of the density of the stationary distribution of the Markov chain and of the density of the alternative that gives the smaller power on our experiments.

\subsubsection{Example 1: AR(1) process}
\label{sec:AR1}

Let us consider some~$\theta \in (0,1)$. Then, we define the AR(1) process~$(X_i)_{i\geq 1}$ starting from~$X_1=0$ with for any~$n \geq 1$,
\[X_{n+1} = \theta X_n + \xi_{n+1},\]
where~$(\xi_n)_n$ are i.i.d. random variables with distribution~$\mathcal N(0,\tau^2)$ with~$\tau >0$. From Example 1 of~\cite[Section 2.6]{duchemin}, we know that Assumptions~\ref{assumption1} and~\ref{assumption2} hold. The stationary measure~$\pi$ of the Markov chain~$(X_i)_{i \geq 1}$ is~$\mathcal N \left(0, \frac{\tau^2}{1-\theta^2}\right)$, i.e.~$\pi$ has density~$f$ with respect to the Lebesgue measure on~$\mathds R$ with
\[\forall y \in \mathds R, \quad f(y)=\frac{\sqrt{1-\theta^2}}{\sqrt{2 \pi  \tau^2}} \exp \left( - \frac{(1-\theta^2) y^2}{2\tau^2}  \right).\]

We focus on the following alternatives
\[f_{\mu,\sigma^2}(x) = \frac{1}{\sqrt{2 \pi \sigma^2}} \exp\left( - \frac{(x-\mu)^2}{2 \sigma^2} \right).\]

\noindent
Table~\ref{simu-1} shows the estimated powers for our test, the KS test and the~$\chi^2$ test.

\renewcommand*{\arraystretch}{1.2}
\begin{table}[!ht]
\centering
{\small 
\begin{tabular}{|c||ccc||c|} \hline
$\mathbf{(\mu,\sigma^2)}$ &  \textbf{Our test} &~$\mathbf{\chi^2}$ test&~$\mathbf{KS}$ test & $\mathbf{\|f-f_{\mu,\sigma^2}\|_2}$\\ \hline \hline
$(2,1.5)$        &      0.99     &   0.85  & 0.98 &  \cellcolor[gray]{0.5} 0.39   \\ 
$(0,1)$      &     0.97      &      0.9    &0.8  &  \cellcolor[gray]{0.6} 0.2  \\
$(-0.2,1.2)$        &      0.86     &      0.63  & 0.84   &  \cellcolor[gray]{0.7} 0.17  \\
$(0,1.2)$        &     0.81     &      0.64  & 0.82 &   \cellcolor[gray]{0.8} 0.16     \\
$(0,2)$        &      0.1    &      0.03  &  0.29     &  \cellcolor[gray]{0.85} 0.06\\
\hline
\end{tabular}
}
\caption{Estimated powers of the tests for Markov chains with size~$n=100$. We worked with $\tau=1$,~$\theta=0.8$ and~$\mathcal M = \left\{(1,i) \; : \; i \in \{ 1, \dots , 10 \} \right\}.$ Hence, the stationary distribution of the chain is approximately~$\mathcal N(0,2.8)$. For the~$\chi^2$ test, we work on the interval~$[-5,5]$ that we split into~$20$ regular parts.}
\label{simu-1}
\end{table}

\subsubsection{Example 2: Markov chain generated from independent Metropolis Hasting algorithm}

Let us consider the probability measure~$\pi$ with density~$f$ with respect to the Lebesgue measure on~$[-3,3]$ where
\[\forall x \in [-3,3], \quad f(x) = \frac{1}{Z}e^{-x^2} \left( 3 +\sin(5x) + \sin(2x)\right),\]
with~$Z$ a normalization constant such that~$\int_{-3}^3 f(x)dx = 1$. To construct a Markov chain with stationary measure~$\pi$, we use an independent Metropolis-Hasting algorithm with proposal density~$q(x) \propto \exp(- x^2/6).$ Using Proposition~\ref{prop:inde-hasting}, we get that the above built Markov chain~$(X_i)_{i \geq 1}$ satisfies Assumptions~\ref{assumption1} and~\ref{assumption2}. We focus on the following alternatives
\[g_{\mu,\sigma^2}(x) = \frac{1}{Z(\mu,\sigma^2)} \exp\left( - \frac{(x-\mu)^2}{2 \sigma^2} \right) \mathds 1_{[-3,3]}(x),\]
where~$Z(\mu,\sigma^2)$ is a normalization constant such that~$\int g_{\mu,\sigma^2}(x) dx =1.$ Table~\ref{simu-2} shows the estimated powers for our test, the KS test and the~$\chi^2$ test.

\renewcommand*{\arraystretch}{1.2}
\begin{table}[!ht]
\centering
{\small 
\begin{tabular}{|c||ccc||c|} \hline
$\mathbf{(\mu,\sigma^2)}$ &  \textbf{Our test} &~$\mathbf{\chi^2}$ test&~$\mathbf{KS}$ test & $\mathbf{\|f-g_{\mu,\sigma^2}\|_2}$\\ \hline \hline
$(0,1)$      &    0.96      &      0.91    &0.9 &  \cellcolor[gray]{0.5}0.29   \\
$(0,0.7^2)$        &     0.95    &      0.84  &  0.93  &  \cellcolor[gray]{0.65} 0.23   \\
$(0.3,0.7^2)$        &    0.92     &      0.87  & 0.93  & \cellcolor[gray]{0.8} 0.19   \\\hline
\end{tabular}
}
\caption{Estimated powers of the tests for Markov chains with size~$n=100$. We used~$\mathcal M = \left\{(1,i) \; : \; i \in \{ 1, \dots , 10 \} \right\}.$ For the~$\chi^2$ test, we work on the interval~$[-3,3]$ that we split into~$20$ regular parts.}
\label{simu-2}
\end{table}

\subsubsection{Example 3: ARCH process }

Let us consider some~$\theta \in (-1,1)$. We are interested in the simple threshold auto-regressive model~$(X_n)_{n \geq 1}$ defined by~$X_1=0$ and for any~$n \geq 1$,
\[X_{n+1} = \theta |X_n| + (1-\theta^2)^{1/2} \xi_{n+1},\]
where the random variables~$(\xi_n)_{n \geq 2}$ are i.i.d. with standard Gaussian distribution. From Example 3 of~\cite[Section 2.6]{duchemin}, we know that Assumptions~\ref{assumption1} and~\ref{assumption2} hold. The transition kernel of the Markov chain~$(X_i)_{i \geq 1}$ is 
\[\forall x,y \in \mathds R, \quad P(x,y) = \frac{1}{\sqrt{2 \pi}} \exp\left( - \frac{(y-\theta |x|)^2}{2(1-\theta^2)} \right).\]
The stationary distribution~$\pi$ of the Markov chain has density~$f$ with respect to the Lebesgue measure on~$\mathds R$ with
\[\forall y \in \mathds R, \quad f(y) = \frac{1}{\sqrt{2\pi}} \exp\left(-\frac{y^2}{2}  \right) \Phi\left( \frac{\theta y}{(1-\theta^2)^{1/2}}\right),\]where~$\Phi$ is the standard normal cumulative distribution function. 
We focus on the following alternatives
\[f_{\mu,\sigma^2}(x) = \frac{1}{\sqrt{2 \pi \sigma^2}} \exp\left( - \frac{(x-\mu)^2}{2 \sigma^2} \right).\]
Table~\ref{simu-3} shows the estimated powers for our test, the KS test and the~$\chi^2$ test.

\renewcommand*{\arraystretch}{1.2}
\begin{table}[!ht]
\centering
{\small 
\begin{tabular}{|c||ccc||c|} \hline
$\mathbf{(\mu,\sigma^2)}$ &  \textbf{Our test} &~$\mathbf{\chi^2}$ test&~$\mathbf{KS}$ test & $\mathbf{\|f-f_{\mu,\sigma^2}\|_2}$\\ \hline \hline
$(0,1)$      &   0.98    &      0.85   &0.95 &  \cellcolor[gray]{0.6}0.3  \\
$(1,0.8^2)$      &    0.95     &      0.79    &0.88 &  \cellcolor[gray]{0.7} 0.22  \\
$(0.5,1)$        &      0.3    &      0.07  &   0.5    &  \cellcolor[gray]{0.8} 0.14  \\
$(0.6,0.8^2)$        &      0.35     &      0.16  & 0.4 &  \cellcolor[gray]{0.85} 0.036    \\\hline
\end{tabular}
}
\caption{Estimated powers of the tests for Markov chains with size~$n=100$. We used~$\theta=0.8$ and~$\mathcal M = \left\{(1,i) \; : \; i \in \{ 1, \dots , 10 \} \right\}.$  For the~$\chi^2$ test, we work on the interval~$[-20,20]$ that we split into~$20$ regular parts.}
\label{simu-3}
\end{table}

\subsubsection{Comments on our numerical experiments}

Our experiments show that the~$\chi^2$ goodness-of-fit test give in general the smaller power compared to our method and to the KS test. The~$\chi^2$ test is better suited to deal with discrete probability distributions and it seems to suffer to small power in our continuous setting. Note that using the~$\chi^2$ test with continuous densities on~$\mathds R$ require to specify some hyperparameters (such as a compact interval and the number of bins to discretize it). In practice, the test results can change drastically by modifying these hyperparameters, making the test unreliable. Our experiments also show that when the~$L^2$ norm between the true density $f$ and the alternative one~$f_0$ is large enough, our method reaches higher power compared to the two other procedures considered. Nevertheless, our approach seems less powerful compared to the KS test when the~$L^2$ norm~$\|f-f_0\|_2$ is getting smaller. This is not surprising since our testing procedure is based on the~$L^2$ norm while the KS test relies on the sup norm between cumulative distribution functions (CDFs). We conduct a final experiment to better stress this distinction between our procedure and the KS test. We consider the notations and the framework of the example from Section~\ref{sec:AR1} with the following alternatives
\[f^{(L,\delta)}(x)=\left\{
    \begin{array}{lll}
        f_{0,\frac{\tau^2}{1-\theta^2}}(x) & \mbox{if } |x|\geq \delta \\
        f_{0,\frac{\tau^2}{1-\theta^2}}(x) -L & \mbox{if }  -\delta <x\leq 0 \\
        f_{0,\frac{\tau^2}{1-\theta^2}}(x) +L & \mbox{if }  0 <x<\delta
    \end{array}
\right.,\]
where~$L,\delta>0$ are chosen so that~$f^{(L,\delta)}(x)\geq0$ for any~$x\in \mathds R$. We work with~$\tau=1,\theta=0.8$ and~$\mathcal M = \left\{(1,i) \; : \; i \in \{ 1, \dots , 10 \} \right\}.$ Figure~\ref{fig:Ld} shows the alternatives considered. The sup norm between the CDFs of~$f$ and~$f^{(L,\delta)}$ is equal to~$L\delta$ while the squared~$L^2$ norm between~$f$ and~$f^{(L,\delta)}$ is~$2L^2\delta^3/3$. Hence, we expect that powers will increase for both tests when~$L$ and/or~$\delta$ are increasing. Moreover, we expect the power of our method to be more sensitive to the parameters~$L$ and~$\delta$. Those intuitions are confirmed with the numerical experiments presented in Table~\ref{tab:Ld}.

  \begin{minipage}[b]{0.43\textwidth}
    \centering
    \includegraphics[scale=0.43]{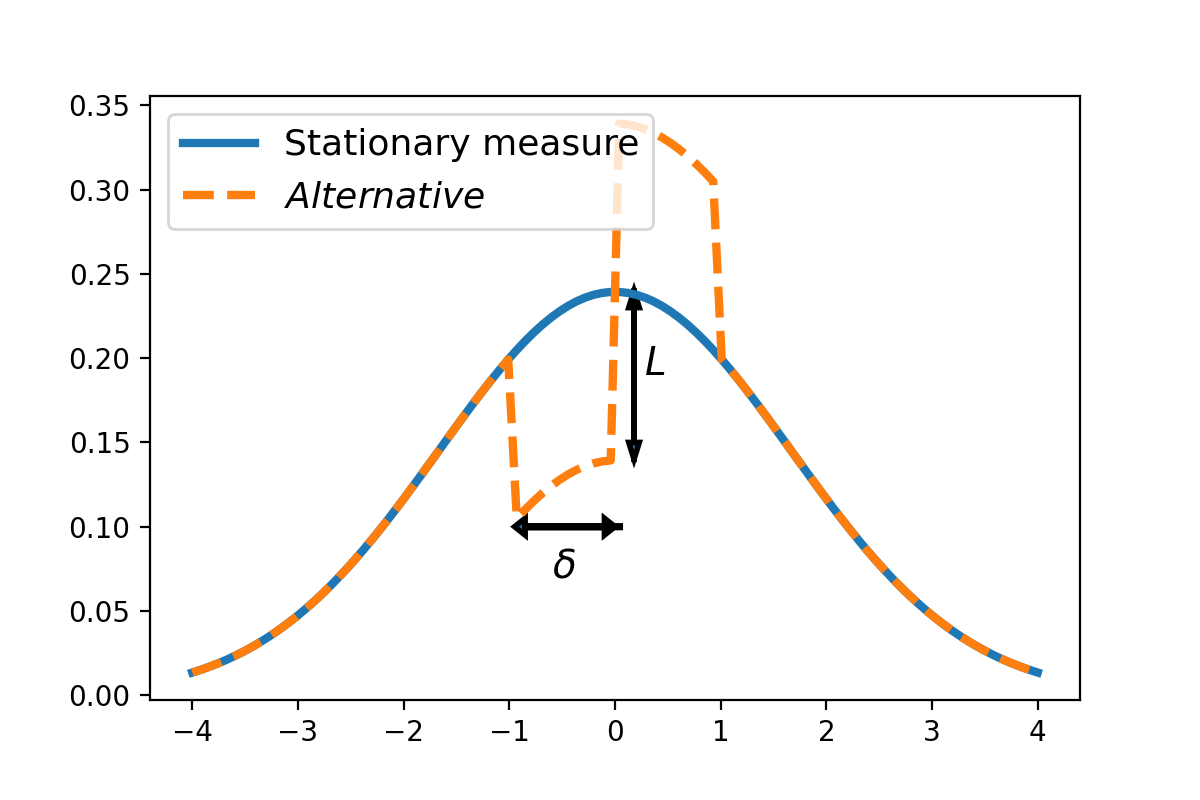}
    \captionof{figure}{Alternative considered.}
    \label{fig:Ld}
  \end{minipage}
  \begin{minipage}[b]{0.5\textwidth}
    \centering
    {\small 
\begin{tabular}{|c||cccc|} \hline
$\mathbf{(L,\delta)}$ &  0.25 & 0.5 & 0.75 & 1\\ \hline \hline
\multirow{2}{*}{0.05}     &  \cellcolor[gray]{0.7} 0.06    &   \cellcolor[gray]{0.7}   0.12   &\cellcolor[gray]{0.7}0.2 &  \cellcolor[gray]{0.7}0.21  \\
 &   0.1    &      0.15   &0.2 &  0.22  \\\hline
\multirow{2}{*}{0.05}     &  \cellcolor[gray]{0.7} 0.16    &   \cellcolor[gray]{0.7}   0.33   &\cellcolor[gray]{0.7}0.36 &  \cellcolor[gray]{0.7}0.4  \\
 &   0.23    &      0.26   &0.33 &  0.37  \\\hline
\multirow{2}{*}{0.1}     &  \cellcolor[gray]{0.7} 0.33    &    \cellcolor[gray]{0.7}  0.66   &\cellcolor[gray]{0.7}0.8 &  \cellcolor[gray]{0.7}0.83  \\
 &   0.26    &      0.35   &0.46 & 0.48  \\\hline
  \multirow{2}{*}{0.15}     &   \cellcolor[gray]{0.7}0.82    &      \cellcolor[gray]{0.7}0.87   &\cellcolor[gray]{0.7}0.9 &  \cellcolor[gray]{0.7}0.95  \\
 &   0.35    &      0.45   &0.55 &  0.67  \\\hline
 \multirow{2}{*}{0.2}     &   \cellcolor[gray]{0.7}0.9    &    \cellcolor[gray]{0.7}  0.93   &\cellcolor[gray]{0.7}0.95 &  \cellcolor[gray]{0.7}0.98  \\
 &   0.46    &      0.54   &0.72 &  0.87  \\\hline
\end{tabular}
}
      \captionof{table}{{\small Estimated powers of the tests for Markov chains with size~$n=100$. Gray cells are the powers of our method while blank cells are the ones obtained with the KS test.}}
      \label{tab:Ld}
    \end{minipage}

\vspace{1cm}

\textbf{Acknowledgments.} We would like to acknowledge support for this project from the Région Ile-de-France.

\newpage

\begin{center}
{\Huge Supplementary Material}
\end{center}

\appendix

{\bf Guidelines for the supplementary material.}

\begin{itemize}
\item \underline{Sections~~\ref{sec:proof-appli1},~\ref{sec:proof-appli2} and~\ref{sec:proof-appli3}: Proofs}\\ Sections~~\ref{sec:proof-appli1},~\ref{sec:proof-appli2} and~\ref{sec:proof-appli3} provide respectively the proofs of our main results from Sections~\ref{sec:integral-ope}, \ref{sec:online-pairwise} and~\ref{sec:adaptivegof}.
\item \underline{Section~\ref{sec:add-proofs}: Technical Lemmas}\\
This section contains some Lemmas useful for our proofs.
\end{itemize}

\setlength\parindent{0pt}

\section{Proofs for Section~\ref{sec:integral-ope}}
\label{sec:proof-appli1}

Let us explain in a nutshell the structure of our proof. For any natural number~$R$, we denote~$\mathbf H^R$ the integral operator with kernel function~$h_R$ at resolution~$R$, namely
\[h_R(x,y):= \sum_{r\in I, r\leq R} \lambda_r \phi_r(x) \phi_r(y), \quad \mathbf H^R f(x) := \int_{E} h_R(x,y) f(y) d\pi(y).\]
We define~$\widetilde {\mathbf H}^R_n$ and~$\mathbf H^R_n$ analogously by using the kernel~$h_R$ in Eq.\eqref{eq:def-Hn}. Using the triangle inequality, we split the distance~$\delta_2(\lambda(\mathbf H),\lambda(\mathbf H_n))$ into four terms.
\begin{enumerate}
\item $\delta_2(\lambda(\mathbf H),\lambda(\mathbf H^R))$ is a bias term induced by working at resolution $R$.
\item A non-trivial preliminary work allows to prove that $\delta_2(\lambda(\mathbf H^R),\lambda(\widetilde {\mathbf H}^R_n))$ can be written as an empirical process of the Markov chain $(X_i)_{i\geq1}$ whose tail can be controlled by applying concentration inequalities for sums of functions of uniformly ergodic Markov chains (this is where we use the assumption that $\Upsilon$ is finite). We refer to Eq.\eqref{1st-approach}.
\item Since the matrices $\mathbf H^R_n$ and $\widetilde{\mathbf H}^R_n$ only differ at diagonal elements,  $\delta_2(\lambda(\widetilde {\mathbf H}^R_n),\lambda(\mathbf H^R_n))$ can be coarsely bounded by $n^{-1/2}\|h_R\|_{\infty}$ (cf. Eq.\eqref{eq:diago}).
\item Applying the Hoffman-Wielandt inequality, one can notice that $\delta_2(\lambda(\mathbf H^R_n),\lambda(\mathbf H_n))$ can be upper-bounded by a U-statistic of order two of the Markov chain $(X_i)_{i\geq1}$ (cf. Eq.\eqref{eq:hoffman-Ustat}). We control the tail behaviour of this U-statistic by applying Theorem~\ref{mainthm2}.
\end{enumerate}
The proof is then concluded by choosing the resolution level $R$ so that $R^2=\ceil{\sqrt n}$.

\subsection{Deviation inequality for the spectrum of signed integral operators}

As shown in Section~\ref{proof-cor-integral-operator-bounded}, Theorem~\ref{cor-integral-operator-bounded} is a direct consequence of the concentration result provided by Theorem~\ref{integral-operator-bounded}.
\begin{theorem}\label{integral-operator-bounded} We keep notations of Section~\ref{sec:integral-ope}. Assume that~$(X_n)_{n \geq 1}$ is a Markov chain on~$E$ satisfying  Assumptions~\ref{assumption1} and~\ref{assumption2} described in Section~\ref{sec:concentration} with stationary distribution~$\pi$. Let us consider some symmetric kernel~$h:E\times E \to \mathds R$, square integrable with respect to~$\pi \otimes \pi$.  Let us consider some~$R \in \mathds N^*$. We assume that there exist continuous functions~$\phi_r:E \to \mathds R$,~$r\in I$ (where~$I=\mathds N$ or~$I=1, \dots,N$) that form an orthonormal basis of~$L^2(\pi)$ such that it holds pointwise \[h(x,y) = \sum_{r \in I} \lambda_r \phi_r(x) \phi_r(y),\]
with \[ \Lambda_R := \sup_{r\in I, \; r\leq R} |\lambda_r| \text{  and  } \Upsilon_R:=\sup_{r\in I, \; r\leq R} \|\phi_r\|_{\infty}^2.\]
We define~$h_R(x,y) = \sum_{r \in I, \; r \leq R} \lambda_r \phi_r(x) \phi_r(y)$ and we assume that~$\|h_R\|_{\infty},\|h-h_R\|_{\infty}<\infty$.
 Then there exists a universal constant $K>0$ such that for any~$t>0$, it holds
\begin{align*}&\mathds P \left(\frac14 \delta_2( \lambda(\mathbf{H}),\lambda( \mathbf{ H}_n))^2   \geq \left( \|h_R\|^2_{\infty} + \kappa \|h-h_R\|_{\infty}^2 \right) \frac{\log n}{n}+2\sum_{i>R, i\in I} \lambda_i^2+t \right)\\
\leq \quad &  16\exp\left(- n \frac{t^2}{Km^2\tau^2\|h-h_R\|_{\infty}^2} \right) +\beta \log(n) \exp\left(  - \frac{n}{16\log n} \left\{\left[\frac{t}{c}\right]  \wedge\left[ \frac{t}{c}\right]^{1/2} \right\} \right)\\
+ \quad & 16R^2\exp\left(   -  \frac{nt}{Km^2\tau^2R^2\Lambda_R^2 \Upsilon_R^2}\right).
\end{align*}
where~$c = \kappa \|h-h_R\|_{\infty}$ with~$\kappa>0$ depending on~$ \delta_M$,~$\tau, L, m$ and~$\rho$.~$\beta$ depends only on~$\rho$.
\end{theorem}

\underline{Proof of Theorem~\ref{integral-operator-bounded}.}
For any integer~$R \geq 1$, we denote
\begin{align*}
X_{n,R} &:= \frac{1}{\sqrt{n}}\left( \phi_r(X_i) \right)_{ 1\leq i \leq n, \; 1\leq r \leq R} \in \mathds R^{n\times R}\\
A_{n,R} &:= \left(X_{n,R}^{\top} X_{n,R}\right)^{1/2} \in \mathds R^{R \times R}\\
K_R &:= \mathrm{Diag}(\lambda_1, \dots , \lambda _R)\\
\mathbf{\widetilde{H}}_{n}^R &:= X_{n,R} K_R X_{n,R}^{\top}\\
\mathbf{H}_{n}^R &:= \left( (1-\delta_{i,j})\left(\mathbf{\widetilde{H}}_{n}^R\right)_{i,j}\right)_{1\leq i,j\leq n}.
\end{align*}

We remark that~$A_{n,R}^2 = I_R + E_{R,n}$ where~$\left(E_{R,n}\right)_{r,s} =(1/n)  \sum_{i=1}^n \left(\phi_r(X_i)\phi_s(X_i) - \delta_{r,s}\right)$ for all~$r,s \in [R].$ Denoting~$\lambda(\mathbf H^R)=(\lambda_1, \dots, \lambda_R)$, we have
\begin{align*}
\delta_2(\lambda(\mathbf H) , \lambda(\mathbf{H}_n))^2 & \leq 4\big[ \delta_2(\lambda(\mathbf H ),\lambda(\mathbf H^R))^2 + \delta_2(\lambda(\mathbf H^R) , \lambda(\mathbf{\widetilde{H}}_{n}^R))^2 +\delta_2( \lambda(\mathbf{\widetilde{H}}_{n}^R),\lambda(\mathbf{H}_{n}^R))^2\\
&\quad + \delta_2( \lambda(\mathbf{H}_{n}^R),\lambda(\mathbf{H}_{n}))^2 \big].
\end{align*}

{\bf Bounding~$\delta_2\big(\lambda(\mathbf H^R) , \lambda(\mathbf{\widetilde{H}}_{n}^R)\big)^2$.}

Let us consider some $\epsilon>0.$

Using a singular value decomposition of~$X_{n,R}$, one can show that~$\lambda (X_{n,R}K_R X_{n,R}^{\top} ) = \lambda( A_{n,R} K_R A_{n,R}) $ which leads to
\begin{align*}
\delta_2\left(\lambda(\mathbf H^R) , \lambda(\mathbf{\widetilde{H}}_{n}^R)\right) & = \delta_2\left(\lambda(K_R) , \lambda(X_{n,R} K_R X_{n,R}^{\top})\right)\\
&= \delta_2\left(  \lambda( K_R),\lambda(A_{n,R}K_RA_{n,R})\right)\\
&\leq \|K_R-A_{n,R} K_R A_{n,R}\|_F,
\end{align*}

Using Equation~$(4.8)$ from \cite[page 127]{gine2000}, we get 
\begin{equation} \label{1st-approach} \delta_2\left(\lambda(\mathbf H^R) , \lambda(\mathbf{\widetilde{H}}_{n}^R)\right)^2  \leq 2\|K_R E_{R,n}\|_F^2 = 2 \sum_{1 \leq r,s \leq R} \lambda_s^2\left( \frac{1}{n} \sum_{i=1}^n \phi_r(X_i) \phi_s(X_i) - \delta_{r,s}\right)^2.\end{equation}

Hence,

\begin{align*}
&\mathds P\left(  \delta_2\left(\lambda(\mathbf H^R) , \lambda(\mathbf{\widetilde{H}}_{n}^R)\right)^2 \geq t \right)\\
 \leq \quad &  \sum_{1\leq s,r\leq R} \mathds P\left(  \sqrt 2 |\lambda_s| \left| \frac{1}{n} \sum_{i=1}^n \phi_r(X_i) \phi_s(X_i) - \delta_{r,s}\right|  \geq \sqrt{t}/R \right)\\
 \leq \quad &   \sum_{1\leq s,r\leq R, \lambda_s\neq 0} \mathds P\left(   \left| \frac{1}{n} \sum_{i=1}^n \phi_r(X_i) \phi_s(X_i) - \delta_{r,s}\right|  \geq \sqrt{t}/(\sqrt 2 R|\lambda_s|)\right)\\
 \leq \quad &  \sum_{1\leq s,r\leq R, \lambda_s\neq 0} 16\exp\left(   -\left(Km^2\tau^2\right)^{-1} \frac{nt}{R^2|\lambda_s|^2\Upsilon_R^4}\right)\\
 = \quad & 16R^2  \exp\left(   - \left(Km^2\tau^2\right)^{-1} \frac{nt}{R^2\Lambda_R^2 \Upsilon_R^4}\right),
\end{align*}
where the last inequality follows from Proposition~\ref{prop:bernstein-orlicz} and where $K>0$ is a universal constant.

\medskip

{\bf Bounding~$\delta_2(\lambda(\mathbf{\widetilde  H}_n^R), \lambda(\mathbf{ H}_n^R))^2$.}

\begin{align}\label{eq:diago}
\delta_2(\lambda(\mathbf{\widetilde H}_n^R), \lambda(\mathbf{ H}_n^R))^2&\leq \|\mathbf{\widetilde H}_n^R-\mathbf{ H}_n^R\|_F^2 = \frac{1}{n^2} \left( \sum_{i=1}^n h_R^2(X_i,X_i) \right) \leq \frac{\|h_R\|^2_{\infty}}{n}.
\end{align}

{\bf Bounding~$\delta_2(\lambda(\mathbf{ H}_n^R),\lambda( \mathbf{ H}_n))^2$.}

\begin{align}\label{eq:hoffman-Ustat}
\delta_2(\lambda(\mathbf{ H}_n^R),\lambda( \mathbf{ H}_n))^2&\leq \|\mathbf{\widetilde H}_n^R-\mathbf{\widetilde H}_n\|_F^2 = \frac{1}{n^2} \left( \sum_{1\leq i,j \leq n, \; i \neq j} (h-h_R)(X_i,X_j)^2 \right).
\end{align}
Let us consider,
\[\forall x,y\in E, \quad m_R(x,y) := (h-h_R)^2(x,y)- s_R(x) - s_R(y)   - \mathds E_{\pi \otimes \pi}[(h-h_R)^2(X,Y)],\]
where~$s_R(x)=  \mathds E_{\pi}[(h-h_R)^2(x,X)]- \mathds E_{\pi \otimes \pi}[(h-h_R)^2(X,Y)]$. One can check that for any~$x \in E$, $\mathds E_{\pi}[m_R(x,X)]=\mathds E_{\pi}[m_R(X,x)]=0$. Hence,~$m_R$ is~$\pi$-canonical. 

\begin{align}&\frac{1}{n(n-1)} \left( \sum_{1\leq i,j \leq n, \; i \neq j} (h-h_R)(X_i,X_j)^2 \right)\\
= \quad &  \frac{1}{n(n-1)}  \sum_{1\leq i,j \leq n, \; i \neq j} m_R(X_i,X_j)+\frac{2}{n}\sum_{i=1}^n s_R(X_i) + \mathds E_{\pi\otimes \pi}[(h-h_R)^2(X,Y)]. \label{hoeffding-decomposition}
\end{align}

Using Theorem~\ref{mainthm2}, we get that there exist two constants~$\beta, \kappa >0$ such that for any~$u\geq1$, it holds with probability at least~$1-\beta e^{-u} \log(n),$
\begin{equation*}\frac{1}{n(n-1)}  \sum_{1\leq i,j \leq n, \; i \neq j} m_R(X_i,X_j) \leq \kappa \|h-h_R\|_{\infty}\log n \left\{   \frac{u}{n} \vee \left(\frac{u}{n}\right)^{2}  \right\}. \end{equation*}

Let us now consider some~$t> 0$ such that \begin{equation} \label{thm2-spectrum}\kappa \|h-h_R\|_{\infty}\log n \left\{   \frac{u}{n}  \vee \left(\frac{u}{n}\right)^{2}  \right\} \leq t.\end{equation} The condition \eqref{thm2-spectrum} is equivalent to
\[u\leq n \left\{\frac{t}{\kappa \|h-h_R\|_{\infty}\log n} \wedge\left(\frac{t}{\kappa \|h-h_R\|_{\infty}\log n}\right)^{1/2} \right\}, \]
which is satisfied in particular if~$t$ and $u$ are such that 
\[ u = \frac{n}{\log n} \left\{  \left[\frac{t}{c}\right]  \wedge\left[ \frac{t}{c}\right]^{1/2} \right\},\]
where~$c = \kappa \|h-h_R\|_{\infty}$. One can finally notice that for this choice of~$u$, the condition~$u \geq 1$ holds in particular for~$n$ large enough in order to have~$n/\log n \geq \kappa \|h-h_R\|_{\infty} t^{-1}$.

We deduce from this analysis that for any~$t>0$, we have for~$n$ large enough to satisfy~$n/\log n  \geq \kappa \|h-h_R\|_{\infty} t^{-1}$,
\[\mathds P\left( \frac{1}{n(n-1)}  \sum_{1\leq i,j \leq n, \; i \neq j} m_R(X_i,X_j) \geq t \right) \leq  \beta \log(n) \exp\left( - \frac{n}{\log n} \left\{  \left[\frac{t}{c}\right]  \wedge\left[ \frac{t}{c}\right]^{1/2} \right\} \right).\]
Using Proposition~\ref{prop:bernstein-orlicz}, we get that for some universal constant $K>0$,
\[\mathds P\left( \frac{2}{n}  \left|\sum_{i=1}^n s_R(X_i)\right| \geq t \right) \leq 16\exp\left(- n \frac{t^2}{K m^2 \tau^2\|h-h_R\|_{\infty}^2} \right).\]
We deduce that for some universal constant $K>0$ it holds \begin{align*}&\mathds P \left(\frac{1}{n^2} \left( \sum_{1\leq i,j \leq n, \; i \neq j} (h-h_R)(X_i,X_j)^2 \right) - \mathds E_{\pi \otimes \pi}\left[ (h-h_R)^2 \right] \geq t  \right)\\
\leq \quad & 16 \exp\left(- n \frac{t^2}{Km^2 \tau^2\|h-h_R\|_{\infty}^2} \right) +\beta \log(n) \exp\left(  - \frac{n}{4\log n} \left\{ \left[\frac{t}{c}\right]  \wedge\left[ \frac{t}{c}\right]^{1/2} \right\} \right).  
\end{align*}
Since~$\mathds E_{\pi \otimes \pi}\left[ (h-h_R)^2(X,Y) \right] = \sum_{i>R, i\in I} \lambda_i^2$, we deduce that 
\begin{align*}&\mathds P \left(\delta_2(\lambda(\mathbf{  H}_n^R), \lambda(\mathbf{ H}_n))^2 - \sum_{i>R, i\in I} \lambda_i^2 \geq t \right)\\
\leq \quad & 16 \exp\left(- n \frac{t^2}{Km^2\tau^2\|h-h_R\|_{\infty}^2} \right) +\beta \log(n) \exp\left(  - \frac{n}{4\log n} \left\{ \left[\frac{t}{c}\right]  \wedge\left[ \frac{t}{c}\right]^{1/2} \right\} \right).\end{align*}
Hence we proved that for any~$u>0$ such that~$n/\log n  \geq \kappa \|h-h_R\|_{\infty} u^{-1}$,
\begin{align*}&\mathds P \left( \frac14 \delta_2(\lambda(\mathbf{H}), \lambda(\mathbf{ H}_n))^2  \geq  \frac{\|h_R\|^2_{\infty}}{n}+2\sum_{i>R, i\in I} \lambda_i^2+ u \right)\\
\leq \quad &  16 \exp\left(- n \frac{u^2}{Km^2\tau^2\|h-h_R\|_{\infty}^2} \right) +\beta \log(n) \exp\left(  - \frac{n}{16\log n} \left\{ \left[\frac{u}{c}\right]  \wedge\left[ \frac{u}{c}\right]^{1/2} \right\} \right)\\
+ \quad & 16R^2 \exp\left(   -  \frac{nu}{Km^2\tau^2R^2\Lambda_R^2 \Upsilon_R^2}\right).
\end{align*}
Considering~$t>0$ and applying the previous inequality with~$u = t+ \frac{\kappa \|h-h_R\|_{\infty} \log n}{n}$, we get 
\begin{align*}&\mathds P \left( \frac14 \delta_2(\lambda(\mathbf{H}), \lambda(\mathbf{ H}_n))^2  \geq  \left( \|h_R\|^2_{\infty} + \kappa \|h-h_R\|_{\infty}^2 \right)\frac{ \log n}{n}+2\sum_{i>R, i\in I} \lambda_i^2+ t \right)\\
\leq \quad & 16 \exp\left(- n \frac{t^2}{Km^2\tau^2\|h-h_R\|_{\infty}^2} \right) +\beta \log(n) \exp\left(  - \frac{n}{16\log n} \left\{ \left[\frac{t}{c}\right]  \wedge\left[ \frac{t}{c}\right]^{1/2} \right\} \right)\\
+ \quad & 16R^2   \exp\left(   -  \frac{nt}{Km^2\tau^2 R^2\Lambda_R^2 \Upsilon_R^2}\right).
\end{align*}
This concludes the proof of Theorem~\ref{integral-operator-bounded}.

\subsection{Proof of Theorem~\ref{cor-integral-operator-bounded}.} \label{proof-cor-integral-operator-bounded}

We consider any~$R \in \mathds N^*$. We remark that for any~$x,y \in E,$
\begin{align*}
|h_R(x,y)|& =\left| \sum_{r=1}^R \lambda_r \phi_r(x) \phi_r(y) \right|\\
&\leq  \left(\sum_{r=1}^R |\lambda_r| \phi_r(x)^2\right)^{1/2}\times \left(\sum_{r=1}^R |\lambda_r| \phi_r(y)^2\right)^{1/2}\text{ (Using Cauchy-Schwarz inequality)}\\
& \leq  S,
\end{align*}
which proves that~$\|h_R\|_{\infty} \leq S.$ Similar computations lead to~$\|h-h_R\|_{\infty} \leq  S.$

Using Theorem~\ref{integral-operator-bounded} we get for any~$t>0$,
\begin{align*}&\mathds P \left(\frac14 \delta_2(\lambda(\mathbf{H}), \lambda(\mathbf{ H}_n))^2  \geq  \frac{S^2(1+\kappa) \log n}{n}+2\sum_{i>R, i\in I} \lambda_i^2+ t \right)\\
\leq \quad &  16 \exp\left(- n \frac{t^2}{Km^2\tau^2 S^2 } \right) +\beta \log(n) \exp\left(  - \frac{n}{16\log n} \left\{ \left[\frac{t}{\kappa S}\right] \wedge\left[ \frac{t}{\kappa S}\right]^{1/2} \right\} \right)\\
+ \quad & 16R^2   \exp\left(   - \frac{nt}{Km^2\tau^2R^2\Lambda^2 \Upsilon^2}\right),
\end{align*}
where~$\displaystyle \Lambda := \sup_{r\geq 1}|\lambda_r|<\infty.$ Choosing~$R^2=\ceil{\sqrt n}$, we get

\begin{align*}&\mathds P \left(\frac14 \delta_2(\lambda(\mathbf{H}),\lambda( \mathbf{ H}_n))^2  \geq  \frac{S^2(1+\kappa) \log n}{n}+2\sum_{i>\ceil{ n^{1/4}}, i\in I} \lambda_i^2+ t \right)\\
\leq \quad &  32 \sqrt{n} \exp\left(- \mathcal C \min\left( n t^2 , \sqrt n t \right)\right) +\beta \log(n) \exp\left( - \frac{n}{\log n}\min\left( \mathcal B t, \left(\mathcal B t\right)^{1/2}\right) \right),
\end{align*}
where~$\mathcal B = \left( K \kappa S\right)^{-1}$ and~$\mathcal C = K^{-1} \left( m^2\tau^2 (S + \Lambda \Upsilon)  \right)^{-2} $.

\section{Proofs for Section~\ref{sec:online-pairwise}}
\label{sec:proof-appli2}

In this section, for any~$k\geq 0$ we denote~$\mathds E_k$ the conditional expectation with respect to the~$\sigma$-algebra $\sigma(X_1, \dots, X_k)$.

\subsection{Proof of Theorem~\ref{pairwise-thm1}} \label{proof-pairwise-thm1}

By definition of~$\mathcal M ^n$, we want to bound \[\mathds P\left( \frac{1}{n-c_n}\sum_{t=c_n}^{n-1}\mathcal R(h_{t-b_n}) -\frac{1}{n-c_n}\sum_{t=c_n}^{n-1} M_t \geq \epsilon \right),\]
which takes the form
\begin{align}  
&\mathds P\left( \frac{1}{n-c_n}\sum_{t=c_n}^{n-1} \left[\mathcal R(h_{t-b_n}) -\mathds E_{t-b_n} [M_t]\right]-\frac{1}{n-c_n}\sum_{t=c_n}^{n-1} \left[ M_t-\mathds E_{t-b_n} [M_t] \right] \geq \epsilon \right)\nonumber\\
\leq \quad & \mathds P\left( \frac{1}{n-c_n}\sum_{t=c_n}^{n-1} \left[\mathcal R(h_{t-b_n}) -\mathds E_{t-b_n} [M_t]\right] \geq \epsilon/2 \right)+ \mathds P\left( \frac{1}{n-c_n}\sum_{t=c_n}^{n-1} \left[\mathds E_{t-b_n} [M_t] -M_t\right] \geq \epsilon/2 \right). \label{eq:pairwise-thm1-decompo}
\end{align}

\subsubsection{Step 1: Martingale difference}

We first deal with the second term of Eq.\eqref{eq:pairwise-thm1-decompo}. Note that we can write
\[\sum_{t=c_n}^{n-1} \left[\mathds E_{t-b_n} [M_t] -M_t\right] = \sum_{t=c_n}^{n-1} \sum_{k=1}^{b_n}\left[\mathds E_{t-k} [M_t] - \mathds E_{t-k+1}[M_t]\right]= \sum_{k=1}^{b_n} \sum_{t=c_n}^{n-1}\left[\mathds E_{t-k} [M_t] - \mathds E_{t-k+1}[M_t]\right].\]
Let us consider some~$k \in \{1,\dots,b_n\}$, then we have that~$V_t^{(k)}=  (\mathds E_{t-k}[M_t]-\mathds E_{t-k+1}[M_t])/(n-c_n)$ is a martingale difference sequence, i.e.~$\mathds E_{t-k}[V_t^{(k)}] = 0$. Since the loss function is bounded in~$[0,1]$, we  have~$|V_t^{(k)}|\leq 2/(n-c_n)$, ~$t=  1,\dots, n$. Therefore by the Hoeffding-Azuma inequality,~$\sum_t V_t^{(k)}$ can be bounded such that
\[\mathds P \left( \frac{1}{n-c_n}\sum_{t=c_n}^{n-1}\left[ \mathds E_{t-k}[M_t]-\mathds E_{t-k+1}[M_t]  \right] \geq \frac{\epsilon}{2b_n}  \right) \leq \exp \left( - \frac{(1-c)n \epsilon^2}{8b_n^2}\right).\]
We deduce that 

\begin{equation}\label{pairwise-azuma}\mathds P\left( \frac{1}{n-c_n}\sum_{t=c_n}^{n-1} \left[\mathds E_{t-b_n} [M_t] -M_t\right] \geq \epsilon/2 \right) \leq b_n \exp \left( - \frac{(1-c)n \epsilon^2}{8b_n^2}\right).\end{equation}

\subsubsection{Step 2: Symmetrization by a ghost sample}

In this step we bound the first term in Eq.\eqref{eq:pairwise-thm1-decompo}. Let us start by introducing a ghost sample~$\{\xi_j\}_{1\leq j \leq n}$, where the random variables~$\xi_j$ i.i.d with distribution~$\pi$. Recall the definition of~$M_t$ and define~$\widetilde M_t$ as
\[M_t=\frac{1}{t-b_n} \sum_{i=1}^{t-b_n}\ell(h_{t-b_n},X_t,X_i), \qquad\widetilde M_t=\frac{1}{t-b_n} \sum_{i=1}^{t-b_n}\ell(h_{t-b_n},X_t,\xi_i).\]

The difference between~$\widetilde M_t$ and~$M_t$ is that~$M_t$ is the sum of the loss incurred by~$h_{t-b_n}$ on the current instance~$X_t$ and all the previous examples~$X_j, \; j= 1,\dots, t-b_n$ on which~$h_{t-b_n}$ is trained, while~$\widetilde M_t$ is the loss incurred by the same hypothesis~$h_{t-b_n}$ on the current instance~$X_t$ and an independent set of examples~$\xi_j, \; j=1,\dots,t-b_n.$

First remark that we have
\begin{align}&\frac{1}{n-c_n} \sum_{t=c_n}^{n-1} \left[\mathcal R(h_{t-b_n})-\mathds E_{t-b_n}[M_t]\right]\notag\\
= \quad & \frac{1}{n-c_n} \sum_{t=c_n}^{n-1} \left[\mathcal R(h_{t-b_n})-\mathds E_{t-b_n}[\widetilde M_t]\right]+ \frac{1}{n-c_n} \sum_{t=c_n}^{n-1} \left[\mathds E_{t-b_n}[\widetilde M_t]-\mathds E_{t-b_n}[M_t]\right].\label{pairwise-step1}
\end{align}
\dnote{The first term of Eq.\eqref{pairwise-step1} is handled in~\cite{online-pairwise} by relying heavily on the assumption that samples are i.i.d \citep[see][Claim 1]{online-pairwise}. Hence, the approach of Wang and al. cannot adapted in our framework. To overcome this difficulty, we use the uniform ergodicity of the Markov chain. This is where the use of the burning parameter $b_n$ is essential.}
\smallskip

Since~$\ell$ is in~$[0,1]$, the first term can be bounded directly using the uniform ergodicity of the Markov chain~$(X_i)_i$ as follows
\begin{align*}& \frac{1}{n-c_n} \sum_{t=c_n}^{n-1} \left[\mathcal R(h_{t-b_n})-\mathds E_{t-b_n}[\widetilde M_t]\right] \\
= \quad & \frac{1}{n-c_n} \sum_{t=c_n}^{n-1} \int_{x \in E} \left(d\pi(x) \mathds E_{X \sim \pi}[\ell(h_{t-b_n},x,X)]-P^{b_n}(X_{t-b_n},dx) \mathds  E_{X \sim \pi}[\ell(h_{t-b_n},x,X)]\right)\\
= \quad & \frac{1}{n-c_n} \sum_{t=c_n}^{n-1} \int_{x \in E}  \mathds E_{X \sim \pi}[\ell(h_{t-b_n},x,X)]\left( d\pi(x)-P^{b_n}(X_{t-b_n},dx)\right)\\
\leq \quad & \frac{1}{n-c_n} \sum_{t=c_n}^{n-1} \int_{x \in E}  \left| d\pi(x)-P^{b_n}(X_{t-b_n},dx)\right|\\
\leq \quad &  L\rho^{b_n},
\end{align*}
where we used Eq.\eqref{uni-ergodicity}.

It remains to control \[\frac{1}{n-c_n} \sum_{t=c_n}^{n-1} \left[\mathds E_{t-b_n}[\widetilde M_t]-\mathds E_{t-b_n}[M_t]\right],\]
and we follow an approach similar to \cite{online-pairwise}. Let us remind that~$M_t$ and~$\widetilde M_t$ depend on the hypothesis~$h_{t-b_n}$ and let us define~$L_t(h_{t-b_n}) = \left[\mathds E_{t-b_n}[\widetilde M_t]-\mathds E_{t-b_n}[M_t]\right]~$. We have
\begin{align}
&\mathds P\left( \frac{1}{n-c_n} \sum_{t=c_n}^{n-1} L_t(h_{t-b_n}) \geq \epsilon \right)\nonumber\\
\leq \quad & \mathds P\left(  \underset{\widehat{h}_{c_n-b_n},\dots,\widehat{h}_{n-1-b_n}}{\sup}  \frac{1}{n-c_n} \sum_{t=c_n}^{n-1} L_t(\widehat h_{t-b_n}) \geq \epsilon \right)\nonumber\\
\leq \quad & \sum_{t=c_n}^{n-1} \mathds P\left( \underset{\widehat{h}\in \mathcal H}{\sup} \; L_t(\widehat h) \geq \epsilon \right).\label{eq:pairwise-eq14}
\end{align}

To bound the right hand side of Eq.\eqref{eq:pairwise-eq14} we give first the following Lemma.
\begin{lemma}\label{pairwise-lemma5} Given any function~$f \in \mathcal H$ and any~$t \geq c_n$,
\[\forall \epsilon >0, \quad \mathds P \left(L_t(f) \geq \epsilon \right)\leq 16\exp\left( -(t-b_n) C(m,\tau) \epsilon^2  \right).\]
\end{lemma}

\dnote{In the i.i.d. framework, the counterpart of Lemma~\ref{pairwise-lemma5} follows from a straightforward application of McDiarmid’s inequality \citep[see][Lemma 5]{online-pairwise}. In our work, we consider uniformly ergodic Markov chains and the proof of Lemma~\ref{pairwise-lemma5} requires extra work. We apply a concentration inequality for Markov chains (see Proposition~\ref{prop:bernstein-orlicz}) which needs to hold for any initial distribution. We apply Proposition~\ref{prop:bernstein-orlicz} by considering the time-reversed sequence and this is where we use the reversibility of the chain.}

\underline{Proof of Lemma~\ref{pairwise-lemma5}.}

Note that 
\begin{align*}
L_t(f) &= \mathds E_{t-b_n}[\widetilde M_t]-\mathds E_{t-b_n}[M_t]\\
&= \frac{1}{t-b_n} \sum_{i=1}^{t-b_n}\left( \mathds E_{t-b_n}[\ell(f,X_t,\xi_i)]-\mathds E_{t-b_n}[\ell(f,X_t,X_i)]\right)\\
&= \frac{1}{t-b_n} \sum_{i=1}^{t-b_n} \mathds E_{\xi \sim \pi}  \left[ \mathds  E_{X_t\sim P^{b_n}(X_{t-b_n},\cdot)}\{\ell(f,X_t,\xi)\} \right]- \mathds  E_{X_t\sim P^{b_n}(X_{t-b_n},\cdot)}\{\ell(f,X_t,X_i)\} .
\end{align*}
Hence, denoting~$m(f,X_{t-b_n},x) = \mathds  E_{X_t\sim P^{b_n}(X_{t-b_n},\cdot)}\{\ell(f,X_t,x)\}~$, we get 
\begin{align*}
L_t(f) &\leq \frac{1}{t-b_n} \sum_{i=1}^{t-b_n}\left\{ \mathds E_{\xi \sim \pi}  \left[ m(f,X_{t-b_n},\xi) \right]- m(f,X_{t-b_n},X_i)\right\} .
\end{align*}

By the reversibility of the chain~$(X_i)_{i\geq1}$, we know that the sequence~$(X_{t-b_n},X_{t-b_n-1}, \dots,X_1)~$ conditionally on~$X_{t-b_n}$ is a Markov chain with stationary distribution~$\pi$. Applying Proposition~\ref{prop:bernstein-orlicz} (see Section~\ref{sec:add-proofs}) we get that
\begin{align*}&\mathds P\left( L_t(f) \geq \epsilon \;|\; X_{t-b_n}  \right)\\
\leq \quad & \mathds P\left( \frac{1}{t-b_n} \sum_{i=1}^{t-b_n}\left\{ \mathds E_{\xi_i \sim \pi}  \left[ m(f,X_{t-b_n},\xi_i) \right]- m(f,X_{t-b_n},X_i)\right\}\geq \epsilon \;|\; X_{t-b_n}  \right)\\ 
\leq \quad & 16 \exp\left( -(t-b_n) C(m,\tau) \epsilon^2  \right),
\end{align*}
for some constant $C(m,\tau)>0$ depending only on $m$ and $\tau$. Then we deduce that 
\begin{align*}
\mathds P\left( L_t(f) \geq \epsilon\right) &= \mathds E\left[ \mathds E\left\{  \mathds 1_{L_t(f) \geq \epsilon} \; |\; X_{t-b_n} \right\} \right]\\
&=  \mathds E\left[ \mathds P\left\{  L_t(f) \geq \epsilon \; |\; X_{t-b_n} \right\} \right]\\
&\leq  16\exp\left( -(t-b_n) C(m,\tau) \epsilon^2  \right),
\end{align*}
which concludes the proof of Lemma~\ref{pairwise-lemma5}.\hfill $\blacksquare$

\medskip

The following two Lemmas are key elements to prove Lemma~\ref{pairwise-lemma8}. Their proofs are strictly analogous to the proofs of Lemmas 6, 7 and 8 from \cite{online-pairwise}.

\begin{lemma}\label{pairwise-lemma6}  \citep[cf.][Lemma 6]{online-pairwise}
For any two functions~$h_1, h_2\in \mathcal H$, the following equation holds \[|L_t(h_1)-L_t(h_2)| \leq 2 \mathrm{Lip}(\phi)\|h_1-h_2\|_{\infty}.\]
\end{lemma}

\begin{lemma}\label{pairwise-lemma7}
Let~$\mathcal H = S_1 \cup \dots \cup S_l$ and~$\epsilon >0$. Then
\[\mathds P \left( \underset{h \in \mathcal H}{\sup} L_t(h) \geq \epsilon \right) \leq \sum_{j=1}^l \mathds P \left( \underset{h \in  S_j}{\sup} L_t(h) \geq \epsilon \right).\]
\end{lemma}

\begin{lemma}\label{pairwise-lemma8} \citep[cf.][Lemma 6]{online-pairwise}
For any~$c_n \leq t \leq n$, it holds
\[\mathds P \left( \underset{h \in \mathcal H}{\sup} L_t(h) \geq \epsilon \right) \leq 16\mathcal N\left(\mathcal H, \frac{\epsilon}{4 \mathrm{Lip}(\phi)}\right) \exp \left(  - \frac{(t-b_n) C(m,\tau)\epsilon^2}{4} \right).\]
\end{lemma}

Combining Lemma~\ref{pairwise-lemma8} and Eq.\eqref{eq:pairwise-eq14}, we have
\[\mathds P\left( \frac{1}{n-c_n} \sum_{t=c_n}^{n-1} L_t(h_{t-b_n}) \geq \epsilon \right) \leq 16 \mathcal N\left(\mathcal H, \frac{\epsilon}{4 \mathrm{Lip}(\phi)}\right)   n \exp \left(  - \frac{(c_n-b_n)  C(m,\tau)\epsilon^2}{4} \right).\]

We deduce that 
\begin{align*}
&\mathds P\left( \frac{1}{n-c_n} \sum_{t=c_n}^{n-1} \left[\mathcal R(h_{t-b_n})-\mathds E_{t-b_n}[M_t]\right]\geq \epsilon/2 \right) \\
\leq \quad &  \mathds P\left( L\rho^{b_n} + \frac{1}{n-c_n} \sum_{t=c_n}^{n-1} \left[\mathds E_{t-b_n}[\widetilde M_t]-\mathds E_{t-b_n}[M_t]\right]\geq \epsilon/2 \right) \\
\leq \quad & 16 \mathcal N\left(\mathcal H, \frac{\epsilon}{8 \mathrm{Lip}(\phi)}\right)  n \exp \left(  - \frac{(c_n-b_n)  C(m,\tau) \left(\epsilon/2-L\rho^{b_n}\right)^2}{4} \right).
\end{align*}

\subsubsection{Step 3: Conclusion of the proof}

\dnote{By considering dependent random variables, we needed the introduction of the burning parameter $b_n$ (see Eq.\eqref{pairwise-step1}). This situation brings extra technicalities to conclude the proof.}
\bigskip

From the previous inequality and \eqref{pairwise-azuma}, we get 

\begin{align*}
&\mathds P\left( \frac{1}{n-c_n}\sum_{t=c_n}^{n-1}\mathcal R(h_{t-b_n}) -\frac{1}{n-c_n}\sum_{t=c_n}^{n-1} M_t \geq \epsilon \right) \\
\leq \quad &  b_n \exp \left( - \frac{(1-c)n \epsilon^2}{8b_n^2}\right)+ 16\mathcal N\left(\mathcal H, \frac{\epsilon}{8 \mathrm{Lip}(\phi)}\right)  n \exp \left(  - \frac{(c_n-b_n)  C(m,\tau) \left(\epsilon/2-L\rho^{b_n}\right)^2}{4} \right).
\end{align*}

Note that~$(c_n-b_n)\epsilon \rho^{b_n}  \underset{n\to \infty}{=} o\left( n \epsilon n^{q \log (\rho)}\right) \underset{n\to \infty}{=} o\left(  n^{1+\xi+q \log (\rho)}\right)$ because by assumption~$\epsilon  \underset{n\to \infty}{=} o\left( n^{\xi}\right)$. However, by choice of~$q$ we have \[1+\xi+q \log (\rho) = 1+\xi+\frac{1+\xi}{\log(1/\rho)} \log (\rho)=0,\]
and we finally get that~$(c_n-b_n)\epsilon \rho^{b_n} \underset{n\to \infty}{=} o\left( 1\right)$. We deduce that for~$n$ large enough it holds
\[\exp \left(  - \frac{(c_n-b_n)  C(m,\tau) \left(\epsilon/2-L\rho^{b_n}\right)^2}{4} \right)\leq 2\exp \left(  - \frac{(c_n-b_n)  C(m,\tau) \epsilon^2}{16} \right).\]

Then, noticing that \[  \exp \left( - \frac{(1-c)n \epsilon^2}{8b_n^2}\right) \underset{n\to \infty}{=} \mathcal O\left(  \exp \left(  - \frac{(c_n-b_n)  C(m,\tau)\epsilon^2}{16b_n^2} \right)  \right),\]
we finally get for~$n$ large enough

\begin{align*}
&\mathds P\left( \frac{1}{n-c_n}\sum_{t=c_n}^{n-1}\mathcal R(h_{t-b_n}) -\frac{1}{n-c_n}\sum_{t=c_n}^{n-1} M_t \geq \epsilon \right) \\
\leq \quad & \left[  32\mathcal N\left(\mathcal H, \frac{\epsilon}{8 \mathrm{Lip}(\phi)}\right) +1 \right]   b_n\exp \left(  - \frac{(c_n-b_n)  C(m,\tau)\epsilon^2}{16b_n^2} \right).
\end{align*}

\subsection{Proof of Theorem~\ref{pairwise-thm1:regret}} \label{proof-pairwise-thm1:regret}

Theorem~\ref{pairwise-thm1} shows that
\begin{align}
&\mathds P\left( \left|\frac{1}{n-c_n}\sum_{t=c_n}^{n-1}\mathcal R(h_{t-b_n}) -\mathcal M^n \right|\geq \epsilon \right) \notag\\
\leq \quad & \left[  32\mathcal N\left(\mathcal H, \frac{\epsilon}{8 \mathrm{Lip}(\phi)}\right) +1 \right]   b_n\exp \left(  - \frac{(c_n-b_n)  C(m,\tau)\epsilon^2}{16b_n^2} \right),\label{eq:proba-regret}
\end{align}
and the assumption on the space $\mathcal H$ gives that for some $\theta>0$, it holds for any $\eta>0$, $\log \mathcal N(\mathcal H, \eta)=\mathcal O(\eta^{-\theta})$. By taking $\epsilon = \frac{\log(n)\log(\log n)}{n^{\frac1{2+\theta}}}$ it is straightforward to prove that the logarithm of the right hand side of Eq.\eqref{eq:proba-regret} goes to~$-\infty$ as~$n\to+\infty$. This concludes the proof of the first part of Theorem~\ref{pairwise-thm1:regret}.

Since the result from Theorem~\ref{pairwise-thm1} trivially holds by considering $h_1= \dots=h_{n-1} = h^*$, the previous computations show that for any~$\delta>0$ there exists some~$N\in \mathds N$ such that for any~$n\geq N$ it holds with probability at least $1-\delta$,
\[ \left| \frac{1}{n-c_n}\sum_{t=c_n}^{n-1}\mathcal R(h_{t-b_n}) -\mathcal M^n \right| \vee \left|  \mathcal M^n(h^*,\dots,h^*)-\mathcal M^n\right| \leq \frac{\log(n)\log(\log n)}{n^{\frac1{2+\theta}}}. \]

Hence, by considering that the online learner has a regret bound~$\mathfrak R_n$ (cf. Definition~\ref{def:regret-bound}), we get that for any~$\delta>0$ there exists some~$N\in \mathds N$ such that for any~$n\geq N$ it holds with probability at least $1-\delta$,
\begin{align*}
&
\frac{1}{n-c_n}\sum_{t=c_n}^{n-1}\mathcal R(h_{t-b_n}) - \mathcal R(h^*)
\\
\leq \quad & 
\frac{1}{n-c_n}\sum_{t=c_n}^{n-1}\mathcal R(h_{t-b_n}) - \mathcal M^n + \mathcal M^n - \mathcal M^n(h^*,\dots,h^*) +  \mathcal M^n(h^*,\dots,h^*)- \mathcal R(h^*)
\\
\leq \quad &2 \frac{\log(n)\log(\log n)}{n^{\frac1{2+\theta}}}+   \mathcal M^n-\inf_{h \in \mathcal H} \mathcal M^n(h,\dots,h)\leq  2 \frac{\log(n)\log(\log n)}{n^{\frac1{2+\theta}}}+ \mathfrak R_n,
\end{align*}
which concludes the proof of Theorem~\ref{pairwise-thm1:regret}.

\subsection{Proof of Theorem~\ref{pairwise-thm9}} \label{proof-pairwise-thm9}

\dnote{The proof of Theorem~\ref{pairwise-thm9} has two main steps. First, we show that~$\mathcal R(\widehat{h})$ is close to $\underset{c_n\leq t \leq n-1}{\min} \mathcal R(h_{t-b_n})+2c_{\gamma}(n-t)$ with high probability. Then we show that $\underset{c_n\leq t \leq n-1}{\min} \mathcal R(h_{t-b_n})+2c_{\gamma}(n-t)$ is close to~$\mathcal M^n$ with high probability. The second step is similar to the proof of~\cite{online-pairwise}. For the first step, we need a concentration inequality for U-statistics of order two for uniformly ergodic Markov chains. This is where we use the Hoeffding decomposition and Theorem~\ref{mainthm2} (see Section~\ref{sec:concentration}).}

Let us recall that for any~$1\leq t \leq n-2$,~$\widehat{\mathcal R}(h_{t-b_n},t+1)= \binom{n-t}{2}^{-1} \sum_{k>i, i \geq t+1}^{n} \ell(h_{t-b_n},X_i,X_k).$ We define 
\[\ell(h,x):= \mathds E_{\pi}[\ell(h,X,x)]-\mathcal R(h), \text{ and }\widetilde \ell(h,x,y)=\ell(h,x,y)-\ell(h,x) - \ell(h,y)-\mathcal R(h).\]

Then for any~$t \in \{b_n+1, \dots, n-2\}$ we have the following decomposition
\begin{align}\widehat{\mathcal R}(h_{t-b_n},t+1)-\mathcal R(h_{t-b_n})=\binom{n-t}{2}^{-1} \sum_{k>i, i \geq t+1}^{n}\widetilde \ell(h_{t-b_n},X_i,X_k)+\frac{2}{n-t}\sum_{i=t+1}^{n} \ell(h_{t-b_n},X_i).\label{eq:hoeffding-decompo} \end{align}
One can check that for any~$x\in E$,~$\mathds E_{\pi}\left[\widetilde \ell(h,X,x)\right]=\mathds E_{\pi}\left[\widetilde \ell(h,x,X)\right]=0$. Moreover, for any hypothesis~$h \in \mathcal H$,~$\|\widetilde \ell(h,\cdot,\cdot)\|_{\infty} \leq 4$ (because the loss function~$\ell$ takes its value in~$[0,1]$). Hence, for any fixed hypothesis~$h \in \mathcal H$, the kernel~$\widetilde \ell(h,\cdot,\cdot)$ satisfies Assumption~\ref{assumption3}. Applying Theorem~\ref{mainthm2}, we know that there exist constants~$\beta,\kappa>0$ such that for any~$t\in \{b_n+1,\dots, n-2\}$ and for any~$\gamma \in (0,1)$, it holds with probability at least~$1-\gamma$,
\[ \left|\binom{n-t}{2}^{-1} \sum_{k>i, i \geq t+1}^{n}\widetilde \ell(h_{t-b_n},X_i,X_k) \right|\leq \kappa \frac{\log (n-t-1)}{n-t-1} \log ( (\beta \vee e^1) \log(n-t+1)/\gamma)^{2}. \]

Note that we used that for ~$u = \log\left( (\beta \vee e^1)\log(n-t+1)/\gamma  \right) \geq 1$ it holds \[\log n \left\{\frac{ u}{n}  \vee  \left( \frac{u}{n}\right)^{2} \right\} \leq  \frac{\log n}{n}u^2.\]

Using Proposition~\ref{prop:bernstein-orlicz}, we also have that for any~$t \in \{b_n+1,\dots, n-2\}$ and any~$\epsilon > 0~$, \[\mathds P\left( \left|\frac{2}{n-t}\sum_{i=t+1}^{n} \ell(h_{t-b_n},X_i)\right| > \epsilon \right)\leq 32 \exp\left(- C(m,\tau)(n-t)\epsilon^2 \right) ,\]
where $C(m,\tau)=(Km^2\tau^2)^{-1}>0$ for some universal constant $K$ (one can check from the proof of Proposition~\ref{prop:bernstein-orlicz} that $K=7\times 10^3$ fits).
We get that for any~$t \in \{b_n+1,\dots, n-2\}$ and any~$\gamma \in (0,1)$, it holds with probability at least~$1-\gamma$,
\[\left|\frac{2}{n-t}\sum_{i=t+1}^{n} \ell(h_{t-b_n},X_i)\right| \leq \frac{\log(32/\gamma)^{1/2} C(m,\tau)^{-1/2}}{\sqrt{n-t}} .\]

We deduce that for any~$t\in \{b_n+1,\dots, n-2\}$ and any fixed~$\gamma \in (0,1)$, it holds with probability at least~$1-\gamma$,
\[ \left|\widehat{\mathcal R}(h_{t-b_n},t+1)-\mathcal R(h_{t-b_n})\right| \leq C(m,\tau)^{-1/2}\sqrt{\frac{\log(64/\gamma)}{n-t}},\]
i.e. \begin{equation}\mathds P\left( \left|\widehat{\mathcal R}(h_{t-b_n},t+1)-\mathcal R(h_{t-b_n})\right| \geq c_{\gamma}(n-t)\right) \leq \frac{\gamma}{(n-c_n)(n-c_n+1)}.\label{eq:pairwise-eq32}\end{equation}

Based on the selection procedure of the hypothesis~$\widehat h$ defined in Eq.\eqref{h-selection}, the concentration result Eq.\eqref{eq:pairwise-eq32} allows us to show that~$\mathcal R(\widehat{h})$ is close to $\underset{c_n\leq t \leq n-1}{\min} \mathcal R(h_{t-b_n})+2c_{\gamma}(n-t)$ with high probability. This is stated by Lemma~\ref{pairwise-lemma25} which is proved in Section~\ref{proof-pairwise-lemma25}. 

\begin{lemma} \label{pairwise-lemma25}
Let~$h_0,\dots, h_{n-1}$ be the set of hypotheses generated by an arbitrary online algorithm~$\mathcal A$ working with a pairwise loss~$\ell$ which satisfies the conditions given in Theorem~\ref{pairwise-thm1}. Then for any~$\gamma \in (0,1)$, we have \[\mathds P\left(\mathcal R(\widehat{h})>\min_{c_n\leq t<n-1}(\mathcal R(h_{t-b_n}) + 2c_{\gamma}(n-t))\right)\leq \gamma.\]
\end{lemma}

To conclude the proof, we need to show that~$\underset{c_n\leq t \leq n-1}{\min} \mathcal R(h_{t-b_n})+2c_{\gamma}(n-t)$ is close to~$\mathcal M^n$. 

First we remark that
\begin{align*}
&\underset{c_n\leq t \leq n-1}{\min} \mathcal R(h_{t-b_n})+2c_{\gamma}(n-t) \\
= \quad & \underset{c_n\leq t \leq n-1}{\min}\quad  \underset{t\leq i\leq n-1}{\min} \mathcal R(h_{i-b_n})+2c_{\gamma}(n-i)\\
\leq \quad & \underset{c_n\leq t \leq n-1}{\min}\; \frac{1}{n-t} \sum_{i=t}^{n-1}\left( \mathcal R(h_{i-b_n})+2c_{\gamma}(n-i)\right)\\ 
\leq \quad & \underset{c_n\leq t \leq n-1}{\min} \left( \frac{1}{n-t} \sum_{i=t}^{n-1} \mathcal R(h_{i-b_n})+ \frac{2}{n-t}\sum_{i=t}^{n-1} \sqrt{\frac{C(m,\tau)^{-1}}{n-i}\log\frac{ 64(n-c_n)(n-c_n+1)}{\gamma}}\right)\\ 
\leq \quad & \underset{c_n\leq t \leq n-1}{\min} \left( \frac{1}{n-t} \sum_{i=t}^{n-1} \mathcal R(h_{i-b_n})+ \frac{2}{n-t}\sum_{i=t}^{n-1} \sqrt{\frac{2C(m,\tau)^{-1}}{n-i}\log\frac{ 64(n-c_n+1)}{\gamma}}\right)\\ 
\leq \quad & \underset{c_n\leq t \leq n-1}{\min} \left( \frac{1}{n-t} \sum_{i=t}^{n-1} \mathcal R(h_{i-b_n})+ 4 \sqrt{\frac{2C(m,\tau)^{-1}}{n-t}\log\frac{ 64(n-c_n+1)}{\gamma}}\right),
\end{align*}
where the last inequality holds because~$\sum_{i=1}^{n-t}\sqrt{1/i} \leq 2 \sqrt{n-t}.$ Indeed,~$x\mapsto 1/\sqrt{x}$ is a decreasing and continuous function and a classical serie/integral approach leads to \[\sum_{i=1}^{n-t}\sqrt{1/i} \leq 1+\int_1^{n-t} \frac{1}{\sqrt{x}}dx = 1+ \left[ 2\sqrt{x}\right]_1^{n-t}\leq 2 \sqrt{n-t}.\]

We define~$\mathcal M_t^n:= \frac{1}{n-t} \sum_{m=t}^{n-1} M_m.$ From Theorem~\ref{pairwise-thm1}, one can see that for each~$t=c_n, \dots, n-1$, 

\[\mathds P\left(  \frac{1}{n-t} \sum_{i=t}^{n-1} \mathcal R(h_{i-b_n}) \geq \mathcal M_t^n +\epsilon \right) \leq \left[  32\mathcal N\left(\mathcal H, \frac{\epsilon}{8 \mathrm{Lip}(\phi)}\right)+1 \right] b_n\exp \left(  - \frac{(t-b_n)  C(m,\tau)\epsilon^2}{16b_n^2} \right).\]

Let us set \[K_t=\mathcal M_t^n+ 4 \sqrt{\frac{2C(m,\tau)^{-1}}{n-t}\log\frac{ 64(n-c_n+1)}{\gamma}}+\epsilon.\]

Using the fact that if~$\min(a_1, a_2)\leq \min(b_1, b_2)$ then either~$a_1\leq b_1$ or~$a_2\leq b_2$, we can write

\begin{align*}
& \mathds P \left( \underset{c_n\leq t \leq n-1}{\min} \mathcal R(h_{t-b_n})+2c_{\gamma}(n-t) \geq \underset{c_n\leq t \leq n-1}{\min} K_t \right)\\
\leq \quad & \mathds P \left( \underset{c_n\leq t \leq n-1}{\min} \left( \frac{1}{n-t} \sum_{i=t}^{n-1} \mathcal R(h_{i-b_n})+ 4 \sqrt{\frac{2C(m,\tau)^{-1}}{n-t}\log\frac{ 64(n-c_n+1)}{\gamma}}\right) \geq  \underset{c_n\leq t \leq n-1}{\min} K_t  \right) \\
\leq \quad & \sum_{t=c_n}^{n-1} \mathds P \left(  \frac{1}{n-t} \sum_{i=t}^{n-1} \mathcal R(h_{i-b_n})+ 4 \sqrt{\frac{2C(m,\tau)^{-1}}{n-t}\log\frac{ 64(n-c_n+1)}{\gamma}} \geq  K_t  \right) \\
= \quad & \sum_{t=c_n}^{n-1} \mathds P \left(  \frac{1}{n-t} \sum_{i=t}^{n-1} \mathcal R(h_{i-b_n})\geq \mathcal M_t^n+ \epsilon \right) \\
\leq \quad & (n-c_n) \left[  32\mathcal N\left(\mathcal H, \frac{\epsilon}{8 Lip(\phi)}\right)+1 \right]b_n\exp \left(  - \frac{(c_n-b_n)  C(m,\tau)\epsilon^2}{16b_n^2} \right)\\
\leq \quad & \left[  32\mathcal N\left(\mathcal H, \frac{\epsilon}{8 Lip(\phi)}\right)+1 \right] \exp \left(  - \frac{(c_n-b_n)  C(m,\tau)\epsilon^2}{16b_n^2}+2\log n \right).
\end{align*}

Using Lemma~\ref{pairwise-lemma25}, we get

\begin{align*}&\mathds P\left( \mathcal R(\widehat{h})\geq \underset{c_n\leq t \leq n-1}{\min} \mathcal M_t^n+ 4 \sqrt{\frac{2C(m,\tau)^{-1}}{n-t}\log\frac{ 64(n-c_n+1)}{\gamma}}+\epsilon\right) \\
\leq \quad & \gamma + \left[  32\mathcal N\left(\mathcal H, \frac{\epsilon}{8 Lip(\phi)}\right)+1 \right] \exp \left(  - \frac{(c_n-b_n)  C(m,\tau)\epsilon^2}{16b_n^2}+2\log n \right),
\end{align*}
which gives in particular
\begin{align*}&\mathds P\left( \mathcal R(\widehat{h})\geq \mathcal M^n+ 4 \sqrt{\frac{2C(m,\tau)^{-1}}{n-c_n}\log\frac{ 64(n-c_n+1)}{\gamma}}+\epsilon\right) \\
\leq \quad & \gamma + \left[  32\mathcal N\left(\mathcal H, \frac{\epsilon}{8 Lip(\phi)}\right)+1 \right]  \exp \left(  - \frac{(c_n-b_n)  C(m,\tau)\epsilon^2}{16b_n^2}+2\log n \right).
\end{align*}

We substitute~$\epsilon$ with~$\epsilon/2$  and  we choose~$\gamma$ such that~$4 \sqrt{\frac{2C(m,\tau)^{-1}}{n-c_n}\log\frac{ 64(n-c_n+1)}{\gamma}}=\epsilon/2$ with~$n$ large enough to ensure that~$\gamma <1$. We have for any~$c >0$,
\begin{align*}&\mathds P\left( \mathcal R(\widehat{h})\geq \mathcal M^n+ 4 \sqrt{\frac{2C(m,\tau)^{-1}}{n-c_n}\log\frac{ 64(n-c_n+1)}{\gamma}}+\frac{\epsilon}{2}\right) \\
\leq \quad &64(n-c_n+ 1) \exp\left(-\frac{(n-c_n)C(m,\tau)\epsilon^2}{128}\right) + \left[  32\mathcal N\left(\mathcal H, \frac{\epsilon}{16 Lip(\phi)}\right)+1 \right]\times \\
&\exp \left(  - \frac{(c_n-b_n)  C(m,\tau)\epsilon^2}{(16b_n)^2}+2\log n \right)\\
\leq \quad &  32\left[  \mathcal N\left(\mathcal H, \frac{\epsilon}{16 Lip(\phi)}\right)+1 \right]  \exp \left(  - \frac{(c_n-b_n)  C(m,\tau)\epsilon^2}{(16b_n)^2}+2\log n \right),
\end{align*}
where these inequalities hold for~$n$ large enough.

\subsection{Proof of Lemma~\ref{pairwise-lemma25}}
\label{proof-pairwise-lemma25}

Let \[T^*:=   \arg \min_{c_n\leq t<n-1}(\mathcal R(h_{t-b_n}) + 2c_{\gamma}(n-t)),\] and~$h^*=h_{T^*-b_n}$ is the corresponding hypothesis that minimizes the penalized true risk and let~$\widehat{\mathcal R}^*=\widehat{\mathcal R}(h^*, T^*+1)$ to be the penalized empirical risk of~$h_{T^*-b_n}$. Set, for brevity \[\widehat{\mathcal R}_{t-b_n}=\widehat{\mathcal R}(h_{t-b_n}, t+ 1),\] and let \[\widehat{T}:= \arg \min_{c_n\leq t<n-1}(\widehat{\mathcal R}_{t-b_n}+c_{\gamma}(n-t)),\] where~$\widehat{h}$ coincides with~$h_{\widehat{T}-b_n}$. Using this notation and since \[ \widehat{\mathcal R}_{\widehat T-b_n} + c_{\gamma}(n-\widehat{T}) \leq \widehat{\mathcal R}^* + c_{\gamma}(n-T^*),\]
holds with certainty, we have
\begin{align*}
&\mathds P\left(\mathcal R(\widehat{h}) > \mathcal R(h^*)+\mathcal{E}  \right) \\
= \quad & \mathds P\left(\mathcal R(\widehat{h}) > \mathcal R(h^*)+\mathcal{E}  , \widehat{\mathcal R}_{\widehat T-b_n} + c_{\gamma}(n-\widehat{T}) \leq \widehat{\mathcal R}^* + c_{\gamma}(n-T^*)\right) \\
\leq \quad & \mathds P\left(  \bigcup_{c_n \leq t \leq n-1} \left\{\mathcal R(h_{t-b_n}) > \mathcal R(h^*)+\mathcal{E}  , \widehat{\mathcal R}_{t-b_n}+ c_{\gamma}(n-t) \leq \widehat{\mathcal R}^* + c_{\gamma}(n-T^*) \right\}\right) \\
\leq \quad &  \sum_{t=c_n}^{n-1} \mathds P\left( \mathcal R(h_{t-b_n}) > \mathcal R(h^*)+\mathcal{E} , \widehat{\mathcal R}_{t-b_n}+ c_{\gamma}(n-t) \leq \widehat{\mathcal R}^* + c_{\gamma}(n-T^*)\right),
\end{align*}
where~$\mathcal E$ is a positive-valued random variable to be specified. Now we remark that if 
\begin{equation}\widehat{\mathcal R}_{t-b_n}+ c_{\gamma}(n-t) \leq \widehat{\mathcal R}^* + c_{\gamma}(n-T^*), \label{eq:error}\end{equation}
holds, then at least one of the following three conditions must hold
\begin{center}
\begin{tabular}{rcl}
$(i)$ &~$\widehat{\mathcal R}_{t-b_n}$ &$\leq \mathcal R(h_{t-b_n})-c_{\gamma}(n-t)$\\
$(ii)$ &$\widehat{\mathcal R}^*$ &$> \mathcal R(h^*)+c_{\gamma}(n-T^*)$\\
$(iii)$ &~$\mathcal R(h_{t-b_n})-\mathcal R(h^*)$ &~$\leq 2c_{\gamma}(n-T^*)$.
\end{tabular}
\end{center}

Stated otherwise, if Eq.\eqref{eq:error} holds for some~$t\in \{c_n, \dots,n-1\}$ then
\begin{itemize}
\item either~$t=T^*$ and~$(iii)$ holds trivially.
\item or~$t\neq T^*$ which can occur because
\begin{itemize}
\item~$\widehat{\mathcal R}_{t-b_n}$ underestimates~$\mathcal{R}(h_{t-b_n})$ and~$(i)$ holds.
\item~$\widehat{\mathcal R}^*$ overestimates~$\mathcal R(h^*)$ and~$(ii)$ holds.
\item~$n$ is too small to statistically distinguish~$\mathcal R(h_{t-b_n})$ and~$\mathcal R(h^*)$, and~$(iii)$ holds.
\end{itemize}
\end{itemize}
Therefore, for any fixed~$t$, we have
\begin{align*}
&\mathds P\left( \mathcal R(h_{t-b_n}) > \mathcal R(h^*)+\mathcal{E} , \widehat{\mathcal R}_{t-b_n}+ c_{\gamma}(n-t) \leq \widehat{\mathcal R}^* + c_{\gamma}(n-T^*)\right)\\
\leq \quad & \mathds P\left( \widehat{\mathcal R}_{t-b_n} \leq \mathcal R(h_{t-b_n})-c_{\gamma}(n-t) \right)+  \mathds P\left( \widehat{\mathcal R}^* > \mathcal R(h^*)+c_{\gamma}(n-T^*) \right)\\
\qquad &+ \mathds P\left( \mathcal{R}(h_{t-b_n})-\mathcal{R}(h^*) \leq 2c_{\gamma}(n-T^*) \;,\; \mathcal R(h_{t-b_n}) > \mathcal R(h^*)+\mathcal{E} \right).
\end{align*}
By choosing~$\mathcal E= 2c_{\gamma}(n-T^*)$, the last term in the previous inequality is zero and we can write
\begin{align*}
&\mathds P\left(  \mathcal R(\widehat h) > \mathcal R(h^*)+2c_{\gamma}(n-T^*)\right)\\
\leq \quad & \sum_{t=c_n}^{n-1} \mathds P\left( \widehat{\mathcal{R}}_{t-b_n} \leq \mathcal{R}(h_{t-b_n}) -c_{\gamma}(n-t)\right)+ (n-c_n) \mathds P\left( \widehat{\mathcal R}^* >\mathcal R(h^*)+c_{\gamma}(n-T^*) \right)\\
\leq \quad & (n-c_n) \frac{\gamma}{(n-c_n)(n-c_n+1)} +(n-c_n)\left\{  \sum_{t=c_n}^{n-1} \mathds P \left(  \widehat{\mathcal R}_{t-b_n} > \mathcal R(h_{t-b_n})+c_{\gamma}(n-t)\right)\right\} \text{ (Using \eqref{eq:pairwise-eq32})}\\
\leq \quad & \frac{\gamma}{n-c_n+1} +(n-c_n)^2 \frac{\gamma}{(n-c_n)(n-c_n+1)} \text{ (Using Eq.\eqref{eq:pairwise-eq32})}\\
\leq \quad & \frac{\gamma}{n-c_n+1} +(n-c_n) \frac{\gamma}{n-c_n+1} = \gamma.
\end{align*}

\subsection{Proof of Corollary~\ref{pairwise-cor9}}

The proof of Corollary~\ref{pairwise-cor9} is analogous to the proof of Theorem~\ref{pairwise-thm1:regret} by applying Theorem~\ref{pairwise-thm9} (instead of Theorem~\ref{pairwise-thm1}) and by choosing $\epsilon = \frac{\log^{2}n}{n^{\frac1{2+\theta}}}$.

\section{Proofs for Section~\ref{sec:adaptivegof}}
\label{sec:proof-appli3}

\subsection{Proof of Theorem~\ref{hypo-test-appli3}}
\label{proof-hypo-test-appli3} 
In the following,~$\mathds P_g$ will denote the distribution of the Markov chain if the stationary distribution of the chain is assumed to have a density~$g$ with respect to the Lebesgue measure on~$\mathds R$. We consider~$q=q_1\vee q_2$ where~$q_1,q_2 \in [1,\infty)$ are such that~$\frac{1}{p_1}+\frac{1}{q_1}=1$ and~$\frac{1}{p_2}+\frac{1}{q_2}=1.$

The main tool of the proof is the Hoeffding (also called canonical) decomposition of the~U-statistics~$\widehat \theta_m$. We introduce the processes~$U_n$ and~$P_n$ defined by \[U_n(h)=\frac{1}{n(n-1)}\sum_{i\neq j =1}^n h(X_i,X_j), \quad P_n(h)=\frac{1}{n}\sum_{i=1}^n h(X_i).\] We also define~$P (h)=\langle h, f\rangle.$ By setting, for all~$m \in \mathcal M$, \[H_m(x, y)=\sum_{l \in \mathcal L_m} (p_l(x)-a_l)(p_l(y)-a_l),\]
with~$a_l=\langle f, p_l\rangle,$ we obtain the decomposition \[\widehat \theta_m=U_n(H_m)+(P_n-P)(2\Pi_{S_m}(f ))+ \|\Pi_{S_m}(f )\|_2^2.\]
Let us consider~$\beta$ in~$]0,1[$. Since \[\mathds P_f(T_{\alpha}\leq 0)=\mathds P_f \left(\sup_{m\in\mathcal M}(\widehat \theta_m+\|f_0\|_2^2-\frac{2}{n}\sum_{i=1}^n f_0(X_i)-t_m(u_{\alpha}))\leq 0\right),\] we have \[\mathds P_f(T_{\alpha}\leq 0)\leq \inf_{m \in \mathcal M} \mathds P_f\left(\widehat \theta_m+\|f_0\|^2_2-\frac{2}{n}\sum_{i=1}^n f_0(X_i)-t_m(u_{\alpha})\leq 0 \right).\]
Since~$\|f-\Pi_{S_m}(f )\|^2_2=\|f\|^2_2-\| \Pi_{S_m}(f )\|^2_2$, it holds 
\begin{align*}&\widehat \theta_m+\|f_0\|^2_2-\frac{2}{n}\sum_{i=1}^n f_0(X_i)\\
=\quad & U_n(H_m)+(P_n-P)(2\Pi_{S_m}(f ))- \|f-\Pi_{S_m}(f )\|^2_2+\|f\|^2_2 +\|f_0\|^2_2-2P_n(f_0)\\
=\quad &  U_n(H_m)+(P_n-P)(2\Pi_{S_m}(f ))- \|f-\Pi_{S_m}(f )\|^2_2+\|f-f_0\|^2_2 +2P(f_0)-2P_n(f_0), 
\end{align*}
which leads to
\begin{align}\mathds P_f(T_{\alpha} \leq 0)& \leq \inf_{m\in \mathcal M} \mathds P_f\Bigg(U_n(H_m)+(P_n-P)(2 \Pi_{S_m}(f )-2f)+(P_n-P)(2f-2f_0)+\|f-f_0\|^2_2 \notag\\
& \qquad  \qquad \leq \|f- \Pi_{S_m}(f ) \|^2_2+t_m(u_{\alpha})\Bigg). \label{sup-2-inf}\end{align} We then need to control~$U_n(H_m),(P_n-P)(2 \Pi_{S_m}(f )-2f)$,~$(P_n-P)(2f-2f_0)$ for every~$m\in \mathcal M$.

{\bf Control of~$U_n(H_m)$.}

$H_m$ is~$\pi$-canonical and a direct application of Theorem~\ref{mainthm2} leads to the following Lemma (the proof of Lemma~\ref{lemma:appli3} is postponed to Section~\ref{proof-lemma:appli3}).

\begin{lemma} \label{lemma:appli3} Let us assume that the stationary distribution of the Markov chain~$(X_i)_{i\geq 1}$ has density~$f$ with respect to the Lebesgue measure on~$\mathds R$.
For all~$m=(l, D)$ with~$l\in \{1,2,3\}$ and~$D\in \mathds D_l$, introduce~$\{p_l,l\in \mathcal L_m\}$ defined as in page~\pageref{page-subspaces} and~$Z_m=\frac{1}{n(n-1)} \sum_{i \neq j=1}^n H_m(X_i,X_j),$ with~$H_m(x, y)=\sum_{l \in \mathcal L_m} (p_l(x)- \langle f, p_l\rangle )(p_l(y)-\langle f, p_l \rangle )$. There exist some constants~$C,\beta>0$ (both depending on the Markov chain~$(X_i)_{i \geq 1}$ while~$C$ also depends on~$\phi$) such that, for all~$l \in \{1,2,3\}$,~$ D\in \mathds D_l$ and~$u \geq 1$, it holds with probability at least~$1-\beta e^{-u}\log n$,
\[ |Z_{(l,D)} |\leq C \left( \|f\|_{\infty}+1\right) D R(n,u),\]
where~$R(n,u) =\log n \left\{   \frac{ u}{n}+ \left(\frac{u}{n}\right)^{2}  \right\}.$
\end{lemma}

We deduce that there exist~$C, \beta>0$ such that for any~$\gamma \in (0,1\wedge (e^{-1}3\beta \log n))$ and any~$m=(l,D) \in \mathcal M$,
\begin{equation} \mathds P_f\left(U_n(H_m) \leq -C  \left( \|f\|_{\infty}+1\right) D R\left(n,\log \left\{\frac{3 \beta \log n}{\gamma}\right\}\right) \right) \leq \gamma/3. \label{test-ustat} \end{equation}

From Eq.\eqref{sup-2-inf} and Eq.\eqref{test-ustat} we get that
\begin{align}
\mathds P_f(T_{\alpha} \leq 0)& \leq  \frac{\gamma}{3}+ \inf_{m\in \mathcal M} \mathds P_f\Bigg((P_n-P)(2 \Pi_{S_m}(f )-2f)+(P_n-P)(2f-2f_0)+\|f-f_0\|^2_2 \notag\\
& \qquad  \qquad \leq \|f- \Pi_{S_m}(f ) \|^2_2+t_m(u_{\alpha})+C  \left( \|f\|_{\infty}+1\right) D R\left(n,\log \left\{\frac{3 \beta \log n}{\gamma}\right\}\right)\Bigg). \label{eq:5.6}
\end{align}

{\bf Control of~$(P_n-P)(2 \Pi_{S_m}(f )-2f)$.}

It is easy to check that there exists some constant~$C'>0$ such that for all~$l$ in~$\{1,2\}$,~$D$ in~$\mathds D_l$, 
\[ \left|2 \Pi_{S_{(l,D)}}(f )(X_i)-2f(X_i)\right|\leq C'\|f\|_{\infty}.\] 

Indeed,

\begin{itemize}
\item when~$l=1$, for any~$k \in \mathds Z$, \[\langle \sqrt D \mathds 1_{[k/D,(k+1)/D[} ,f \rangle = \int \sqrt D \mathds 1_{[k/D,(k+1)/D[}(x)f(x)dx \leq D^{-1/2}\|f\|_{\infty}.\]
Hence, \begin{align*}
&\sup_x |\Pi_{S_{(1,D)}}(f )(x)| \leq \sup_x \sum_{k \in \mathds Z } \left|\langle \sqrt D \mathds 1_{[k/D,(k+1)/D[} ,f \rangle \right|  \sqrt D \mathds 1_{[k/D,(k+1)/D[}(x)\\
\leq\quad&  D^{-1/2}\|f\|_{\infty} \sup_x \sum_{k \in \mathds Z}\sqrt D \mathds 1_{[k/D,(k+1)/D[}(x)= \|f\|_{\infty}  .\end{align*}
\item when~$l=2$,~$D=2^J$ for some~$J \in \mathds N$ and we have for any~$k \in \mathbb Z$, \[\langle \phi_{J,k},f \rangle = \int 2^{J/2} \phi(2^Jx-k)f(x)dx \leq \|f\|_{\infty}\int 2^{J/2} |\phi(2^Jx)|dx \leq 2^{-J/2}\|f\|_{\infty} \|\phi\|_1.\]
Hence, \begin{align*}
&\sup_x |\Pi_{S_{(2,D)}}(f )(x)| \leq \sup_x \sum_{k \in \mathds Z} \left|\langle\phi_{J,k} ,f \rangle \right|\times |  \phi_{J,k}(x)|\\
\leq\quad& 2^{-J/2} \|f\|_{\infty} \|\phi\|_1 \sup_x \sum_{k \in \mathds Z}|2^{J/2}\phi(2^Jx-k)|\leq c\|f\|_{\infty}  \|\phi\|_1 ,\end{align*}
where~$c>0$ is a constant depending only on~$\phi$ since~$\phi$ is bounded and compactly supported. Stated otherwise, there is only a finite number of integers~$k \in \mathds Z$ (which is independent of~$x$ and~$J$) such that for any~$x \in \mathds R$ and any~$J\in \mathds Z$,~$ 2^Jx-k$ falls into the support of~$\phi$.  \end{itemize}

\medskip

Moreover, it is proved in \cite[Page 269]{devore1993constructive}, that one can take~$C'$ such that for all~$D$ in~$\mathds D_3$, \[ |2 \Pi_{S_{(3,D)}}(f )(X_i)-2f(X_i)|\leq C' \|f\|_{\infty} \log(D+1).\] Since \[\mathds E_{X \sim \pi}\left(2\Pi_{S_m}(f )(X)-2f(X)\right)^2\leq 4 \|f\|_{\infty} \| \Pi_{S_m}(f )-f\|_2^2,\]
we can deduce using Proposition~\ref{bernstein-markov} (see Section~\ref{apdx-bernstein}) that for all~$m=(l, D)\in \mathcal M$,
\begin{align*}\mathds P_f\Bigg(  (P_n-P)(2 \Pi_{S_m}(f )-2f) &< - \frac{2 C'\log(3C_{\chi}/\gamma) qA_1 \|f\|_{\infty} \log(D+1)}{n}\\
&- 2\sqrt{\frac{2\log(3C_{\chi}/\gamma)q A_2 \|f\|_{\infty}}{n}} \|\Pi_{S_m}(f )-f\|_2 \Bigg) \leq  \frac{\gamma}{3}.\end{align*}

Considering some~$\epsilon \in ]0,2[$, we use the inequality~$\forall a,b \in \mathds R,$~$2ab\leq 4a^2/\epsilon+\epsilon b^2/4$ and we obtain that for any~$m=(l, D)\in \mathcal M$,
\begin{align}\mathds P_f\Bigg(  (P_n-P)(2 \Pi_{S_m}(f )-2f) + \frac{\epsilon}{4}\|\Pi_{S_m}(f )-f\|_2^2 &< - \frac{2 C'\log(3C_{\chi}/\gamma)qA_1 \|f\|_{\infty} \log(D+1)}{n} \notag\\
&- \frac{8\log(3C_{\chi}/\gamma)q A_2 \|f\|_{\infty}}{\epsilon n}  \Bigg) \leq  \frac{\gamma}{3}.\label{eq:5.7}\end{align}

The control of~$(P_n-P)(2f-2f_0)$ is computed in the same way and we get
\begin{align}\mathds P_f\Bigg(  (P_n-P)(2 f-2f_0) + \frac{\epsilon}{4}\|f-f_0\|_2^2 &< - \frac{4 \log(3C_{\chi}/\gamma)qA_1 (\|f\|_{\infty}+\|f_0\|_{\infty}) }{n} \notag\\
&- \frac{8\log(3C_{\chi}/\gamma)q A_2 \|f\|_{\infty}}{\epsilon n}  \Bigg) \leq  \frac{\gamma}{3}. \label{eq:5.8}\end{align}

Finally, we deduce from Eq.\eqref{eq:5.6}, Eq.\eqref{eq:5.7} and Eq.\eqref{eq:5.8} that if there exists some~$m=(l, D)$ in~$\mathcal M$ such that
\begin{align*}
\left(1-\frac{\epsilon}{4}\right) \|f-f_0\|^2_2 &> \left(1+\frac{\epsilon}{4}\right)\|f-\Pi_{S_m}(f)\|^2_2  +\frac{8\log(3C_{\chi}/\gamma) qA_2 \|f\|_{\infty}}{\epsilon n} \\
&+\frac{4 \log(3C_{\chi}/\gamma)qA_1 (\|f\|_{\infty}+\|f_0\|_{\infty}) }{n}\\
& + \frac{8\log(3C_{\chi}/\gamma)q A_2 \|f\|_{\infty}}{\epsilon n} +\frac{2 C'\log(3C_{\chi}/\gamma)qA_1 \|f\|_{\infty} \log(D+1)}{n}\\
& + t_m(u_{\alpha})+C  \left( \|f\|_{\infty}+1\right) D R\left(n,\log \left\{\frac{3 \beta \log n}{\gamma}\right\}\right),
\end{align*}

i.e. such that
\begin{align*}
\left(1-\frac{\epsilon}{4}\right) \|f-f_0\|^2_2 &> \left(1+\frac{\epsilon}{4}\right)\|f-\Pi_{S_m}(f)\|^2_2  +\frac{16\log(3C_{\chi}/\gamma) qA_2 \|f\|_{\infty}}{\epsilon n}\\
&  +4\left(\|f\|_{\infty}( C' \log(D+1) +1) + \|f_0\|_{\infty}\right)  \frac{\log(3C_{\chi}/\gamma)qA_1  }{n}\\
& + t_m(u_{\alpha})+C \left( \|f\|_{\infty}+1\right) D R\left(n,\log \left\{\frac{3 \beta \log n}{\gamma}\right\}\right),
\end{align*}
then
\[\mathds P_f\left( T_{\alpha}\leq 0  \right) \leq \gamma.\]
To conclude the proof of Theorem~\ref{hypo-test-appli3}, it suffices to notice that for any~$\epsilon \in ]0,2[$, choosing~$\eta>0$ such that~$1+\eta = \frac{1+\frac{\epsilon}{4}}{1-\frac{\epsilon}{4}}$ leads to~$\epsilon=\frac{4 \eta}{2+\eta}$. One can immediately check that the condition~$\epsilon \in ]0,2[$ is equivalent to~$\eta \in ]0,2[$. Noticing further that~$\frac{1}{\epsilon} = \frac{2+\eta}{4\eta} < \frac{2+2}{4\eta}=\frac{1}{\eta}$, we deduce that for any~$\eta \in ]0,2[$, if
\begin{align*}
\|f-f_0\|^2_2 &> (1+\eta) \Bigg\{\|f-\Pi_{S_m}(f)\|^2_2  +\frac{16\log(3C_{\chi}/\gamma) q A_2 \|f\|_{\infty}}{\eta n}\\
&  +4\left(\|f\|_{\infty}( C' \log(D+1) +1) + \|f_0\|_{\infty}\right)  \frac{\log(3C_{\chi}/\gamma)q A_1  }{n}\\
& + t_m(u_{\alpha})+C \left( \|f\|_{\infty}+1\right) D R\left(n,\log \left\{\frac{3 \beta \log n}{\gamma}\right\}\right) \Bigg\},
\end{align*}
then
\[\mathds P_f\left( T_{\alpha}\leq 0  \right) \leq \gamma.\]

\subsection{Proof of Lemma~\ref{lemma:appli3}} \label{proof-lemma:appli3}

Lemma~\ref{lemma:appli3} will follow from Theorem~\ref{mainthm2} if we can show that the function~$H_m$ is bounded. Let us denote~$m=(l,D)$ for some~$l\in \{1,2,3\}$ and~$D \in \mathds D_l$. Let us first remark that the Bessel's inequality states that 
\begin{equation} \label{bessel}\sum_{k \in \mathcal L_m} |\langle p_k,f\rangle |^2 \leq \|f\|_2^2 = \int f(x)f(x)dx \leq \|f\|_{\infty},\end{equation}
since~$\int f(x)dx = 1$ and~$f(x) \geq 0,  \quad \forall x$.

$\bullet~$ If~$l=1$, then we notice that for any~$k \in \mathds Z$,
\begin{align*}| \langle \sqrt D \mathds {1}_{]k/D,(k+1)/D[},f\rangle |&= \left|\int \sqrt D \mathds {1}_{]k/D,(k+1)/D[}(x)f(x)dx \right| \\
&\leq \|f\|_{\infty} \sqrt D \int \mathds {1}_{]k/D,(k+1)/D[}(x)dx \\
&\leq D^{-1/2} \|f\|_{\infty}.
\end{align*}
Then for any~$x,y \in \mathds R$ it holds
\begin{align*}
|H_m(x,y)| &\leq \sum_{k \in \mathcal L_m} |p_k(x)p_k(y)| + \sum_{k \in \mathcal L_m} |p_k(x)\langle p_k,f\rangle |+\sum_{k \in \mathcal L_m} |p_k(y)\langle p_k,f\rangle |+\sum_{k \in \mathcal L_m} |\langle p_k,f\rangle |^2\\
&\leq  \sum_{k\in \mathds Z} D \mathds{1}_{]k/D,(k+1)/D[}(x)\mathds{1}_{]k/D,(k+1)/D[}(y) \\ &\qquad +2\sup_z \sum_{k \in \mathds Z} \sqrt D|\mathds{1}_{]k/D,(k+1)/D[}(z)| \times | \langle \sqrt D \mathds{1}_{]k/D,(k+1)/D[},f\rangle |+\sum_{k \in \mathcal L_m} |\langle p_k,f\rangle |^2\\
&\leq D + 2\|f\|_{\infty} + \|f\|_{\infty},
\end{align*}
where in the last inequality we used Eq.\eqref{bessel}.

$\bullet~$ If~$l=2$ then~$D=2^J$ for some~$J \in \mathds N$ and we have for any~$k \in \mathbb Z$, \[\langle \phi_{J,k},f \rangle = \int 2^{J/2} \phi(2^Jx-k)f(x)dx \leq \|f\|_{\infty}\int 2^{J/2} |\phi(2^Jx)|dx \leq 2^{-J/2}\|f\|_{\infty} \|\phi\|_1.\]
We get that for any~$ x,y \in \mathds R$,
\begin{align*}
|H_m(x,y)| &\leq \sum_{k \in \mathcal L_m} |p_k(x)p_k(y)| + \sum_{k \in \mathcal L_m} |p_k(x)\langle p_k,f\rangle |+\sum_{k \in \mathcal L_m} |p_k(y)\langle p_k,f\rangle |+\sum_{k \in \mathcal L_m} |\langle p_k,f\rangle |^2\\
&\leq  \sum_{k\in \mathds Z} 2^J \phi(2^Jx-k)\phi(2^Jy-k) +2\sup_z \sum_{k \in \mathds Z} 2^{-J/2} \|f\|_{\infty} \|\phi\|_1 2^{J/2} |\phi(2^{J/2}z-k)| +\sum_{k \in \mathcal L_m} |\langle p_k,f\rangle |^2\\
&\leq c 2^J+c'\|\phi\|_{1}\|f\|_{\infty} + \|f\|_{\infty}\\
&=  c D+c'\|\phi\|_{1}\|f\|_{\infty} + \|f\|_{\infty},
\end{align*}
for some constants~$c,c'>0$. In the last inequality we used Eq.\eqref{bessel} and the fact~$\phi$ is bounded and compactly supported. Indeed, this implies that there is only a finite number of integers~$k \in \mathds Z$ (which is independent of~$x$ and~$J$) such that for any~$x \in \mathds R$ and any~$J\in \mathds Z$,~$ 2^Jx-k$ falls into the support of~$\phi$. 

$\bullet~$ If~$l=3$ then we easily get for any~$ x,y \in [0,1]$,
\begin{align*}
|H_m(x,y)| &\leq \sum_{k \in \mathcal L_m} |p_k(x)p_k(y)| + \sum_{k \in \mathcal L_m} |p_k(x)\langle p_k,f\rangle |+\sum_{k \in \mathcal L_m} |p_k(y)\langle p_k,f\rangle |+\sum_{k \in \mathcal L_m} |\langle p_k,f\rangle |^2\\
&\leq  2D + 4 D \|f\|_{\infty} +\|f\|_{\infty}.
\end{align*}

We deduce that in any case,~$H_m$ is bounded~$c(1+\|f\|_{\infty})D$ for some constant~$c>0$ (depending only on~$\phi$) which concludes the proof of Lemma~\ref{lemma:appli3}.

\subsection{Proof of Corollary~\ref{cor:appli3}} \label{proof-cor:appli3}

\underline{Step 1}: We start by providing an upper bound on~$t_m(u_{\alpha})$ with Lemma~\ref{lemma2:appli3}.

\begin{lemma} \label{lemma2:appli3}There exists a constant~$C(\alpha)>0$ such that for any~$m=(l,D) \in \mathcal M$ it holds,
\[t_m(u_{\alpha}) \leq W_m(\alpha),\]
where
 \[W_m(\alpha)= C(\alpha)\left( \|f_0\|_{\infty}+1\right) \left[  D R\left(n,\log \log n\right)+  \frac{\log \log n }{n} \right].\]
\end{lemma}

\underline{Proof of Lemma~\ref{lemma2:appli3}.}

Let us recall that~$t_m(u)$ denotes the~$(1-u)$ quantile of the distribution of~$\widehat T_m$ under the null hypothesis. One can easily see that~$|\mathcal M|\leq  3 (1+\log_2 n)$. So, setting~$\alpha_n=\alpha/(3 (1+\log_2 n))$, \begin{align*}\mathds P_{f_0}(\sup_{m \in \mathcal M} (\widehat T_m-t_m(\alpha_n))>0)&\leq \sum_{m \in \mathcal M} \mathds P_{f_0}(\widehat T_m-t_m(\alpha_n)>0)\\
&\leq \sum_{m \in \mathcal M} \alpha /(3 (1+\log_2 n))\\
&\leq \alpha.\end{align*}
By definition of~$u_{\alpha}$, this implies that~$\alpha_n\leq u_{\alpha}$ and for all~$m \in \mathcal M$, \[t_m(u_{\alpha})\leq t_m(\alpha_n).\]
Hence it suffices to upper bound~$t_m(\alpha_n)$. Let~$m=(l, D) \in \mathcal M$. We use the same notation as in the proof of Theorem~\ref{hypo-test-appli3} to obtain that 
\[\widehat T_m=U_n(H_m)+(P_n-P)(2\Pi_{S_m}(f ))-2P_n(f_0)+\|f_0\|_2^2+\|\Pi_{S_m}(f )\|_2^2.\]
Under the null hypothesis, this reads as 
\begin{align*}\widehat T_m&=U_n(H_m)+(P_n-P)(2\Pi_{S_m}(f_0 )-2f_0)-\|f_0\|_2^2+\|\Pi_{S_m}(f_0 )\|_2^2\\
&=U_n(H_m)+(P_n-P)(2\Pi_{S_m}(f_0 )-2f_0)-\|f_0-\Pi_{S_m}(f_0 )\|_2^2.
\end{align*}

We control~$U_n(H_m)$ and~$(P_n-P)(2\Pi_{S_m}(f_0 )-2f_0)$ exactly like in the proof of Theorem~\ref{hypo-test-appli3}.

From Lemma~\ref{lemma:appli3}, there exist~$C, \beta>0$ such that for any~$m=(l,D) \in \mathcal M$, it holds 
\begin{equation} \mathds P_{f_0}\left(U_n(H_m) \leq C \left( \|f_0\|_{\infty}+1\right) D R\left(n,\log \left\{\frac{2 \beta \log n}{\alpha_n}\right\}\right) \right) \leq \alpha_n/2.  \end{equation}
Moreover, since
 \[ |2 \Pi_{S_{(l,D)}}(f_0 )(X_i)-2f_0(X_i)|\leq C' \|f_0\|_{\infty} \log(D+1),\] and \[\mathds E_{X \sim \pi}\left(2\Pi_{S_m}(f_0 )(X)-2f_0(X)\right)^2\leq 4 \|f_0\|_{\infty} \| \Pi_{S_m}(f_0 )-f_0\|_2^2,\]
 we get using Proposition~\ref{bernstein-markov} (see Section~\ref{apdx-bernstein}) that for all~$m=(l, D)\in \mathcal M$,
\begin{align*}\mathds P_{f_0}\Bigg(  (P_n-P)(2 \Pi_{S_m}(f_0 )-2f_0) &>  \frac{2 C'\log(2 C_{\chi}/\alpha_n)qA_1 \|f_0\|_{\infty} \log(D+1)}{n}\\
& +2\sqrt{\frac{2\log(2 C_{\chi}/\alpha_n) qA_2 \|f_0\|_{\infty}}{n}} \|\Pi_{S_m}(f_0 )-f_0\|_2 \Bigg) \leq  \frac{\alpha_n}{2}.\end{align*}
Using the inequality~$\forall a,b \in \mathds R, \; 2ab\leq a^2+b^2$, and the fact that for~$n\geq 16, \;  \log(D+1)\leq \log(n^2+1)$, we obtain that there exists~$C''>0$ such that
\begin{align*}\mathds P_{f_0}\Bigg(  (P_n-P)(2 \Pi_{S_m}(f_0 )-2f_0)-\|\Pi_{S_m}(f_0 )-f_0\|_2^2 &>  \frac{ C''\|f_0\|_{\infty} \log(2C_{\chi}/\alpha_n) \log(n) }{n}\Bigg) \leq  \frac{\alpha_n}{2}.\end{align*}

We deduce that it holds
\begin{align*}\mathds P_{f_0}\Bigg(  \widehat T_m> C \left( \|f_0\|_{\infty}+1\right) D R\left(n,\log \left\{\frac{2 \beta \log n}{\alpha_n}\right\}\right)+  \frac{ C''\|f_0\|_{\infty} \log(2C_{\chi}/\alpha_n) \log(n) }{n}\Bigg) \leq  \alpha_n.\end{align*}

Noticing that there exists some constant~$c(\alpha)>0$ such that 
\[ \log \left\{\frac{2 \beta \log n}{\alpha_n}\right\} \vee \log(2C_{\chi}/\alpha_n) \leq c(\alpha) \log \log n,\]
we deduce by definition of~$t_m(\alpha_n)$ that for some $c(\alpha)>0$, \[t_m(\alpha_n) \leq c(\alpha) C\left( \|f_0\|_{\infty}+1\right) D R\left(n,\log \log n\right)+  c(\alpha)\frac{ C''\|f_0\|_{\infty} \log \log n }{n}.\]
\hfill $\blacksquare$

\underline{Step 2}: Proof of Corollary~\ref{cor:appli3}.

Let us fix~$\gamma \in]0,1[$ and~$l \in \{1,2,3\}$.  From  Theorem~\ref{hypo-test-appli3} and Lemma~\ref{lemma2:appli3}, we deduce that if~$f$ satisfies
\[\|f-f_0\|_2^2>(1+\epsilon) \inf_{D\in \mathcal D_l}{\|f-\Pi_{S_{(l,D)}}(f )\|^2_2+W_{(l,D)}(\alpha)+V(l,D)(\gamma)},\]
then 
\[\mathds P_f(T_{\alpha}\leq 0)\leq \gamma.\]
It is thus a matter of giving an upper bound for 
\[\inf_{D\in \mathcal D_l} \left\{\|f-\Pi_{S_{(l,D)}}(f )\|^2_2+W_{(l,D)}(\alpha)+V_{(l,D)}(\gamma)\right\},\]
when~$f$ belongs to some specified classes of functions. Recall that
\[\mathcal B^{(l)}_s(P, M)=\{f\in L_ 2(\mathds R)\; |\; \forall D\in \mathcal D_l,\|f-\Pi_{S_{(l,D)}}(f )\|_2^2\leq P^2D^{-2s},\|f\|_{\infty} \leq M\}.\]
We now assume that~$f$ belongs to~$\mathcal B^{(l)}_s(P, M)$. Since~$\|f-\Pi_{S_{(l,D)}}(f )\|^2_2\leq P^2D^{-2s}$, we only need an upper bound for

\begin{align*}
 \inf_{D \in \mathcal D_l} \Bigg\{ & P^2 D^{-2s} + C(\alpha)\left( \|f_0\|_{\infty}+1\right) \left[  D R\left(n,\log \log n\right)+  \frac{\log \log n }{n} \right] +C_1\|f\|_{\infty}\frac{\log(3C_{\chi}/\gamma)  }{\epsilon n} \\
 &+C_2\left(\|f\|_{\infty}  \log(D+1)+ \|f_0\|_{\infty}\right)  \frac{\log(3C_{\chi}/\gamma)  }{n} +C_3 \left( \|f\|_{\infty}+1\right)D R\left(n,\log \left\{\frac{3 \beta \log n}{\gamma}\right\}\right) \Bigg\}.
\end{align*}

Using that~$f$ belongs to~$\mathcal B^{(l)}_s(P, M)$ and the fact that 
\[R\left(n,\log \log n\right) \vee R\left(n,\log \left\{\frac{3 \beta \log n}{\gamma}\right\}\right) \lesssim \log (n) \frac{\log \log n}{n}  ,\]
where~$\lesssim$ states that the inequality holds up to some multiplicative constant independent of~$n$,~$D$ and~$P$, we deduce that we want to upper bound 
\begin{align*}
 \inf_{D \in \mathcal D_l} \Bigg\{ & P^2 D^{-2s} +  D \log (n) \frac{\log \log n}{n} +  \frac{\log \log n }{n} 
+ \frac{\log(D+1) }{n}  \Bigg\}.
\end{align*}
Since~$\log(D+1) \leq D$ for all~$D \in \mathcal D_l$, we only need to focus on
\begin{align*}
 \inf_{D \in \mathcal D_l} \Bigg\{ & P^2 D^{-2s} +  D \log (n) \frac{\log \log n}{n}  \Bigg\}.
\end{align*}

$P^2 D^{-2s}<D \log (n) \frac{\log \log n}{n}$ if and only if~$D> \left( \frac{P^4 n^2}{\log ^2(n) (\log \log n)^2} \right)^{\frac{1}{4s+2}}$. Hence we define~$D_*$ by
\[\log_2(D_*) :=  \floor{\log_2 \left( \left( \frac{P^4 n^2}{\log ^2(n) (\log \log n)^2} \right)^{\frac{1}{4s+2}}   \right)} +1.\]
We consider three cases.
\begin{itemize}
\item If~$D_*<1$, then~$P^2 D^{-2s}<D \log (n) \frac{\log \log n}{n}$ for any~$D \in  \mathcal D_l$ and by choosing~$D_0=1$ to upper bound the infimum we get 
\[\inf_{D\in \mathcal D_l} \left\{\|f-\Pi_{S_{(l,D)}}(f )\|^2_2+W_{(l,D)}(\alpha)+V_{(l,D)}(\gamma)\right\}\leq \log (n) \frac{\log \log n}{n}.\]
\item If~$D_*> 2^{\floor{\log_2( n/ (\log(n) \log \log n)^2)}}$, then~$P^2 D^{-2s}>D \log (n) \frac{\log \log n}{n}$ for any~$D \in  \mathcal D_l$ and by choosing~$D_0=2^{\log_2(\floor{ n/ (\log(n) \log \log n)^2)})}$ to upper bound the infimum we get 
\[\inf_{D\in \mathcal D_l} \left\{\|f-\Pi_{S_{(l,D)}}(f )\|^2_2+W_{(l,D)}(\alpha)+V_{(l,D)}(\gamma)\right\}\lesssim 2P^2D_0^{-2s}\leq  2^{2s+1} P^2 \left( \frac{(\log(n)\log \log n)^2}{n} \right)^{2s}.\]
\item Otherwise~$D_*~$ belongs to~$\mathcal D_l$ and we upper bound the infimum by choosing~$D_0=D_*$ and we get 
\[\inf_{D\in \mathcal D_l} \left\{\|f-\Pi_{S_{(l,D)}}(f )\|^2_2+W_{(l,D)}(\alpha)+V_{(l,D)}(\gamma)\right\}\lesssim 4 P^{\frac{2}{2s+1}} \left(\frac{\log (n) \log \log n}{n}\right)^{\frac{2s}{2s+1}}.\]
\end{itemize}

The proof of Corollary~\ref{cor:appli3} ends with simple computations that we provide below for the sake of completeness. Since
\begin{align*}
  \log (n) \frac{\log \log n}{n} & \leq P^{\frac{2}{2s+1}} \left(\frac{\log (n) \log \log n}{n}\right)^{\frac{2s}{2s+1}}\\
  \Leftrightarrow    \left( \log (n) \frac{\log \log n}{n}\right)^{1/2} & \leq  P.
\end{align*}
and since 

\begin{tabular}{crl}
 &~$\displaystyle P^2 \left( \frac{(\log(n)\log \log n)^2}{n} \right)^{2s}$ &~$\displaystyle \leq P^{\frac{2}{2s+1}} \left(\frac{\log (n) \log \log n}{n}\right)^{\frac{2s}{2s+1}}$\\
 ~$\Leftrightarrow$&  ~$ \displaystyle P \left( \frac{(\log(n)\log \log n)^2}{n} \right)^{s}$ &~$\displaystyle \leq P^{\frac{1}{2s+1}} \left(\frac{\log (n) \log \log n}{n}\right)^{\frac{s}{2s+1}}$\\
 ~$\Leftrightarrow$ & ~$\displaystyle  P^{2s} \left( \frac{(\log(n)\log \log n)^2}{n} \right)^{s(2s+1)}$ &~$\displaystyle \leq  \left(\frac{\log (n) \log \log n}{n}\right)^{s}$\\
 ~$\Leftrightarrow$ & ~$ \displaystyle P \left( \frac{(\log(n)\log \log n)^2}{n} \right)^{s+1/2}$ &~$\displaystyle \leq  \left(\frac{\log (n) \log \log n}{n}\right)^{1/2}$\\
~$\Leftrightarrow$ &  ~$P$&~$\displaystyle \leq  \frac{n^s}{(\log(n)\log \log n)^{2s+1/2}}$,
\end{tabular}

we deduce that if~$P$ is chosen such that 
\begin{equation} \label{final-1} \left( \log (n) \frac{\log \log n}{n}\right)^{1/2} \leq P \leq  \frac{n^s}{(\log(n)\log \log n)^{2s+1/2}},\end{equation}
then the uniform separation rate of the test~$\mathds 1_{T_{\alpha}>0}$ over~$\mathcal B^{(l)}_s(P, M)$ satisfies
\begin{equation} \label{final-2}\rho\left(\mathds{1}_{T_{\alpha}>0},\mathcal B^{(l)}_s(P, M), \gamma \right) \leq C' P^{\frac{1}{2s+1}}  \left( \frac{\log(n) \log \log n}{n}  \right)^{\frac{s}{2s+1}}.\end{equation}

{\bf Remark} This final statement can allow the reader to understand our choice for the size of the model~$|\mathcal M|$ that we considered. Indeed, we chose for any~$l \in \{1,2,3\}$,~$\mathcal D_l=\{2^J,0\leq J\leq \log_2\left( n/ (\log(n) \log \log n)^2  \right)\}$ in order to ensure that for values of~$P$ saturing the right inequality in \eqref{final-1} (i.e. for~$P \approx \frac{n^s}{(\log(n)\log \log n)^{2s+1/2}}$), the upper-bound in Eq.\eqref{final-2} still tends to zero as~$n$ goes to~$+\infty$ for any possible values of the smoothness parameter~$s$.

\section{Concentration Lemmas for Markov chains}

\label{sec:add-proofs}
\subsection{Hoeffding inequality for uniformly ergodic Markov chains}

\label{proof-prop:bernstein-orlicz}

Proposition~\ref{prop:bernstein-orlicz} is an Hoeffding bound for uniformly ergodic Markov chains.

\begin{proposition}
\label{prop:bernstein-orlicz}
Let $(X_i)_{i\geq 1}$ be a Markov chain on $E$ uniformly ergodic (namely satisfying Assumption~\ref{assumption1}) with stationary distribution $\pi$ and let us consider some function $f:E \to \mathds R$ such that $\mathds E_{X \sim \pi}[f(X)]=0$. Then it holds for any $t\geq0$

\begin{align*}
\mathds P\left( \left|\sum_{i=1}^{n}f(X_i)\right|\geq t \right)
&\leq 16\exp \left( -\frac{1}{K(m,\tau)}  \frac{t^2}{ n\|f\|_{\infty}^2} \right),
\end{align*}
where $K(m,\tau)= 2 Km^2\tau^2$ for some universal constant $K>0.$ We refer to Assumption~\ref{assumption1} and the following remark (or to \cite[Section 2]{duchemin}) for the definitions of the constants $m$ and $\tau$.
\end{proposition}

\underline{Proof of Proposition~\ref{prop:bernstein-orlicz}.} Let us first recall that under Assumption~\ref{assumption1}, the  $1$-Orlicz norm of the regeneration times of the split chain are bounded by some finite constant $\tau>0$ (see the remark after Assumption~\ref{assumption1}). In this proof, we will use the notations introduced in \cite[Section 2.3]{duchemin}. Since the chain $(\widetilde{X}_n)_n$ is distributed as $(X_i)_{i\geq1}$, we will identify $(\widetilde X_i)_{i\geq1}$ and $(X_i)_{i\geq1}$ in the proof.

Let us consider~$N = \sup \{i \in \mathbb N \;:\; mS_{i+1}+m-1\leq n\}$. Then,
\begin{align}
\big| \sum_{i=1}^{n}f(X_i)\big|&=\big| \sum_{l=0}^{N}Z_l+\sum_{i=m(S_{N}+1)}^{n}f(X_i)\big| \leq \big|\sum_{l=0}^{\floor{N/2}}Z_{2l}\big|+\big|\sum_{l=0}^{\floor{(N-1)/2}}Z_{2l+1}\big|+\big|\sum_{i=m(S_{N}+1)}^{n}f(X_i)\big|. \label{apdx:bern-Tau}
\end{align}

We have~$|\sum_{i=m(S_{N}+1)}^{n}f(X_i)|\leq AmT_{N+1}$. So using the definition of the Orlicz norm and the fact that the random variables~$(T_i)_{i\geq 2}$ are i.i.d., it holds for any~$t\geq 0$,
\begin{align*}
\mathds P\big(\big|\sum_{i=m(S_{N}+1)}^{n}f(X_i)\big|\geq t\big)&\leq \mathds P(T_{N+1}\geq \frac{t}{Am})\leq \mathds P(\max(T_1,T_2)\geq \frac{t}{Am})
\leq 4 \exp\big(-\frac{t}{Am\tau}\big).
\end{align*}

In order to control the first two terms in~\eqref{apdx:bern-Tau}, we need to describe the tail behaviour of the random variable~$N$ with Lemma~\ref{adam-lemma5}.
\begin{lemma} \label{adam-lemma5} \citep[cf.][Lemma 5]{adamczak:ustats} \\
We denote~$R=\floor{3n/(\mathds E T_2)}.$ If~$\|T_1\|_{\psi_1},\|T_2\|_{\psi_1}\leq \tau$, then $\mathds P (N> R) \leq 2 \exp\left(-  \frac{n \mathds E T_2}{8\tau^2}  \right).$
\end{lemma}
The random variable~$Z_{2l}$ is~$\sigma(X_{m(S_{2l}+1)},\dots,X_{m(S_{2l+1}+1)-1})$-measurable. Hence the random variables~$(Z_{2l})_l$ are independent \citep[see][Section 2.3]{duchemin}. Moreover, one has that for any~$l$,~$\mathds E[Z_{2l} ]=0$. This is due to~\cite[Eq.(17.23) Theorem~17.3.1]{tweedie} together with the assumption that~$\mathds E_{X\sim\pi}[f(X)]=0.$ Let us finally notice for any~$l\geq0$,~$|Z_{2l}| \leq Am  T_{2l+1}$, so~$\|Z_{2l}\|_{\psi_1}\leq Am\max(\|T_1\|_{\psi_1},\|T_2\|_{\psi_1})\leq Am \tau$. One can similarly get that~$(Z_{2l+1})_l$ are independent with~$\mathds E[Z_{2l+1} ]=0$ and~$\|Z_{2l+1}\|_{\psi_1}\leq Am \tau$ for all~$l \in \mathds N$. Using these facts we have for any~$t\geq 0$,
\begin{align*}
&\mathds P\big(\big|\sum_{l=0}^{\floor{N/2}}Z_{2l}\big|+\big|\sum_{l=0}^{\floor{(N-1)/2}}Z_{2l+1}\big|\geq t\big)\\
&\leq \mathds P\big(\big|\sum_{l=0}^{\floor{N/2}}Z_{2l}\big|+\big|\sum_{l=0}^{\floor{(N-1)/2}}Z_{2l+1}\big|\geq t\;,N\leq R\big) +2 \exp\big(- \frac{n \mathds E T_2}{8\tau^2}  \big)\\
&\leq \mathds P\big( \max_{0\leq s \leq \floor{R/2}}\big|\sum_{l=0}^{s}Z_{2l}\big|\geq t/2\big) + \mathds P\big(\max_{0\leq s \leq \floor{(R-1)/2}}\big|\sum_{l=0}^{s}Z_{2l+1}\big|\geq t/2\big) +2 \exp\big(-  \frac{n \mathds E T_2}{8\tau^2}  \big)\\
&\leq 3\mathds P\big(\big|\sum_{l=0}^{\floor{R/2}}Z_{2l}\big|\geq t/6\big)+3\mathds P \big(\big|\sum_{l=0}^{\floor{(R-1)/2}}Z_{2l+1}\big|\geq t/6\big) +2 \exp\big(- \frac{n \mathds E T_2}{8\tau^2}  \big) \quad \text{(Using Lemma~\ref{adam-lemma4})}\\
&\leq 12\exp \big( -\frac{1}{8} \min \big( \frac{t^2}{36RA^2m^2\tau^2},\frac{t}{6Am\tau}\big) \big)+2 \exp\big(-  \frac{n \mathds E T_2}{8\tau^2}\big),
\end{align*}
where we used Lemma~\ref{bernstein-psi1} in the last inequality.
\begin{lemma} (Bernstein's~$\psi_1$ inequality,~\cite[Lemma 2.2.11]{VW00} and the subsequent remark).\label{bernstein-psi1} \\If~$Y_1,\dots,Y_n$ are independent random variables such that~$\mathds EY_i= 0$ and~$\|Y_i\|_{\psi_1}\leq \tau$, then for every~$t >0,$
\[ \mathds P \left( \left|\sum_{i=1}^n Y_i\right|>t\right)\leq 2\exp \left( -\frac{1}{K} \min \left( \frac{t^2}{n\tau^2},\frac{t}{\tau}\right) \right) ,\]
for some universal constant~$K>0$ ($K=8$ fits).
\end{lemma}

\begin{lemma} \label{adam-lemma4} \citep[cf.][Proposition 1.1.1]{KW92} 
If~$X_1, X_2, \dots$ are independent Banach space valued random variables (not necessarily identically distributed), and if~$S_k=\sum_{i=1}^kX_i$, then \[\mathds P \left(  \underset{1\leq j \leq k}{\max}\|S_j\|>t \right)\leq 3 \underset{1\leq j \leq k}{\max}\mathds P \left(  \|S_j\|>t/3\right).\]
\end{lemma}

Gathering the previous results, we obtain that for any~$t\geq 0$
\begin{align*}
\mathds P\big( \big|\sum_{i=1}^{n}f(X_i)\big|\geq t \big)
&\leq 12\exp \big( -\frac{1}{8} \min \big( \frac{t^2\big(\mathds E T_2\big)}{36\times 12 \times nA^2m^2\tau^2},\frac{t}{12Am\tau}\big) \big)\\
&\quad +2 \exp\big(-  \frac{n \mathds E T_2}{8\tau^2}\big) + 4 \exp\big(-\frac{t}{2Am\tau}\big).
\end{align*}
Since the left hand side of the previous inequality is zero for~$t\geq nA,$ and since~$m\geq 1$, we obtain Proposition~\ref{prop:bernstein-orlicz}.

\subsection{Bernstein's inequality for non-stationary Markov chains}
\label{apdx-bernstein}

Proposition~\ref{bernstein-markov} is an extension of the Bernstein type concentration inequality from~\cite{jiang2018bernsteins} to non-stationary Markov chains. A proof can be found in the Appendix of \cite{duchemin}. Let us highlight that Proposition~\ref{bernstein-markov} is only used in the proofs of the main results from Section~\ref{sec:adaptivegof}. One could have used other concentration results such as the one from~\cite{paulin15} (by using jointly Theorem~3.4 and Proposition~3.10) which would give strictly analogous results.

\begin{proposition}  \label{bernstein-markov}
Suppose that the sequence~$(X_i)_{i\geq1}$ is a Markov chain satisfying Assumptions~\ref{assumption1} and~\ref{assumption0} with stationary distribution~$\pi$ and with an absolute spectral gap~$1-\lambda>0$. Let us consider some~$n \in \mathds N \backslash \{0\}$ and bounded real valued functions~$(f_i)_{1\leq i \leq n}$ such that for any~$i \in \{1, \dots,n\}$,~$\int f_i(x) d\pi(x)=0$ and~$\|f_i\|_{\infty}\leq c$ for some~$c>0$. Let~$\sigma^2 = \sum_{i=1}^n \int f_i^2(x) d\pi(x) /n$.
Then for any~$\epsilon \geq 0$ it holds
\[\mathds P \left( \sum_{i=1}^n f_i(X_i) \geq \epsilon \right) \leq  \left\| \frac{d\chi}{d\pi} \right\|_{\pi,p} \exp\left(- \frac{\epsilon^2/(2q)}{A_2\sigma^2+A_1c\epsilon}\right) ,\]
where~$A_2 := \frac{1+\lambda}{1-\lambda}$ and~$A_1:=\frac13 \mathds1_{\lambda=0}+ \frac{5}{1-\lambda}\mathds 1_{\lambda >0}$.~$q$ is the constant introduced in Assumption~\ref{assumption0}. Stated otherwise, for any~$u>0$ it holds
\[\mathds P \left( \frac{1}{n} \sum_{i=1}^n f_i(X_i)> \frac{2quA_1 c}{n} + \sqrt{\frac{2quA_2\sigma^2}{n}} \right) \leq  \left\| \frac{d\chi}{d\pi} \right\|_{\pi,p} e^{-u}.\]

\end{proposition}

\newpage
\vskip 0.2in
\bibliography{sample}

\end{document}